\def\DateTime{25/February/1998, 3:30PM}
\def\Version{Version 3.0}

\documentclass[leqno,12pt]{amsart}
\usepackage{amssymb}
\usepackage{amscd}
\usepackage{latexsym}
\usepackage{verbatim}

\theoremstyle{plain}
\newtheorem{Theorem}{Theorem}[section]
\newtheorem{Proposition}[Theorem]{Proposition}
\newtheorem{Lemma}[Theorem]{Lemma}
\newtheorem{Corollary}[Theorem]{Corollary}
\newtheorem{Claim}{Claim}[Theorem]

\theoremstyle{definition}

\newtheorem{Remark}[Theorem]{Remark}
\newtheorem{Example}[Theorem]{Example}

\renewcommand{\theTheorem}{\arabic{section}.\arabic{subsection}.\arabic{Theorem}}
\renewcommand{\theClaim}{\arabic{section}.\arabic{subsection}.\arabic{Theorem}.\arabic{Claim}}
\renewcommand{\theequation}{\arabic{section}.\arabic{subsection}.\arabic{Theorem}.\arabic{Claim}}

\def\rom{\textup}
\newcommand{\NN}{{\mathbb{N}}}
\newcommand{\ZZ}{{\mathbb{Z}}}
\newcommand{\QQ}{{\mathbb{Q}}}
\newcommand{\RR}{{\mathbb{R}}}
\newcommand{\CC}{{\mathbb{C}}}
\newcommand{\PP}{{\mathbb{P}}}

\newcommand{\OO}{{\mathcal{O}}}
\newcommand{\rank}{\operatorname{rk}}
\newcommand{\codim}{\operatorname{codim}}

\newcommand{\diag}{\operatorname{diag}}
\newcommand{\GL}{\operatorname{GL}}
\newcommand{\Aut}{\operatorname{Aut}}
\newcommand{\HEnd}{\operatorname{End}}

\newcommand{\dis}{\operatorname{dis}}
\newcommand{\adis}{\widehat{\operatorname{dis}}}
\newcommand{\adeg}{\widehat{\operatorname{deg}}}

\newcommand{\Ker}{\operatorname{Ker}}

\newcommand{\Image}{\operatorname{Image}}
\newcommand{\Sym}{\operatorname{Sym}}
\newcommand{\Supp}{\operatorname{Supp}}
\newcommand{\End}{\operatorname{\mathcal{E}\textsl{nd}}}
\newcommand{\Hom}{\operatorname{\mathcal{H}\textsl{om}}}
\newcommand{\Spec}{\operatorname{Spec}}
\newcommand{\Gal}{\operatorname{Gal}}
\newcommand{\Rat}{\operatorname{Rat}}

\newcommand{\Cycle}{Z}
\newcommand{\aCycle}{\widehat{Z}}
\newcommand{\aACycle}{\widehat{Z}_A}
\newcommand{\aBCycle}{\widehat{Z}_{L^1}}
\newcommand{\aDCycle}{\widehat{Z}_D}

\newcommand{\aChow}{\widehat{\operatorname{CH}}}
\newcommand{\aAChow}{\widehat{\operatorname{CH}}_A}
\newcommand{\aBChow}{\widehat{\operatorname{CH}}_{L^1}}
\newcommand{\aDChow}{\widehat{\operatorname{CH}}_D}

\newcommand{\Div}{\operatorname{Div}}
\newcommand{\aDiv}{\widehat{\operatorname{Div}}}
\newcommand{\aBDiv}{\widehat{\operatorname{Div}}_{L^1}}
\newcommand{\aDDiv}{\widehat{\operatorname{Div}}_D}

\newcommand{\aPic}{\widehat{\operatorname{Pic}}}
\newcommand{\aBPic}{\widehat{\operatorname{Pic}}_{L^1}}
\newcommand{\aDPic}{\widehat{\operatorname{Pic}}_D}

\newcommand{\chernch}{\operatorname{ch}}
\newcommand{\achernch}{\widehat{\operatorname{ch}}}

\newcommand{\cherncl}{{c}}
\newcommand{\acherncl}{\widehat{{c}}}

\newcommand{\todd}{\operatorname{td}}
\newcommand{\atodd}{\widehat{\operatorname{td}}}

\newcommand{\alev}{\operatorname{a.e.}}
\newcommand{\VERT}{|\!|\!|}
\newcommand{\LocInt}{L^1_{\operatorname{loc}}}
\newcommand{\zero}{\operatorname{div}}
\newcommand{\Lint}{L\!\!\!\int}
\newcommand{\Proof}{{\sl Proof.}\quad}
\newcommand{\QED}{{\unskip\nobreak\hfil\penalty50\quad\null\nobreak\hfil
{$\Box$}\parfillskip0pt\finalhyphendemerits0\par\medskip}}
\newcommand{\rest}[2]{\left.{#1}\right\vert_{{#2}}}

\newcommand{\Wedge}{\operatornamewithlimits{\wedge}}
\newcommand{\Norm}{\operatorname{Norm}}
\newcommand{\Perm}{{\mathfrak{S}}}
\newcommand{\sign}{\operatorname{sign}}
\newcommand{\Tr}{\operatorname{Tr}}
\newcommand{\GM}{\operatorname{g.\!m.}}

\newcommand{\Ch}{\operatorname{Ch}}

\newcommand{\bfGL}{\operatorname{GL}}
\newcommand{\bfSL}{\operatorname{SL}}
\newcommand{\Unitary}{\operatorname{U}}
\newcommand{\Proj}{\operatorname{Proj}}

\begin{document}

\title[Inequalities for semistable families of arithmetic varieties]%
{Inequalities for semistable families \\
of arithmetic varieties}
\author{Shu Kawaguchi}
\author{Atsushi Moriwaki}
\address{Department of Mathematics, Faculty of Science,
Kyoto University, Kyoto, 606-01, Japan}
\email[Shu Kawaguchi]{kawaguch@kusm.kyoto-u.ac.jp}
\email[Atsushi Moriwaki]{moriwaki@kusm.kyoto-u.ac.jp}
\thanks{The first author is partially supported by JSPS Research Fellowships
for Young Scientists.}
\date{\DateTime, (\Version)}
\begin{abstract}
In this paper, we will consider a generalization
of Bogomolov's inequality and Cornalba-Harris-Bost's inequality
to the case of semistable families of arithmetic varieties under the idea
that geometric semistability implies a certain kind of
arithmetic positivity.
The first one is 
an arithmetic analogue of the relative Bogomolov's inequality
in \cite{MoRB}.
We also establish the arithmetic Riemann-Roch formulae
for stable curves over regular arithmetic varieties and
generically finite morphisms of arithmetic varieties.
\end{abstract} 

\maketitle


\section*{Introduction}
\renewcommand{\theTheorem}{\Alph{Theorem}}

In this paper, we will consider a generalization
of Bogomolov's inequality and Cornalba-Harris-Bost's inequality
to the case of semistable families of arithmetic varieties.
An underlying idea of these inequalities as in
\cite{Bo}, \cite{BGS}, \cite{Ga}, \cite{MiBi}, \cite{MoBG},
\cite{MoABG}, \cite{MoBU}, \cite{MorFh}, \cite{SoVan},
and \cite{Zh}
is that geometric semistability implies a certain kind of
arithmetic positivity.
The first one is related to the semistability of vector bundles,
and the second one involves the Chow (or Hilbert) semistability
of cycles.

\medskip
First of all, let us consider Bogomolov's inequality.
Let $X$ and $Y$ be smooth algebraic varieties over
an algebraically closed field of characteristic zero, and
$f : X \to Y$ a semi-stable curve.
Let $E$ be a vector bundle of rank $r$ on $X$,
and $y$ a point of $Y$.
In \cite{MoRB}, the second author proved that if
$f$ is smooth over $y$ and $\rest{E}{X_{\bar{y}}}$
is semistable, then $\dis_{X/Y}(E) = f_* \left(
2r c_2(E) - (r-1)c_1^2(E) \right)$
is weakly positive at $y$.

In the first half of this paper, we would like to consider an
arithmetic analogue of the above result.
Let us fix regular arithmetic varieties $X$ and $Y$,
and a semistable curve $f : X \to Y$.
Since we have a good dictionary for
translation from a geometric case to an arithmetic case,
it looks like routine works.
There are, however, two technical difficulties
to work over the standard dictionary.

The first one is how to define a push-forward
of arithmetic cycles
in our situation. If $f_{\QQ} : X_{\QQ} \to Y_{\QQ}$
is smooth, then,
according to Gillet-Soul\'{e}'s arithmetic intersection theory
\cite{GSArInt}, we can get the push-forward
$f_* : \aChow^{p+1}(X) \to \aChow^p(Y)$. 
We would not like to restrict ourselves to
the case where $f_{\QQ}$ is smooth because
in the geometric case,
the weak positivity of $\dis_{X/Y}(E)$ gives
wonderful applications to analyses of
the boundary of the moduli space of stable curves.
Thus the usual push-forward
for arithmetic cycles is insufficient
for our purpose.
A difficulty in defining the push-forward arises from a fact:
if $f_{\CC} : X(\CC) \to Y(\CC)$ is not smooth,
then $(f_{\CC})_*(\eta)$ is not necessarily $C^{\infty}$
even for a $C^{\infty}$ form $\eta$.
This suggests us that we need to extend the usual arithmetic
Chow groups defined by Gillet-Soul\'{e} \cite{GSArInt}.
For this purpose, we will introduce 
an arithmetic $L^1$-cycle of
codimension $p$, namely, a pair $(Z, g)$ such that
$Z$ is a cycle of codimension $p$,
$g$ is a current of type $(p-1,p-1)$, and
$g$ and $dd^c(g) + \delta_{Z(\CC)}$ are
represented by locally integrable forms.
Thus, dividing by the usual arithmetical rational
equivalence, an arithmetic Chow group,
denoted by $\aBChow^p$, consisting
of arithmetic $L^1$-cycles of codimension $p$ will be defined
(cf. \S\ref{subsec:var:arith:chow}).
In this way, we have the natural push-forward
\[
f_* : \aBChow^{p+1}(X) \to \aBChow^p(Y)
\]
as desired (cf. Proposition~\ref{prop:push:forward:arith:cycle}).

The second difficulty is the existence of
a suitable Riemann-Roch formula in our situation.
As before, if $f_{\QQ} : X_{\QQ} \to Y_{\QQ}$ is smooth,
we have the arithmetic Riemann-Roch theorem due to
Gillet-Soul\'{e} \cite{GSRR}. If we ignore Noether's formula,
then, under the assumption that
$f_{\QQ} : X_{\QQ} \to Y_{\QQ}$ is smooth,
their Riemann-Roch formula can be written in the following form:
\begin{multline*}
\acherncl_1 \left( \det Rf_*(E), h_Q^{\overline{E}} \right) -
\rank (E) \acherncl_1 \left( \det Rf_*(\OO_X), h_Q^{\overline{\OO}_X}
\right) \\ = f_* \left( \frac{1}{2} \left(
\acherncl_1 (\overline{E})^2  -
\acherncl_1 (\overline{E}) \cdot \acherncl_1 (\overline{\omega}_{X/Y}) \right)
- \acherncl_2 (\overline{E})
\right)
\end{multline*}
where $\overline{E} = (E, h)$ is
a Hermitian vector bundle on $X$ and
$\overline{\omega}_{X/Y}$ is the dualizing sheaf of $f : X \to Y$
with a Hermitian metric.
If we consider a general case 
where $f_{\QQ} : X_{\QQ} \to Y_{\QQ}$ 
is not necessarily smooth, 
the right hand side in the
above equation is well defined and
sits in $\aBChow^1(X)_{\QQ}$.
On the other hand, the left hand side is rather complicated.
If we admit singular fibers of $f_{\CC} : X(\CC) \to Y(\CC)$,
then the Quillen metric $h_Q^{\overline{E}}$ is no longer 
$C^{\infty}$. According to \cite{BBQm},
it extends to a generalized metric. Thus, we may define
$\acherncl_1 \left( \det Rf_*(E), h_Q^{\overline{E}} \right)$
(cf. \S\ref{subsec:arith:div:gen:metric}).
In general, this cycle is not an $L^1$-cycle.
However, using Bismut-Bost's formula \cite{BBQm},
we can see that
\[
\acherncl_1 \left( \det Rf_*(E), h_Q^{\overline{E}} \right) -
\rank (E) \acherncl_1 \left( \det Rf_*(\OO_X), h_Q^{\overline{\OO}_X}
\right)
\]
is an element of $\aBChow^1(Y)$. Thus, we have a way to establish
a Riemann-Roch formula in the arithmetic Chow group
$\aBChow^1(Y)_{\QQ}$. Actually,
we will prove the above formula in our situation
(cf. Theorem~\ref{thm:arith:Riemann:Roch:stable:curves}).
The idea of comparing two sides in $\aBChow^1(Y)_{\QQ}$
is the tricky Lemma~\ref{lem:criterion:linear:equiv:B:cycle}.

Let us go back to our problem.
First of all, we need to define an arithmetic analogue of
weak positivity.
Let $\alpha$ be an element of $\aBChow^1(Y)_{\QQ}$, $S$ a subset of $Y(\CC)$,
and $y$ a closed point of $Y_{\QQ}$.
We say $\alpha$ is semi-ample at $y$ with respect to $S$ if
there are an arithmetic $L^1$-cycle $(E, f)$ and a positive integer
$n$ such that 
(1) $dd^c(f) + \delta_{E(\CC)}$ is $C^{\infty}$ around each $z \in S$,
(2) $E$ is effective, 
(3) $y \not\in \Supp(E)$,
(4) $f(z) \geq 0$ for all $z \in S$, and
(5) $n \alpha$ coincides with the class of $(E, f)$ in $\aBChow^1(Y)_{\QQ}$.
Moreover, $\alpha$ is said to be weakly positive at $y$ with respect to $S$
if it is the limit of semi-ample cycles at $y$ with respect to $S$
(for details, see \S\ref{subsec:wp:div}).
For example, if $Y = \Spec(O_K)$, $y$ is the generic point, and
$S = Y(\CC)$, then,
$\alpha$ is weakly positive at $y$ with respect to $S$
if and only if $\adeg(\alpha) \geq 0$,
where $K$ is a number field and $O_K$ is the ring of integers in $K$
(cf. Proposition~\ref{prop:wp:for:curve}).

Let $(E, h)$ be a Hermitian vector bundle of rank $r$ on $X$, and
$\adis_{X/Y}(E, h)$ the arithmetic discriminant divisor of 
$(E, h)$ with respect
to $f : X \to Y$, that is, the element of $\aBChow^1(Y)$ given by
$f_* \left( 2r \acherncl_2(E, h) - (r-1)\acherncl_1(E, h)^2 \right)$.
We assume that $f$ is smooth over $y$ and $\rest{E}{X_{\bar{y}}}$ is
poly-stable.
In the case where $\dim X = 2$ and $Y = \Spec(O_K)$,
Miyaoka \cite{MiBi}, Moriwaki \cite{MoBG,MoABG,MoBU}, and Soul\'{e} \cite{SoVan}
proved that $\adeg \left( \adis_{X/Y}(E, h) \right) \geq 0$, consequently,
$\adis_{X/Y}(E, h)$ is weakly positive at $y$ with respect to $Y(\CC)$.
One of the main theorems of this paper is the following generalization.

\begin{Theorem}[cf. Theorem~\ref{thm:relative:Bogomolov:inequality:arithmetic:case}]
\label{thm:A:relative:Bogomolov:inequality}
Under the above assumptions, $\adis_{X/Y}(E, h)$
is weakly positive at $y$ with respect to any subsets $S$ of $Y(\CC)$ 
with the following
properties: \rom{(1)} $S$ is finite, and
\rom{(2)} $f_{\CC}^{-1}(z)$ is smooth and $\rest{E_{\CC}}{f_{\CC}^{-1}(z)}$ is
poly-stable for all $z \in S$.
In particular, if the residue field of $x$ is $K$,
and the canonical morphism $\Spec(K) \to X$ 
induced by $x$ extends to $\tilde{x} : \Spec(O_K) \to X$,
then $\adeg\left(  \tilde{x}^*\left(\adis_{X/Y}(E, h)\right)\right) \geq 0$.
\end{Theorem}

\medskip
Next, let us consider Cornalba-Harris-Bost's inequality. 
Motivated by the work of Cornalba and Harris \cite{CoHa} 
in the geometric case, 
Bost \cite[Theorem~I]{Bo} proved that, roughly speaking, 
if $X(\overline{\QQ}) \subset \PP^{r-1}(\overline{\QQ})$ has 
the $\bfSL_r(\overline{\QQ})$ semi-stable Chow point, 
then the height of $X$ has a certain kind of positivity. 
We call this result Cornalba-Harris-Bost's inequality. 
Zhang \cite{Zh} then gave precision to it and 
also showed the converse of Bost's result. 
Further, Gasbarri \cite{Ga} considered a wide range of actions 
instead of $\bfSL_r(\overline{\QQ})$-action.

In the second half of this paper, 
we would like to consider a relative version of
Cornalba-Harris-Bost's inequality.
First, let us fix a terminology.
Let $V$ be a set, $\phi$ a non-negative function on $V$, and
$S$ a finite subset of $V$.
We define the geometric mean $\GM(\phi; S)$ of
$\phi$ over $S$ to be
\[
\GM(\phi; S) = \left( \prod_{s \in S} \phi(s) \right)^{1/\#(S)}.
\]
Then, the following is our solution.

\begin{Theorem}[cf. Theorem~\ref{thm:semistability:imply:average:semi-ampleness}]
\label{thm:intro:B}
Let $Y$ be a regular projective arithmetic variety, and 
$\overline{E} = (E,h)$ a Hermitian vector bundle of rank $r$. 
Let $\pi : \PP(E) = \Proj (\bigoplus_{n \geq 0} \Sym^n(E^{\lor})) \to Y$ 
be the projection and $\overline{\OO_{E}(1)}$ 
the tautological line bundle with the quotient metric induced from $f^*(h)$. 
Let $X$ be an effective cycle in $\PP (E)$ 
such that $X$ is flat over $Y$ 
with the relative dimension $d$ 
and degree $\delta$ on the generic fiber. 
For each irreducible component $X_i$ of $X_{red}$,
let $\tilde{X}_i \to X_i$ be a proper birational morphism
such that $(\tilde{X}_i)_{\QQ}$ is smooth over $\QQ$.
Let $Y_0$ be the maximal open set of $Y$ such that
the induced morphism $\tilde{X}_i \to Y$ is smooth over $Y_0$
for every $i$.
Let $(B, h_B)$ be a line bundle equipped with
a generalized metric on $Y$ given by the equality:
\[
\acherncl_1(B, h_B) = 
r \pi_* \left( \acherncl_1(\overline{\OO_{E}(1)})^{d+1} \cdot (X,g_X) \right) 
+ \delta (d+1) \acherncl_1(\overline{E}).
\]
\rom{(}Here we postpone the definition of $g_X$, i.e.,
a suitable compactification of $X$ in the arithmetic sense.\rom{)}
Then, $h_B$ is $C^{\infty}$ over $Y_0$. Moreover,
there are a positive integer
$e=e(r,d,\delta)$, a positive integer $l=l(r,d,\delta)$,
a positive constant $C=C(r,d,\delta)$, and
sections $s_1, \ldots, s_l \in H^0(Y, B^{\otimes e})$
with the following properties.
\begin{enumerate}
\renewcommand{\labelenumi}{(\roman{enumi})}
\item
$e$, $l$, and $C$ depend only on $r$, $d$, and $\delta$.

\item
For a closed point $y$ of $Y_{\QQ}$, if $X_y$ is Chow semistable,
then $s_i(y) \not= 0$ for some $i$.

\item
For all $i$ and all closed points $y$ of $(Y_0)_{\QQ}$,
\[
\GM\left( \left( h_B^{\otimes e} \right)(s_i, s_i);\ 
          O_{\Gal(\overline{\QQ}/\QQ)}(y)\right)
\leq C,
\]
where $O_{\Gal(\overline{\QQ}/\QQ)}(y)$ is the orbit of $y$
by the Galois action in $Y_0(\overline{\QQ})$.
\end{enumerate}
\end{Theorem}

Compared with the geometric analogue 
(cf. Remark~\ref{rem:geom:analog:Cornalba-Harris-Bost}),
the difficult part of this theorem is
the estimate of the geometric mean of the norm over
the Galois orbits of closed points.
We will do this by reducing to the absolute case. 
For this purpose, we have to associate $X$ with a `nice' Green current $g_X$. 
How do we do? One way is to fix a K\"{a}hler metric 
$\mu \in A^{1,1}(\PP(E)_{\RR})$ and 
to attach a $\mu$-normalized Green current for $X$, 
namely, a Green current $g$ such that 
$dd^c g + \delta_X = H(\delta_Y)$ and $H(g_Y) = 0$,  
where $H : D^{p,p}(\PP(E)_{\RR}) \to H^{p,p}(\PP(E)_{\RR})$ is 
the harmonic projection (cf. \cite[2.3.2]{BGS}). 
This construction however is not suitable for our purpose 
because it does not behave well 
when restricted on fibers.  

Thus we are led to define an $\Omega$-normalized Green form  
which is given, roughly speaking, by attaching a Green form fiberwisely 
(Here $\Omega = \cherncl_1(\overline{\OO_{E}(1)})$). 
Precisely, an $\Omega$-normalized Green form $g_X$ for $X$ is 
characterized by the following three conditions;  
(i) $g_X$ is an $L^1$-form on $\PP(E)$,  
(ii) $dd^c([g_X]) + \delta_X 
= \sum_{i=0}^{d} \left[ \pi^*(\gamma_i) \wedge \Omega^i \right]$.
where $\gamma_i$ is a $d$-closed $L^1$-form of type 
$(d-i,d-i)$ 
on $Y$ ($i=0, \ldots, d$). 
(iii) $\pi_*(g_X \wedge \Omega^{r - d}) = 0$ 
(cf. Proposition~\ref{prop:normalized:Green:form}). 
Then we can show that it has a desired property 
when restricted on fibers 
(cf. Remark~\ref{rem:norm:Green:general:fiber}). 

Suppose now $X$ is regular. 
Let $i : X \to \PP(E)$ be the inclusion map 
and $f : X \to Y$ the restriction of $\pi$. 
If we set $\overline{L} = i^*(\overline{\OO_{E}(1)})$, 
then 
$\pi_* ( \acherncl_1(\overline{\OO_{E}(1)})^{d+1} \cdot (X,g_X) )
= f_*(\acherncl_1(\overline{L})^{d+1})$ 
(cf. Proposition~\ref{prop:when:Bost:divisor:smooth}). 
Since $f_*(\acherncl_1(\overline{L})^{d+1})$ is in general  
only an element of $\aChow^1_{L^1}(Y)$, 
the above equality explains why we need to consider $(X,g_X)$ in 
the enlarged arithmetic Chow group $\aChow^{r-d-1}_{L^1}(\PP(E))$. 
Moreover, a similar equality when $X$ is not necessarily regular 
shows that 
$\pi_* ( \acherncl_1(\overline{\OO_{E}(1)}) \cdot (X,g_X) )$ 
is independent of the choice of an $\Omega$-normalized Green 
form $g_X$ for $X$ 
(cf. Proposition~\ref{prop:when:Bost:divisor:smooth}). 

Suppose now $Y=\Spec(O_K)$, $y$ is the generic point, and
$X_y$ is Chow semistable, 
where $K$ is a number field. 
In this case, 
there exists a generic resolution of $X$ 
smooth over $y$. Then Theorem~\ref{thm:intro:B} tells us that 
\[
r \adeg ( \acherncl_1(\overline{L})^{d+1} ) 
+ \delta (d+1) \adeg (\overline{E}) 
+ [K:\QQ] \alpha(r,d,\delta) \geq 0 
\]
for some constant $\alpha(r,d,\delta)$ depending only on
$r$, $d$ and $\delta$,
which is nothing but Theorem~I of Bost \cite{Bo}. 


We can also think a wide range of actions like \cite{Ga}. 
Namely, let $\rho : \GL_r \to \GL_R$ be a morphism of group schemes 
such that there is an integer $k$ with $\rho(t I_r) = t^k I_R$ for any $t$, 
and that $\rho$ commutes with the transposed morphism. 
For a Hermitian vector bundle $\overline{E}$, 
we then get the associated Hermitian vector bundle $\overline{E}^{\rho}$ 
(cf. \S\ref{subsec:associated:herm:vb}). 
If $X$ is a flat cycle on $\PP(E^{\rho})$ 
and $y$ is a closed point of $Y_{\QQ}$, 
then $\bfSL_r(\overline{\QQ})$ acts a Chow form ${\Phi_X}_y$. 
The stability of ${\Phi_X}_y$ under this action yields a similar inequality 
(cf. Theorem~\ref{thm:semistability:imply:average:semi-ampleness}). 

Finally, in \S\ref{section:Bogomolov:to:Bost} 
we make a comparison between 
the relative Bogomolov's inequality 
(Theorem~\ref{thm:relative:Bogomolov:inequality:arithmetic:case}) 
and the relative Cornalba-Harris-Bost's inequality 
(Theorem~\ref{thm:semistability:imply:average:semi-ampleness}).

\renewcommand{\theTheorem}{\arabic{section}.\arabic{subsection}.\arabic{Theorem}}
\section{Locally integrable forms and their push-forward}

\subsection{Locally integrable forms}
\setcounter{Theorem}{0}
Let $M$ be an $n$-dimensional orientable differential manifold.
We assume that $M$ has a countable basis of open sets.
Let $\omega$ be a $C^{\infty}$ volume element of $M$, and
$C_c^0(M)$ the set of all complex valued continuous functions on $M$
with compact supports.
Then, there is a unique Radon measure $\mu_{\omega}$ 
defined on the topological $\sigma$-algebra of $M$ such that
\[
\Lint_M f d\mu_{\omega} = \int_M f \omega
\]
for all $f \in C_c^0(M)$, where
${\displaystyle \Lint_M f d\mu_{\omega}}$ is the Lebesgue integral arising
from the measure $\mu_{\omega}$.

Let $f$ be a complex valued
function on $M$. We say $f$ is {\em locally integrable},
denoted by $f \in \LocInt(M)$, if $f$ is measurable
and, for any compact sets $K$,
\[
\Lint_K |f| d\mu_{\omega} < \infty.
\]
Let $\omega'$ be another $C^{\infty}$ volume form on $M$.
Then, there is a positive $C^{\infty}$ function $a$ on $M$
with $\omega' = a \omega$.
Thus,
\[
\Lint_K |f| d\mu_{\omega'} = \Lint_K |f| a d\mu_{\omega},
\]
which shows us that local integrability does not depend on the choice
of the volume form $\omega$. Moreover, it is easy to see that,
for a measurable complex valued function $f$ on $M$,
the following are equivalent.
\begin{enumerate}
\renewcommand{\labelenumi}{(\alph{enumi})}
\item
$f$ is locally integrable.

\item
For each point $x \in M$, there is an open neighborhood
$U$ of $x$ such that the closure of $U$ is compact and
${\displaystyle \Lint_U |f| d\mu_{\omega} < \infty}$.
\end{enumerate}

Let $\Omega_M^p$ be a $C^{\infty}$ vector bundle consisting
of $C^{\infty}$ complex valued $p$-forms.
Let $\pi_p : \Omega_M^p \to M$ be the canonical map.
We denote $C^{\infty}(M, \Omega_M^p)$
(resp. $C^{\infty}_c(M, \Omega_{M}^p)$) by $A^p(M)$ (resp. $A_c^p(M)$).
Let $\alpha$ be a section of $\pi_p : \Omega_M^p \to M$.
We say $\alpha$ is 
{\em locally integrable}, or simply an {\em $L^1$-form}
if, at any point of $M$,
all coefficients of $\alpha$ in terms of
local coordinates are locally integrable functions.
The set of all locally integrable $p$-forms is denoted by
$\LocInt(M, \Omega^p_M)$.
For an maximal form $\alpha$ on $M$,
there is a unique function $g$ on $M$ with
$\alpha = g \omega$. 
We denote this function $g$ by $c_{\omega}(\alpha)$.

Let us define the Lebesgue integral of locally integrable
$n$-forms with compact support.
Let $\alpha$ be an element of $\LocInt(M, \Omega^n_M)$ such that
the support of $\alpha$ is compact.
Then $c_{\omega}(\alpha) \in \LocInt(M)$ and
$\operatorname{supp}(c_{\omega}(\alpha))$ is compact. Thus, 
${\displaystyle \Lint_{M} c_{\omega}(\alpha) d\mu_{\omega}}$ exists.
Let $\omega'$ be another $C^{\infty}$ volume element of $M$.
Then, there is a positive $C^{\infty}$ function $a$ on $M$
with $\omega' = a \omega$. Here $a c_{\omega'}(\alpha) = c_{\omega}(\alpha)$.
Thus,
\[
\Lint_{M} c_{\omega'}(\alpha) d\mu_{\omega'} =
\Lint_{M} c_{\omega'}(\alpha) a d\mu_{\omega} =
\Lint_{M} c_{\omega}(\alpha) a d\mu_{\omega}.
\]
Hence, ${\displaystyle \Lint_{M} c_{\omega}(\alpha) d\mu_{\omega}}$
does not depend on the choice of the volume form $\omega$.
Thus, the Lebesgue integral of $\alpha$ is defined by
\[
\Lint_M \alpha = \Lint_{M} c_{\omega}(\alpha) d\mu_{\omega}.
\]

Moreover, we denote by $D^p(M)$ the space of
currents of type $p$ on $M$.
Then, there is the natural homomorphism
\[
[\ ] : \LocInt(M, \Omega^p_M) \to D^p(M)
\]
given by ${\displaystyle [\alpha](\phi) = \Lint_M \alpha \wedge \phi}$
for $\phi \in A_c^{n-p}(M)$.
It is well known that the kernel of $[ \ ]$ is 
$\{ \alpha \in \LocInt(M, \Omega^p_M) \mid
\alpha = 0 \ (\alev) \}$.
A topology on $D^p(M)$ is defined in the following way.
For an sequence $\{ T_n \}_{n=1}^{\infty}$ in $D^p(M)$,
$T_n \to T$ as $n \to \infty$ if and only if
$T_n(\phi) \to T(\phi)$ as $n \to \infty$ 
for each $\phi \in A_c^{n-p}(M)$.
For an element $T \in D^n(M)$,
by abuse of notation,
we denote by $c_{\omega}(T)$
a unique distribution $g$ on $M$ given by
$T = g \omega$. 

\begin{Proposition}
\label{prop:criterion:loc:int}
Let $T$ be a current of type $p$ on $M$.
Then, the following are equivalent.
\begin{enumerate}
\renewcommand{\labelenumi}{(\arabic{enumi})}
\item
$T$ is represented by a $L^1$-form.

\item
For any $\phi \in A^{n-p}(M)$,
$c_{\omega}(T \wedge \phi)$ is represented by
a locally integrable function.
\end{enumerate}
\end{Proposition}

\Proof
(1) $\Longrightarrow$ (2):
Let $\phi \in A^{n-p}(M)$.
Then, by our assumption, for any point $x \in M$,
there are an open neighborhood $U$ of $x$,
$C^{\infty}$ functions $a_1, \ldots, a_r$ on $U$, and
locally integrable functions $b_1, \ldots, b_r$ on $U$ such that
\[
\rest{c_{\omega}(T \wedge \phi)}{U} = \sum_{i=1}^r [a_i b_i].
\]
Thus, if $K$ is a compact set in $U$, then
\[
\Lint_K \left|\sum_{i=1}^r a_i b_i \right| d\mu_{\omega} \leq
\Lint_K \sum_{i=1}^r |a_i| |b_i| d\mu_{\omega}  \leq
\max_i \sup_{x \in K} \{ |a_i(x)| \} \sum_{i=1}^r \Lint_K |b_i|
d\mu_{\omega} < \infty.
\]
Thus, we get (2).

\medskip
(2) $\Longrightarrow$ (1):
Before starting the proof, we would like to claim
the following fact.
Let $\{ U_{\alpha} \}_{\alpha \in A}$ be an open covering of $M$
such that $A$ is at most a countable set.
Let $\lambda_{\alpha}$ be a locally integrable form $U_{\alpha}$
with $\lambda_{\alpha} = \lambda_{\beta} \ (\alev)$ on 
$U_{\alpha} \cap U_{\beta}$ for
all $\alpha, \beta \in A$. Then, there is a locally integrable
form $\lambda$ on $M$ such that $\lambda = \lambda_{\alpha} \ (\alev)$
on $U_{\alpha}$ for all $\alpha \in A$.
Indeed, let us fix a map $a : M \to A$ with $x \in U_{a(x)}$
and define a form $\lambda$ by $\lambda(x) = \lambda_{a(x)}(x)$.
Then, $\lambda$ is our desired form because
for each $\alpha \in A$,
\[
\{ x \in U_{\alpha} \mid \lambda(x) \not= \lambda_{\alpha}(x) \}
\subseteq
\bigcup_{\beta \in A \setminus \{ \alpha \}}
\{ x \in U_{\alpha} \cap U_{\beta} \mid
\lambda_{\beta}(x) \not= \lambda_{\alpha}(x) \}
\]
and the right hand side has measure zero.

\medskip
Let $U$ be an open neighborhood of a point $x \in M$ and
$(x_1, \ldots, x_n)$ a local coordinate of $U$ such that
$dx_1 \wedge \cdots \wedge dx_n$ coincides with the orientation
by $\omega$.
Then, there is a positive $C^{\infty}$ function $a$ on $U$
with $\omega = a dx_1 \wedge \cdots \wedge dx_n$ over $U$.
We set
\[
T = \sum_{i_1 < \cdots < i_p }
T_{i_1 \cdots i_p} dx_{i_1} \wedge \cdots \wedge dx_{i_p}
\]
for some distributions $T_{i_1 \cdots i_p}$.
We need to show that $T_{i_1 \cdots i_p}$
is represented by a locally integrable function. 
Since $M$ has a countable basis of open sets, by the above claim, 
it is sufficient to check that
$T_{i_1 \cdots i_p}$ is represented by
an integral function on every compact set $K$ in $U$.
Let $f$ be a non-negative $C^{\infty}$ function on $M$ such that
$f = 1$ on $K$ and $\operatorname{supp}(f) \subset U$.
Choose $i_{p+1}, \ldots, i_{n}$ such that
$\{ i_1, \ldots, i_{n} \} = \{ 1, \ldots, n \}$.
Here we set $\phi = f a dx_{i_{p+1}} \wedge \cdots \wedge dx_{i_n}$.
Then, $\phi \in A^{n-p}(M)$ and
\[
T \wedge \phi =
\epsilon T_{i_1 \cdots i_p} f a dx_{1} \wedge \cdots \wedge dx_{n}
= \epsilon T_{i_1 \cdots i_p} f \omega,
\]
where $\epsilon = 1$ or $-1$ depending on
the orientation of $\{ x_{i_1}, \ldots, x_{i_n} \}$.
By our assumption, there is a locally integrable function
$h$ on $M$ with $c_{\omega}(T \wedge \phi) = [h]$.
Thus, $[\epsilon h] = T_{i_1 \cdots i_p} f$.
Therefore, $T_{i_1 \cdots i_p}$ is represented by
$\epsilon h$ on $K$ because $f = 1$ on $K$.
Thus, we get (2).
\QED

\subsection{Push-forward of $L^1$-forms as current}
\label{subsec:push:forward:L1:current}
\setcounter{Theorem}{0}
First of all, we recall the push-forward of currents.
Let $f : M \to N$ be a proper morphism
of orientable manifolds with the relative dimension $d = \dim M - \dim N$.
Then,
\[
f_* : D^p(M) \to D^{p-d}(N)
\]
is defined by $(f_*(T))(\phi) = T(f^*(\phi))$ for
$\phi \in A_c^{\dim N - p + d}(N)$.
It is easy to see that $f_*$ is a continuous homomorphism.
Let us begin with the following lemma.

\begin{Lemma}
\label{lem:push:forward:product}
Let $F$ be an orientable compact differential manifold and
$Y$ an orientable differential manifold. Let $\omega_F$
\rom{(}resp. $\omega_Y$\rom{)} be
a $C^{\infty}$ volume element of $F$
\rom{(}resp. $Y$\rom{)}.
Let $p : F \times Y \to Y$
be the projection to the second factor.
Then, we have the following.
\begin{enumerate}
\renewcommand{\labelenumi}{(\arabic{enumi})}
\item
If $g$ is a continuous function on $F \times Y$,
then $\int_F g \omega_F$ is a continuous function on $Y$.

\item
If $\alpha$ is a continuous maximal form on $F \times Y$,
then $p_*([\alpha])$ is represented by
a unique continuous from.
This continuous form is denoted by $\int_p \alpha$.

\item
For a continuous function $g$ on $F \times Y$,
\[
\left| c_{\omega_Y}\left( \int_{p} g  \omega_F \wedge \omega_Y \right) 
\right|
\leq c_{\omega_Y} \left( \int_{p} |g| \omega_F \wedge \omega_Y \right).
\]
\end{enumerate}
\end{Lemma}

\Proof
(1) This is standard.

\medskip
(2) Since $\omega_F \wedge \omega_Y$ is a volume form on
$F \times Y$, there is a continuous function $g$ on
$F \times Y$ with $\alpha = g \omega_F \wedge \omega_Y$.
Thus, it is sufficient to show that
\[
p_*([\alpha]) = \left[ \left(\int_F g \omega_F \right) \omega_Y \right].
\]
Indeed, by Fubini's theorem,
for $\phi \in A_c^0(Y)$,
\[
p_*([\alpha])(\phi) = \int_{F \times Y} \phi \alpha 
= \int_Y \left( \int_F g \omega_F \right) \phi \omega_Y 
= \left[ \left(\int_F g \omega_F \right) \omega_Y \right](\phi).
\]

\medskip
(3) This is obvious because
\[
\left| \int_F g \omega_F \right| \leq
\int_F |g| \omega_F.
\]
\QED

\begin{Corollary}
\label{cor:cont:ineq:integral:fiber}
Let $f : X \to Y$ be a proper, surjective and smooth morphism
of connected complex manifolds. Let $\omega_X$ and
$\omega_Y$ be volume elements of $X$ and $Y$ respectively. Then,
\begin{enumerate}
\renewcommand{\labelenumi}{(\arabic{enumi})}
\item
For a continuous maximal form $\alpha$ on $X$,
$f_*([\alpha])$ is represented by a unique continuous form.
We denote this continuous form by $\int_{f} \alpha$.

\item
For any continuous functions $g$ on $X$,
\[
\left| c_{\omega_Y}\left( \int_{f} g \omega_X \right) \right|
\leq c_{\omega_Y} \left( \int_{f} |g| \omega_X \right).
\]
\end{enumerate}
\end{Corollary}

\Proof
(1)
This is a local question on $Y$.
Thus, we may assume that there are a compact complex manifold
$F$ and a differomorphism $h : X \to F \times Y$
such that the following diagram is commutative:
\[
\begin{CD}
X @>{\sim}>{h}> F \times Y \\
@V{f}VV @VV{p}V \\
Y @= Y,
\end{CD}
\]
where $p : F \times Y \to Y$ is the natural projection.
Hence, (1) is a consequence of (2) of Lemma~\ref{lem:push:forward:product}.

\medskip
(2)
First, we claim that if the above inequality holds for some special
volume elements $\omega_X$ and $\omega_Y$, then
the same inequality holds for any volume elements.
Let $\omega'_{X}$ and $\omega'_{Y}$ be another volume elements
of $X$ and $Y$ respectively. We set $\omega'_X = a \omega_X$ and
$\omega'_Y = b \omega_Y$. Then, $a$ and $b$ are positive $C^{\infty}$
functions.
Let $g$ be any continuous function on $X$. Then, by our assumption,
\[
\left| c_{\omega_Y}\left( \int_{f} g \omega'_X \right) \right| =
\left| c_{\omega_Y}\left( \int_{f} g a \omega_X \right) \right| \leq
c_{\omega_Y} \left( \int_{f} |g| a \omega_X \right) =
c_{\omega_Y} \left( \int_{f} |g| \omega'_X \right).
\]
On the other hand, for any maximal forms $\alpha$ on $Y$,
\[
c_{\omega_Y}\left( \alpha \right) = 
b c_{\omega'_Y}\left( \alpha \right).
\]
Thus, we get our claim.

Hence, as in the proof of (1),
using the differomorphism $h$ and
(3) of Lemma~\ref{lem:push:forward:product},
we can see (2).
\QED

\begin{Remark}
In the situation of Corollary~\ref{cor:cont:ineq:integral:fiber}, 
if $\alpha$ is a $C^{\infty}$-form on $X$, 
then $f_*([\alpha])$ is represented by a unique $C^{\infty}$-form. 
\end{Remark}

\begin{Proposition}
\label{prop:loc:int:integral:fiber}
Let $f : X \to Y$ be a proper and surjective morphism
of connected complex manifolds.
Let $U$ be a non-empty Zariski open set of $Y$ such that
$f$ is smooth over $U$.
Let $\alpha$ be a compactly supported continuous maximal form on $X$.
If we set
\[
\lambda =
\begin{cases}
{\displaystyle \int_{f^{-1}(U) \to U} \alpha} 
& \text{on $U$}, \\
{} & {} \\
0 & \text{on $Y \setminus U$,}
\end{cases}
\]
then $\lambda$ is integrable.
Moreover, $f_*([\alpha]) = [\lambda]$.
\end{Proposition}

\Proof
Let $\omega_X$ and $\omega_Y$ be volume forms of $X$ and $Y$
respectively. Let $h$ be a function on $Y$ with
$\lambda = h \omega_Y$. Then, $h$ is continuous on $U$ by
Corollary~\ref{cor:cont:ineq:integral:fiber}.
Moreover, let $g$ be a continuous function on $X$ with
$\alpha = g \omega_X$.
We need to show that
$h$ is an integrable function.
First note that $\int_X |g| \omega_X < \infty$
because $g$ is a compactly supported continuous function.
Let $\{ U_n \}_{n=1}^{\infty}$ be a sequence of
open sets such that $\overline{U_n} \subset U$,
$\overline{U_n}$ is compact,
$U_1 \subseteq U_2 \subseteq \cdots \subseteq U_n \subseteq \cdots$, and
$\bigcup_{n=1}^{\infty} U_n = U$.
Here we set
\[
h_n(y) = \begin{cases}
|h(y)| & \text{if $y \in U_n$} \\
0    & \text{otherwise}.
\end{cases}
\]
Then, $0 \leq h_1 \leq h_2 \leq \cdots \leq h_n \leq \cdots$ and
${\displaystyle \lim_{n \to \infty} h_n(y) = |h(y)|}$. 
By Corollary~\ref{cor:cont:ineq:integral:fiber},
\[
\left| \rest{h}{U} \right| \leq
c_{\omega_{Y}}\left(  
\int_{f^{-1}(U) \to U} |g| \omega_X  \right).
\]
Thus,
\begin{align*}
\int_{U_n} |h| \omega_Y 
& \leq \int_{U_n} 
c_{\omega_{Y}}\left( 
\int_{f^{-1}(U_n) \to U_n} 
|g| \omega_X \right) \omega_Y
= \int_{U_n} \int_{f^{-1}(U_n) \to U_n}
|g| \omega_X \\
& = \int_{f^{-1}(U_n)} 
|g| \omega_X  \leq \int_X |g| \omega_X.
\end{align*}
Therefore,
\[
\Lint_{Y} h_n d\mu_{\omega_Y} = \int_{U_n} |h| \omega_Y \leq
\int_{X} 
|g| \omega_X < \infty.
\]
Thus, by Fatou's theorem,
\[
\Lint_{Y} |h| d\mu_{\omega_Y} = 
\lim_{n \to \infty} \Lint_{Y} h_n d\mu_{\omega_Y}
\leq 
\Lint_{X} 
|g| \omega_X < \infty.
\]
Hence, $h$ is integral.

Let $\phi$ be any element of $A^{0}_c(Y)$.
Then, since 
${\displaystyle \lim_{n \to \infty} \mu_{\omega_Y}(Y \setminus U_n) = 0}$
and $h \phi$ is integrable,
by the absolute continuity of Lebesgue integral,
\[
\lim_{n \to \infty} \Lint_{Y \setminus U_n} h \phi
d\mu_{\omega_Y} = 0.
\]
Thus,
\begin{align*}
\Lint_Y \lambda \phi &=
\lim_{n \to \infty} \left(
\Lint_{U_n} h \phi d\mu_{\omega_{Y}} +
\Lint_{Y \setminus U_n} h \phi d\mu_{\omega_{Y}}
\right) \\
& = \lim_{n \to \infty} 
\Lint_{U_n} h \phi d\mu_{\omega_{Y}} 
= \lim_{n \to \infty} 
\int_{U_n} h \phi \omega_Y 
= \lim_{n \to \infty} \int_{U_n} \lambda \phi.
\end{align*}
In the same way,
\[
\int_X \alpha f^*(\phi) =
\lim_{n \to \infty} 
\int_{f^{-1}(U_n)} \alpha f^*(\phi).
\]
On the other hand, we have
\[
\int_{U_n} \lambda \phi = 
\int_{f^{-1}(U_n)} \alpha f^*(\phi).
\]
Hence
\begin{align*}
f_*([\alpha])(\phi) & = [\alpha](f^*(\phi)) 
= \int_X \alpha \wedge f^*(\phi) 
= \lim_{n \to \infty} \int_{f^{-1}(U_n)} \alpha f^*(\phi) \\
& = \lim_{n \to \infty} \int_{U_n} \lambda \phi 
= \Lint_{Y} \lambda \phi
= [\lambda](\phi)
\end{align*}
Therefore, $f_*([\alpha]) = [\lambda]$.
\QED

Let $X$ be an equi-dimensional complex manifold, i.e.,
every connected component has the same dimension.
We denote by $A^{p,q}(X)$ the space of $C^{\infty}$ complex valued
$(p,q)$-forms on $X$. Let $A^{p,q}_c(X)$ be the subspace
of compactly supported forms.
Let $D^{p,q}(X)$ be the space of currents on $X$ of type $(p,q)$.
As before, there is a natural homomorphism
\[
[\ ] : \LocInt(\Omega_{X}^{p,q}) \to D^{p,q}(X).
\]
Then, as a corollary of Proposition~\ref{prop:loc:int:integral:fiber},
we have the following main result of this section.

\begin{Proposition}
\label{prop:push:forward:B:pq}
Let $f : X \to Y$ be a proper morphism of equi-dimensional complex manifolds.
We assume that every connected component of $X$ maps surjectively to
a connected component of $Y$.
Let $\alpha$ be an $L^1$-form of type $(p+d,q+d)$ on $X$, where
$d = \dim X - \dim Y$.
Then there is a $\lambda \in \LocInt(\Omega_Y^{p,q})$
with $f_*([\alpha]) = [\lambda]$.
\end{Proposition}

\Proof
Clearly we may assume that $Y$ is connected.
Since $f$ is proper, there are finitely many connected components
of $X$, say, $X_1, \ldots, X_e$.
If we set $\alpha_i = \rest{\alpha}{X_i}$ and $f_i = \rest{f}{X_i}$ for
each $i$, then 
$f_*([\alpha]) = (f_1)_*([\alpha_1]) + \cdots + (f_e)_*([\alpha_e])$.
Thus, we may assume that $X$ is connected.
Further, since $f_*([ \alpha \wedge f^*(\phi) ]) = f_*([\alpha]) \wedge \phi$
for all $\phi \in A^{\dim Y - p, \dim Y - q}(Y)$,
we may assume that $\alpha$ is a maximal form
by Proposition~\ref{prop:criterion:loc:int}.

Let $g$ be a locally integrable function on $X$ with
$\alpha = g \omega_X$. Since the question is local with respect to $Y$,
we may assume that $g$ is integrable. Thus, 
since $C_c^0(Y)$ is dense on $L^1(Y)$ (cf. \cite[Theorem~3.14]{Ru}), 
there is a sequence
$\{ g_n \}_{n=1}^{\infty}$ of compactly supported continuous functions on $X$
such that
\[
\lim_{n \to \infty} \Lint_X |g_n - g| d\mu_{\omega_X} = 0.
\]
By Proposition~\ref{prop:loc:int:integral:fiber}, for each $n$,
there is an integrable function $h_n$ on $Y$
such that $f_*([g_n \omega_X]) = [h_n \omega_Y]$.
Moreover, by (2) of Corollary~\ref{cor:cont:ineq:integral:fiber},
\[
|h_n - h_m| \leq c_{\omega_Y} \left(
\int_{f^{-1}(U) \to U} |g_n - g_m | \omega_X \right)
\]
over $U$. Thus, we can see
\[
\Lint_Y |h_n - h_m| d\mu_{\omega_Y} \leq 
\Lint_X |g_n - g_m| d\mu_{\omega_X}
\]
for all $n, m$.
Hence, $\{ h_n \}_{n=1}^{\infty}$ is a Cauchy sequence
in $L^1(Y)$. Therefore, 
by the completeness of $L^1(Y)$,
there is an integrable function $h$ on $Y$
with $h = \lim_{n \to \infty} h_n$ in $L^1(Y)$.
Then, for any $\phi \in A_c^{0,0}(Y)$,
\[
\lim_{n \to \infty} \Lint_Y h_n \phi \omega_Y =
\Lint_Y h \phi \omega_Y
\quad\text{and}\quad
\lim_{n \to \infty} \Lint_X g_n f^*(\phi) \omega_X =
\Lint_X g f^*(\phi) \omega_X.
\]
Thus,
\begin{align*}
f_*([\alpha])(\phi) & = 
\Lint_X g f^*(\phi)\omega_X =
\lim_{n \to \infty} \Lint_X g_n f^*(\phi) \omega_X  \\
& = \lim_{n \to \infty} \Lint_Y h_n \phi \omega_Y =
\Lint_Y h \phi \omega_Y = [h\omega_Y](\phi).
\end{align*}
Therefore, $f_*([\alpha]) = [h\omega_Y]$.
\QED

\section{Variants of arithmetic Chow groups}

\subsection{Notation for arithmetic varieties}
\label{subsec:notation:arith:variety}
\setcounter{Theorem}{0}
An {\em arithmetic variety $X$} is 
an integral scheme which is flat and quasi-projective
over $\Spec(\ZZ)$, and
has the smooth generic fiber $X_{\QQ}$.

Let us consider the $\CC$-scheme $X \otimes_{\ZZ} \CC$.
We denote the underlying analytic space of $X \otimes_{\ZZ} \CC$
by $X(\CC)$. We may view $X(\CC)$ as the set of all
$\CC$-valued points of $X$.
Let $F_{\infty} : X(\CC) \to X(\CC)$ be the anti-holomorphic
involution given by the complex conjugation.
For an arithmetic variety $X$, 
every $(p,p)$-form $\alpha$ on $X(\CC)$ is
always assumed to be compatible with $F_{\infty}$, i.e.,
$F_{\infty}^*(\alpha) = (-1)^p \alpha$.

Let $E$ be a locally free sheaf on $X$ of finite rank, and
$\pi : \pmb{E} \to X$ the vector bundle associated with $E$, i.e.,
$\pmb{E} = \Spec\left( \bigoplus_{n=0}^{\infty} \Sym^n(E) \right)$.
As before, we have the analytic space $\pmb{E}(\CC)$ and
the anti-holomorphic involution $F_{\infty} : \pmb{E}(\CC) \to \pmb{E}(\CC)$.
Then, $\pi_{\CC} : \pmb{E}(\CC) \to X(\CC)$ is a holomorphic
vector bundle on $X(\CC)$, and the following diagram is commutative:
\[
\begin{CD}
\pmb{E}(\CC) @>{F_{\infty}}>> \pmb{E}(\CC) \\
@V{\pi_{\CC}}VV @VV{\pi_{\CC}}V \\
X(\CC) @>>{F_{\infty}}> X(\CC)
\end{CD}
\]
Here note that
$F_{\infty} : \pmb{E}(\CC) \to \pmb{E}(\CC)$
is anti-complex linear at each fiber.
Let $h$ be a $C^{\infty}$ Hermitian metric of $\pmb{E}(\CC)$.
We can think $h$ as a $C^{\infty}$ function on
$\pmb{E}(\CC) \times_{X(\CC)} \pmb{E}(\CC)$.
For simplicity, we denote by $F_{\infty}^*(h)$
the $C^{\infty}$ function 
$\left( F_{\infty} \times_{X(\CC)} F_{\infty} \right)^*(h)$
on $\pmb{E}(\CC) \times_{X(\CC)} \pmb{E}(\CC)$.
Then, $\overline{F_{\infty}^*(h)}$ is
a $C^{\infty}$ Hermitian metric of $\pmb{E}(\CC)$.
We say {\em $h$ is invariant under $F_{\infty}$}
if $F_{\infty}^*(h) = \overline{h}$.
Moreover, the pair $(E, h)$ is called {\em a Hermitian vector bundle on $X$}
if $h$ is invariant under $F_{\infty}$.
Note that even if $h$ is not invariant under $F_{\infty}$,
$h + \overline{F_{\infty}^*(h)}$ is an invariant metric.

\subsection{Variants of arithmetic cycles}
\setcounter{Theorem}{0}
\label{subsec:var:arith:chow}
Let $X$ be an arithmetic variety.
We would like to define three types of arithmetic cycles,
namely, arithmetic $A$-cycles, arithmetic $L^1$-cycles, and
arithmetic $D$-cycles.
In the following definition,
$g$ is compatible with $F_{\infty}$ as mentioned in
\S\ref{subsec:notation:arith:variety}.

\begin{enumerate}
\renewcommand{\labelenumi}{(\alph{enumi})}
\item
(arithmetic $A$-cycle on $X$ of codimension $p$) : 
a pair $(Z, g)$ such that
$Z$ is a cycle on $X$ of codimension $p$ and
$g$ is represented by a Green form $\phi$ of $Z(\CC)$, namely,
$\phi$ is a $C^{\infty}$ form on $X(\CC) \setminus
\Supp(Z(\CC))$ of logarithmic type along $\Supp(Z(\CC))$
with $dd^c([\phi]) + \delta_{Z(\CC)} \in A^{p,p}(X(\CC))$.

\item
(arithmetic $L^1$-cycle on $X$ of codimension $p$) : 
a pair $(Z, g)$ such that
$Z$ is a cycle on $X$ of codimension $p$ and,
there are $\phi \in \LocInt(\Omega_{X(\CC)}^{p-1,p-1})$ and
$\omega \in \LocInt(\Omega_{X(\CC)}^{p,p})$ with
$g = [\phi]$ and $dd^c(g) + \delta_{Z(\CC)} = [\omega]$.

\item
(arithmetic $D$-cycle on $X$ of codimension $p$) : 
a pair $(Z, g)$ such that
$Z$ is a cycle on $X$ of codimension $p$ and
$g \in D^{p-1,p-1}(X(\CC))$.
\end{enumerate}

The set of all arithmetic $A$-cycles (resp.
$L^1$-cycles, $D$-cycles) of codimension $p$
is denoted by $\aACycle^p(X)$ (resp.
$\aBCycle^p(X)$, $\aDCycle^p(X)$).

Let $\widehat{R}^p(X)$ be the subgroup of $\aCycle^p(X)$ generated
by the following elements:
\begin{enumerate}
\renewcommand{\labelenumi}{(\roman{enumi})}
\item 
$((f), - [\log |f|^2])$, 
where $f$ is a rational function on some
subvariety $Y$ of codimension $p-1$ and $[\log |f|^2]$ 
is the current defined by
\[
[\log |f|^2](\gamma) = 
        \Lint_{Y(\CC)} (\log |f|^2)\gamma.
\]

\item
$(0, \partial(\alpha) + \bar{\partial}(\beta))$,
where $\alpha \in D^{p-2, p-1}(X(\CC))$,
$\beta \in D^{p-1, p-2}(X(\CC))$.
\end{enumerate}
Here we define
\[
\begin{cases}
\aAChow^p(X) = \aACycle^{p}(X)/\widehat{R}^p(X) \cap \aACycle^{p}(X), \\
\aBChow^p(X) = \aBCycle^{p}(X)/\widehat{R}^p(X) \cap \aBCycle^{p}(X), \\
\aDChow^p(X) = \aDCycle^{p}(X)/\widehat{R}^p(X).
\end{cases}
\]

\begin{Proposition}
\label{prop:AChow:equal:Chow}
The natural homomorphism $\aAChow^p(X) \to \aChow^p(X)$
is an isomorphism.
\end{Proposition}

\Proof
Let $(Z, g) \in \aCycle^p(X)$.
By \cite[Theorem~1.3.5]{GSArInt}, there is 
a Green form $g_Z$ of $Z(\CC)$.
Then, $dd^c(g - [g_Z]) \in A^{p,p}(X(\CC))$.
Hence, by \cite[Theorem~1.2.2]{GSArInt},
there are $a \in A^{d,d}(X(\CC))$ and
$v \in \Image(\partial) + \Image(\bar{\partial})$
with $g - [g_Z] = [a] + v$.
Since $g - [g_Z]$ is compatible with $F_{\infty}$,
replacing $a$ and $v$ by $(1/2)(a + (-1)^pF_{\infty}^*(a))$ and
$(1/2)(v + (-1)^p F_{\infty}^*(v))$ respectively,
we may assume that $a$ and $v$ are compatible with $F_{\infty}$.
Here, $g_Z + a$ is a Green form of $Z$. Thus,
$(Z, [g_Z + a]) \in \aACycle^p(X)$.
Moreover, since
$(Z, g) - (Z, [g_Z + a]) \in \widehat{R}^p(X)$,
our proposition follows.
\QED

Let $f : X \to Y$ be a proper morphism of arithmetic varieties
with $d = \dim X - \dim Y$. Then, we have a homomorphism
\[
  f_* : \aDCycle^{p+d}(X) \to \aDCycle^{p}(Y)
\]
defined by $f_*(Z, g) = (f_*(Z), f_*(g))$.
In the same way as in the proof of \cite[Theorem~3.6.1]{GSArInt},
we can see $f_*(\widehat{R}^{p+d}(X)) \subseteq  \widehat{R}^p(Y)$.
Thus, the above homomorphism induces
\[
   f_* : \aDChow^{p+d}(X) \to \aDChow^{p}(Y).
\]
Then we have the following. 

\begin{Proposition}
\label{prop:push:forward:arith:cycle}
If $f$ is surjective, then
$f_* : \aDChow^{p+d}(X) \to \aDChow^{p}(Y)$ gives rise to
\[
f_* : \aBChow^{p+d}(X) \to \aBChow^{p}(Y).
\]
In particular, we have the homomorphism
$f_* : \aChow^{p+d}(X) \to \aBChow^{p}(Y)$.
\end{Proposition}

\Proof
Clearly we may assume that $p \geq 1$.
It is sufficient to show that if $(Z, g) \in \aBCycle^{p+d}(X)$,
then $(f_*(Z), f_*(g)) \in \aBCycle^{p}(Y)$.
By the definition of $L^1$-arithmetic cycles,
$g$ and $dd^c(g) + \delta_{Z(\CC)}$ are
represented by $L^1$-forms.
Thus, by Proposition~\ref{prop:push:forward:B:pq},
there is an $\omega \in \LocInt(\Omega_{Y(\CC)}^{p,p})$ with
\[
f_*\left( dd^c(g) + \delta_{Z(\CC)} \right) = [\omega].
\]
On the other hand, 
\[
f_*\left( dd^c(g) + \delta_{Z(\CC)} \right) =
dd^c(f_*(g)) + \delta_{f_*(Z(\CC))}.
\]
Moreover, by Proposition~\ref{prop:push:forward:B:pq},
$f_*(g)$ is represented by an $L^1$-form on $Y(\CC)$.
Thus, $(f_*(Z), f_*(g))$ is an element of $\aBCycle^p(Y)$.
\QED

\subsection{Scalar product for arithmetic $L^1$-cycles and arithmetic $D$-cycles}
\renewcommand{\theequation}{\arabic{section}.\arabic{subsection}.\arabic{Theorem}}
\setcounter{Theorem}{0}

Let $X$ be a regular arithmetic variety. 
The purpose of this subsection is to give a scalar product on
$\aDChow^{*}(X)_{\QQ} = \bigoplus_{p \geq 0} \aDChow^p(X)_{\QQ}$ by
the arithmetic Chow ring 
$\aChow^{*}(X)_{\QQ} = \bigoplus_{p \geq 0} \aChow^p(X)_{\QQ}$.
Roughly speaking, the scalar product is defined by
\[
(Y, f) \cdot (Z, g) =
( Y \cap Z, f \wedge \delta_Z + \omega((Y,f)) \wedge g)
\]
for $(Y, f) \in \aCycle^p(X)$ and $(Z, g) \in \aDCycle^q(X)$.
This definition, however, works only under the assumption
that $Y$ and $Z$ intersect properly.
Usually, by using Chow's moving lemma,
we can avoid the above assumption.
This is rather complicated, so that in this paper we try to 
use the standard arithmetic intersection theory to define
the scalar product.

Let $x \in \aChow^p(X)$, $(Z, g) \in \aDCycle^q(X)$, and
$g_Z$ a Green current for $Z$.
First we shall check that
\[
x \cdot [(Z, g_Z)] + [(0, \omega(x) \wedge (g-g_Z))]
\]
in $\aDChow^{p+q}(X)_{\QQ}$ does not depend on the choice of $g_Z$.
For, let $g'_Z$ be another Green current for $Z$.
Then, there are $\eta \in A^{p-1,p-1}(X(\CC))$, and 
$v \in \Image(\partial) + \Image(\bar{\partial})$ with
$g'_Z = g_Z + [\eta] + v$.
Then, since $[(0, [\eta] + v)] \in \aChow^p(X)$,
\begin{align*}
x \cdot [(Z, g'_Z)] + [(0, \omega(x) \wedge (g-g'_Z))] & =
x \cdot [(Z, g_Z)] + x \cdot [(0, [\eta] + v)] \\
& \qquad\qquad + [(0, \omega(x) \wedge (g-g_Z-[\eta]-v))] \\
& = x \cdot [(Z, g_Z)] + [(0, \omega(x) \wedge ([\eta] + v))] \\
& \qquad\qquad + [(0, \omega(x) \wedge (g-g_Z-[\eta]-v))] \\
& = x \cdot [(Z, g_Z)] + [(0, \omega(x) \wedge (g-g_Z))].
\end{align*}
Thus, we have the bilinear homomorphism
\[
\aChow^p(X) \times \aDCycle^q(X) \to \aDChow^{p+q}(X)_{\QQ}
\]
given by $x \cdot (Z, g) = x \cdot [(Z, g_Z)] + [(0, \omega(x) \wedge (g-g_Z))]$.
Moreover, if $(Z, g) \in \widehat{R}^q(X)$,
then, by \cite[Theorem~4.2.3]{GSArInt}, $x \cdot (Z, g) = 0$ in $\aChow^{p+q}(X)_{\QQ}$.
Thus, the above induces
\addtocounter{Theorem}{1}
\begin{equation}
\label{eqn:scalar:product:1}
\aChow^p(X) \otimes \aDChow^q(X) \to \aDChow^{p+q}(X)_{\QQ},
\end{equation}
which may give rises to a natural scalar product of
$\aDChow^{*}(X)_{\QQ}$ over
the arithmetic Chow ring $\aChow^{*}(X)_{\QQ}$.
To see that this is actually a scalar product, we need to check
that
\[
(x \cdot y) \cdot z = x \cdot (y \cdot z)
\]
for all $x \in \aChow^p(X)$, $y \in \aChow^q(X)$ and $z \in \aDChow^r(X)$.
If $z \in \aChow^r(X)$, then this is nothing more than
the associativity of the product of 
the arithmetic Chow ring (cf. \cite[Theorem~4.2.3]{GSArInt}).
Thus, we may assume that $z = [(0, g)]$ for some $g \in D^{r-1,r-1}(X(\CC))$.
Then, since 
\[
(x \cdot y) \cdot z = [(0, \omega(x \cdot y) \wedge g)] = 
[(0, \omega(x) \wedge \omega(y) \wedge g)]
\]
and
\[
x \cdot (y \cdot z) = x \cdot [(0, \omega(y) \wedge g)] = 
[(0, \omega(x) \wedge \omega(y) \wedge g)],
\]
we can see $(x \cdot y) \cdot z = x \cdot (y \cdot z)$.
Therefore, we get the natural scalar product.

Moreover, (\ref{eqn:scalar:product:1}) induces
\addtocounter{Theorem}{1}
\begin{equation}
\label{eqn:scalar:product:2}
\aChow^p(X) \otimes \aBChow^q(X) \to \aBChow^{p+q}(X)_{\QQ}.
\end{equation}
Indeed, if $(Z, g) \in \aBCycle^{q}(X)$ and $g_Z$ is a Green form of $Z$,
then,
\[
x \cdot [(Z, g)] = x \cdot [(Z, g_Z)] + [(0, \omega(x) \wedge (g - g_Z))].
\]
Thus, in order to see that $x \cdot [(Z, g)] \in \aBChow^{p+q}(X)_{\QQ}$,
it is sufficient to check that
\[
\begin{cases}
\omega(x) \wedge (g - g_Z) 
\in \LocInt(\Omega_{X(\CC)}^{p+q-1,p+q-1}), \\
dd^c \left( \omega(x) \wedge (g - g_Z) \right) 
\in \LocInt(\Omega_{X(\CC)}^{p+q,p+q}).
\end{cases}
\]
The first assertion is obvious because $g$ and $g_Z$ are $L^1$-forms.
Further, we can easily see the second assertion because
\[
dd^c \left( \omega(x) \wedge (g - g_Z) \right) =
\pm \omega(x) \wedge dd^c(g - g_Z) = 
\pm \omega(x) \wedge (\omega(g) - \omega(g_Z)).
\]

Gathering all observations, 
we can conclude the following proposition,
which is a generalization of \cite[Theorem~4.2.3]{GSArInt}.

\begin{Proposition}
\label{prop:module:structure:arith:D:cycle}
$\aBChow^{*}(X)_{\QQ}$ and $\aDChow^{*}(X)_{\QQ}$ has a natural
module structure over the
arithmetic Chow ring $\aChow^{*}(X)_{\QQ}$.
\end{Proposition}

Moreover, we have the following projection formula.

\begin{Proposition}
\label{prop:projection:formula:regular:smooth}
Let $f : X \to Y$ be a proper morphism of regular arithmetic varieties such that
$f_{\QQ} : X_{\QQ} \to Y_{\QQ}$ is smooth.
Then, for any $\alpha \in \aChow^p(Y)$ and $\beta \in \aChow_{L^1}^q(X)$,
\[
f_* (f^*(\alpha) \cdot \beta) = \alpha \cdot f_* (\beta) 
\] 
in $\aChow_{L^1}^{p + q -(\dim X - \dim Y)}(Y)_{\QQ}$.
\end{Proposition}

\Proof
If $\alpha \in \aChow^p(Y)$ and $\beta \in \aChow^q(X)$,
then this is well known (cf. \cite{GSArInt}). Thus,
we may assume that 
$\beta = (0,[\phi]) \in \aCycle_{L^1}^q(Y)$. 
Then 
\begin{align*}
f_* (f^*(\alpha) \cdot \beta) 
& = f_* ((0,\omega (f^*(\alpha)) \wedge [\phi]) \\
& = (0, [f_*\left( \omega (f^*(\alpha) \wedge \phi) \right)]).
\end{align*} 
On the other hand, 
\begin{equation*}
\alpha \cdot f_* (\beta) 
= \alpha \cdot (0,[f_*(\phi)]) 
= (0, \omega(\alpha) \wedge [f_*(\phi)]).
\end{equation*} 
Since $f_* (\omega (f^*(\alpha))) = \omega(\alpha)$, 
we have proven the projection formula. 
\QED
\renewcommand{\theequation}{\arabic{section}.\arabic{subsection}.\arabic{Theorem}.\arabic{Claim}}
\subsection{Scalar product, revisited (singular case)}
\setcounter{Theorem}{0}

Let $X$ be an arithmetic variety. 
Here $X$ is not necessarily regular.
Let $\Rat_X$ be the sheaf of rational functions on $X$.
We denote $H^0(X, \Rat_X^{\times}/\OO_X^{\times})$ by $\Div(X)$.
An element of $\Div(X)$ is called {\em a Cartier divisor on $X$}.
For a Cartier divisor $D$ on $X$, we can assign a divisor $[D] \in \Cycle^1(X)$
in a natural way. This gives rise to the homomorphism 
\[
c_X : \Div(X) \to \Cycle^1(X).
\]
Note that $c_X$ is neither injective nor surjective
in general. An exact sequence
\[
1 \to \OO_X^{\times} \to \Rat_X^{\times} \to \Rat_X^{\times}/\OO_X^{\times} \to 1
\]
induces to the homomorphism $\Div(X) \to H^1(X, \OO_X^{\times})$.
For a Cartier divisor $D$ on $X$, the image of $D$ by
the above homomorphism induces the line bundle on $X$.
We denote this line bundle by $\OO_X(D)$.

Here we set
\[
\aDiv(X) = \{ (D, g) \mid
\text{$D \in \Div(X)$ and $g$ is a Green function for $D(\CC)$ on $X(\CC)$} \}.
\]
Similarly, we can define
$\aBDiv(X)$ and $\aDDiv(X)$.
The homomorphism $c_X : \Div(X) \to \Cycle^1(X)$ gives rise to
the homomorphism $\hat{c}_X : \aDiv(X) \to \aCycle^1(X)$.
Then, we define $\aPic(X)$, $\aBPic(X)$, and $\aDPic(X)$ as follows.
\[
\begin{cases}
\aPic(X) = \aDiv(X)/\hat{c}_X^{-1}(\widehat{R}^1(X)), \\
\aBPic(X) = \aBDiv(X)/\hat{c}_X^{-1}(\widehat{R}^1(X)), \\
\aDPic(X) = \aDDiv(X)/\hat{c}_X^{-1}(\widehat{R}^1(X)).
\end{cases}
\]
Note that if $X$ is regular, then
\[
\aPic(X) = \aChow^1(X),
\quad
\aBPic(X) = \aBChow^1(X)
\quad\text{and}\quad
\aDPic(X) = \aDChow^1(X).
\]

Let $(E, h)$ be a Hermitian vector bundle on $X$. Then,
by virtue of \cite[Theorem~4]{GSRR},
we have a cap product of $\achernch(E, h)$ on $\aChow^{*}(X)_{\QQ}$,
i.e., a homomorphism
$\aChow^{*}(X)_{\QQ} \to \aChow^{*}(X)_{\QQ}$
given by $x \mapsto \achernch(E, h) \cap x$
for $x \in \aChow^{*}(X)_{\QQ}$.
In the same way as before, we can see that
the above homomorphism extends to
\[
\aDChow^{*}(X)_{\QQ} \to \aDChow^{*}(X)_{\QQ}
\quad\text{and}\quad
\aBChow^{*}(X)_{\QQ} \to \aBChow^{*}(X)_{\QQ}
\]
as follows. If $(Z, g) \in \aDCycle^p(X)$ and $g_Z$ is a Green current of $Z$,
then
\[
\achernch(E, h) \cap (Z, g) = \achernch(E, h) \cap (Z, g_Z) +
a(\operatorname{ch}(E, h) \wedge (g - g_Z)).
\]
In particular,
we have intersection pairings
\[
\aPic(X)_{\QQ} \otimes \aDChow^{p}(X)_{\QQ} \to \aDChow^{p+1}(X)_{\QQ}
\quad\text{and}\quad
\aPic(X)_{\QQ} \otimes \aBChow^{p}(X)_{\QQ} \to \aBChow^{p+1}(X)_{\QQ}.
\]
For simplicity, the cap product ``$\cap$'' is denoted by the dot ``$\cdot$''.
Note that
\[
\aPic(X)_{\QQ} \otimes \aDChow^{p}(X)_{\QQ} \to \aDChow^{p+1}(X)_{\QQ}
\]
is actually defined by
\[
(D, g) \cdot (Z, f) = (D \cdot Z, g \wedge \delta_Z + \omega(g) \wedge f)
\]
if $D$ and $Z$ intersect properly.
Then, we have the following projection formula.

\begin{Proposition}
\label{prop:projection:formula:line:bundle}
Let $f : X \to Y$ be a proper morphism of arithmetic varieties. 
Let $(L, h)$ be a Hermitian line bundle on $Y$, and
$z \in \aDChow^p(X)$. Then
\[
f_*(\acherncl_1 (f^*L, f^*h) \cdot z) = \acherncl_1(L, h) \cdot f_*(z).
\]
\end{Proposition}

\Proof
Let $(Z, g)$ be a representative of $z$.
Clearly, we may assume that $Z$ is reduced and irreducible.
We set $T = f(Z)$ and $\pi = \rest{f}{Z} : Z \to T$. 
Let $s$ be a rational section of $\rest{L}{T}$.
Then, $\pi^*(s)$ gives rise to a rational section of $\rest{f^*(L)}{Z} = 
\pi^* \left( \rest{L}{T} \right)$.
Thus, $\acherncl_1 (f^*L, f^*h) \cdot z$ can be represented by
\[
\left( \zero(\pi^*(s)), 
\left[ -\log \pi^*\left( \rest{h}{T} \right)(\pi^*(s), \pi^*(s)) \right]
+ c_1(f^*L, f^*h) \wedge g \right),
\]
where 
$\left[ -\log \pi^*\left( \rest{h}{T} \right)(\pi^*(s), \pi^*(s)) \right]$
is the current given by
\[
\left[ -\log \pi^*\left( \rest{h}{T} \right)(\pi^*(s), \pi^*(s)) \right](\phi)
= \int_{Z(\CC)}
\left( - \log \pi^*\left( \rest{h}{T} \right)(\pi^*(s), \pi^*(s)) \right) \phi.
\]
If we set
\[
\deg(\pi) =
\begin{cases}
0 & \text{if $\dim T < \dim Z$} \\
\deg(Z \to T) & \text{if $\dim T = \dim Z$,}
\end{cases}
\]
then
\begin{align*}
\int_{Z(\CC)}
\left( - \log \pi^*\left( \rest{h}{T} \right)(\pi^*(s), \pi^*(s)) \right) 
f^*(\psi)
& = \int_{Z(\CC)} \pi^* \left( \left(
-\log \left( \rest{h}{T} \right) (s, s) \right) \psi \right) \\
& = \deg(\pi) \int_{T(\CC)} \left( -\log \left( \rest{h}{T} \right) (s, s) \right) \psi
\end{align*}
for a $C^{\infty}$-form $\psi$ on $Y(\CC)$.
Thus, we have
\[
f_* \left[ -\log \pi^*\left( \rest{h}{T} \right)(\pi^*(s), \pi^*(s)) \right]
= \deg(\pi) \left[ -\log \left( \rest{h}{T} \right) (s, s) \right].
\]
Therefore,
\begin{align*}
f_*(\acherncl_1 (f^*L, f^*h) \cdot z)
& =
\left( \deg(\pi) \zero(s), \deg(\pi) \left[ -\log \left( \rest{h}{T} \right) (s, s) \right] +
c_1(L, h) \wedge f_*(g) \right) \\
& = \acherncl_1 (L, h) \cdot (\deg(\pi) T, f_*(g)) = \acherncl_1 (L, h) \cdot f_*(z).
\end{align*}
Hence, we get our proposition.
\QED

Let $Z$ be a quasi-projective integral scheme over $\ZZ$.
Then, by virtue of Hironaka's resolution of singularities \cite{Hiro},
there is a proper birational morphism $\mu : Z' \to Z$
of quasi-projective integral schemes over $\ZZ$ such that
$Z'_{\QQ}$ is non-singular.
The above $\mu : Z' \to Z$ is called a {\em generic resolution of
singularities of $Z$}.
As a corollary of the above projection formula,
we have the following proposition.

\begin{Proposition}
\label{prop:formula:restriction:intersection}
Let $X$ be a arithmetic variety, and
$\overline{L}_1 = (L_1, h_1), \ldots, \overline{L}_{n} = (L_{n}, h_{n})$
be Hermitian line bundles on $X$.
Let $(Z, g)$ be an arithmetic $D$-cycle on $X$, and
$Z = a_1 Z_1 + \cdots + a_r Z_r$ the irreducible decomposition as cycles.
For each $i$, let $\tau_i : Z'_i \to Z_i$ be a proper birational morphism
of quasi-projective integral schemes. We assume that if $Z_i$ is horizontal with respect to
$X \to \Spec(\ZZ)$, then $\tau_i$ is a generic resolution of
singularities of $Z_i$. Then, we have
\[
\acherncl_1(\overline{L}_1) \cdots \acherncl_1(\overline{L}_n) \cdot (Z, g) =
\sum_{i=1}^r a_i {\mu_i}_* \left( 
\acherncl_1(\mu_i^* \overline{L}_1) \cdots \acherncl_1(\mu_i^* \overline{L}_n) 
\right) +
a(c_1(\overline{L}_1) \wedge \cdots \wedge c_1(\overline{L}_n) \wedge g)
\]
in $\aDChow^*(X)_{\QQ}$,
where $\mu_i$ is the composition of
$Z'_i \overset{\tau_i}{\longrightarrow} Z_i \hookrightarrow X$
for each $i$.
\end{Proposition}

\Proof
We prove this proposition by induction on $n$.
First, let us consider the case $n = 1$.
Clearly we may assume that $Z$ is integral, i.e., $Z = Z_1$.
Let $h_1$ be the Hermitian metric of $\overline{L}_1$, and
$s$ a rational section of $\rest{L_1}{Z}$.
Then,
\[
\left(\zero(s), 
   - \log (\rest{h_1}{Z})(s, s) + c_1(\overline{L}_1) \wedge g \right) 
=
\left(\zero(s), 
   - \log (\rest{h_1}{Z})(s, s) \right)  + a(c_1(\overline{L}_1) \wedge g)
\]
is a representative of $\acherncl_1(\overline{L}_1) \cdot (Z, g)$.
Moreover, 
\[
\left(\zero(\tau_1^*(s)), 
   - \log \tau_1^*(\rest{h_1}{Z})(\tau_1^*(s), \tau_1^*(s)) \right)
\]
is a representative of $\acherncl_1(\mu_1^* \overline{L}_1)$.
Hence, we have our assertion in the case $n=1$ because
\[
\left({\mu_1}_* (\zero(\tau_1^*(s)), 
   - \log \tau_1^*(\rest{h_1}{Z})(\tau_1^*(s), \tau_1^*(s)) \right)
= (\zero(s), - \log (\rest{h_1}{Z})(s, s)).
\]
Thus, we may assume that $n > 1$. 
Therefore, using Proposition~\ref{prop:projection:formula:line:bundle} and
hypothesis of induction,
\begin{align*}
\acherncl_1(\overline{L}_1) \cdots \acherncl_1(\overline{L}_n) \cdot (Z, g) & =
\acherncl_1(\overline{L}_1) \cdot 
\left( \acherncl_1(\overline{L}_2) \cdots \acherncl_1(\overline{L}_n) \cdot (Z, g) \right) \\
& = 
\sum_{i=1}^r a_i \acherncl_1(\overline{L}_1) {\mu_i}_* \left( 
\acherncl_1(\mu_i^* \overline{L}_2) \cdots \acherncl_1(\mu_i^* \overline{L}_n) 
\right) + \\
& \qquad\qquad\qquad\qquad
\acherncl_1(\overline{L}_1) a(c_1(\overline{L}_2) \wedge \cdots \wedge c_1(\overline{L}_n) \wedge g) \\
& =
\sum_{i=1}^r a_i {\mu_i}_* \left(
\acherncl_1(\mu_i^* \overline{L}_1) \cdots \acherncl_1(\mu_i^* \overline{L}_n) 
\right) +
a(c_1(\overline{L}_1) \wedge \cdots \wedge c_1(\overline{L}_n) \wedge g).
\end{align*}
\QED

\subsection{Injectivity of $i^*$}
\setcounter{Theorem}{0}
Let $X$ be an arithmetic variety, $U$ a non-empty
Zariski open set of $X$, and
$i : U \to X$ the inclusion map.
Then, there is a natural homomorphism
\[
i^* : \aBCycle^1(X) \to \aBCycle^1(U)
\]
given by $i^*(D, g) = (\rest{D}{U}, \rest{g}{U(\CC)})$.
Since $i^* \left( \widehat{(f)} \right) = \widehat{(\rest{f}{U})}$
for any non-zero rational functions $f$ on $X$, the above
induces the homomorphism
\[
i^* : \aBChow^1(X) \to \aBChow^1(U).
\]
Then, we have the following useful lemma.

\begin{Lemma}
\label{lem:criterion:linear:equiv:B:cycle}
If $X \setminus U$ does not contain any irreducible components of
fibers of $X \to \Spec(\ZZ)$, then
\[
i^* : \aBChow^1(X) \to \aBChow^1(U).
\]
is injective. In particular,
$i^* : \aBChow^1(X)_{\QQ} \to \aBChow^1(U)_{\QQ}$
is injective.
\end{Lemma}

\Proof
Suppose that $i^*(\alpha) = 0$ for some $\alpha \in \aBChow^1(X)$.
Let $(D,g) \in \aBCycle^1(X)$ be a representative of $\alpha$.
Since $i^*(\alpha) = 0$, there is a non-zero rational function
$f$ on $X$ with
\[
(\rest{D}{U}, \rest{g}{U(\CC)}) = (\rest{(f)}{U}, \rest{-[\log |f|^2]}{U(\CC)}).
\]
Pick up $\phi \in \LocInt(X(\CC))$ with $g = [\phi]$.
Then, the above implies that
$\rest{[\phi]}{U(\CC)} = \rest{-[\log |f|^2]}{U(\CC)}$.
Thus, $\phi = - \log |f|^2 \ (\alev)$.
Therefore, we have
\addtocounter{Claim}{1}
\begin{equation}
\label{eqn:1:lem:criterion:linear:equiv:B:cycle}
g = [\phi] = - [\log |f|^2].
\end{equation}
Here, $dd^c(g) + \delta_{D(\CC)} = [h]$ for some 
$h \in \LocInt(\Omega_{X(\CC)}^{1,1})$ and
$dd^c(-[\log |f|^2]) + \delta_{(f)(\CC)} = 0$.
Thus, by (\ref{eqn:1:lem:criterion:linear:equiv:B:cycle}),
$\delta_{D(\CC)} - \delta_{(f)(\CC)} =[h]$.
This shows us that
$h = 0 \ (\alev)$ over
$X(\CC) \setminus \left( \Supp(D(\CC)) \cup \Supp((f)(\CC)) \right)$.
Hence $h = 0 \ (\alev)$ on $X(\CC)$.
Therefore, we have $D(\CC) = (f)(\CC)$, which implies
$D = (f)$ on $X_{\QQ}$. 
Thus, $D - (f)$ is a linear combination of
irreducible divisors lying on finite fibers.
On the other hand, $D = (f)$ on $U$ and
$X \setminus U$ does not contain any irreducible components of
fibers.
Therefore, $D = (f)$. 
Hence $\alpha = 0$ because $(D, g) = \widehat{(f)}$. 
\QED

\section{Weakly positive arithmetic divisors}

\subsection{Generalized metrics}
\setcounter{Theorem}{0}
\label{subsec:gen:metric}
Let $X$ be a smooth algebraic scheme over $\CC$ 
and $L$ a line bundle on $X$. 
We say $h$ is {\em a generalized metric on $L$} if there is a
$C^{\infty}$ Hermitian metric $h_0$ of $L$ over $X$ and
$\varphi \in \LocInt(X)$ with $h = e^{\varphi}h_0$.

To see when a metric of a line bundle defined over 
a dense Zariski open set
extends to a generalized metric, the following lemma is useful.
 
\begin{Lemma}
\label{lem:criterion:gen:metric}
Let $X$ be a smooth algebraic variety over $\CC$ and
$L$ a line bundle on $X$. Let $U$ be a non-empty Zariski open set of $X$ and
$h$ a $C^{\infty}$ Hermitian metric of $L$ over $U$.
We fix a non-zero rational section $s$ of $L$.
Then, $h$ extends to a generalized metric of $L$ on $X$ if and only if
$\log h(s, s) \in \LocInt(X)$.
\end{Lemma}

\Proof
If $h$ extends to a generalized metric of $L$ on $X$, then
$\log h(s, s) \in \LocInt(X)$ by the definition of generalized metrics.
Conversely, we assume that $\log h(s, s) \in \LocInt(X)$.
Let $h_0$ be a $C^{\infty}$ Hermitian metric of $L$ over $X$. Here we consider
the function $\phi$ given by
\[
  \phi = \frac{h(s,s)}{h_0(s, s)}.
\]
Let $y \in U$ and $\omega$ be a local frame of $L$ around $y$.
If we set $s = f \omega$ for some meromorphic function $f$ around $y$, then
\[
\phi = \frac{h(s,s)}{h_0(s, s)} = 
\frac{|f|^2 h(\omega, \omega)}{|f|^2 h_0(\omega, \omega)} =
\frac{h(\omega, \omega)}{h_0(\omega, \omega)}.
\]
This shows us that $\phi$ is a positive $C^{\infty}$ function on $U$ and
$h = \phi h_0$ over $U$.
On the other hand, 
\[
\log \phi = \log h(s, s) - \log h_0(s, s).
\]
Here since $\log h(s, s), \log h_0(s, s) \in \LocInt(X)$, 
we have $\log \phi \in \LocInt(X)$.
Thus, if we set $\varphi = \log \phi$, then $\varphi \in \LocInt(X)$ and
$h = e^{\varphi} h_0$.
\QED

\subsection{Arithmetic $D$-divisors and generalized metrics}
\label{subsec:arith:div:gen:metric}
\setcounter{Theorem}{0}
Let $X$ be an arithmetic variety, $L$ a line bundle on $X$,
and $h$ a generalized metric of $L$ on $X(\CC)$ with 
$F_{\infty}^*(h) = \overline{h} \ (\alev)$.
We would like to define $\acherncl_1(L, h)$ as an element
of $\aDChow^1(X)$.
Let $s, s'$ be two non-zero rational sections of $L$, and
$f$ a non-zero rational function on $X$
with $s' = fs$. 
Then, it is easy to see that
\[
(\zero(s'), [-\log h(s', s')]) =
(\zero(s), [-\log h(s,s)]) + \widehat{(f)}
\]
in $\aDCycle^1(X)$.
Thus, we may define $\acherncl_1(L, h)$ as the class of
$(\zero(s), [-\log h(s,s)])$ in $\aDChow^1(X)$.

Let us consider the homomorphism
\[
\omega : \aDCycle^p(X) \to  D^{p,p}(X(\CC))
\]
given by $\omega(Z, g) = dd^c(g) + \delta_{Z(\CC)}$.
Since $\omega\left( \widehat{R}^{p}(X) \right) = 0$,
the above homomorphism induces the homomorphism
$\aDChow^p(X) \to  D^{p,p}(X(\CC))$.
Hence, we get the homomorphism
$\aDChow^p(X)_{\QQ} \to  D^{p,p}(X(\CC))$
because $D^{p,p}(X(\CC))$ has no torsion.
By abuse of notation, 
we denote this homomorphism by $\omega$.

\begin{Proposition}
\label{prop:B:cycle:produce:hermitian:line:bundle}
Let $X$ be an arithmetic variety,
$(Z, [\phi]) \in \aDDiv(X)$ with
$\phi \in \LocInt(X(\CC))$, and
$1$ a rational section of $\OO_X(Z)$
with $\zero(1) = Z$.
Then, there is a unique generalized
metric $h$ of $\OO_X(Z)$ such that 
$F_{\infty}^*(h) = \overline{h} \ (\alev)$ and
$[-\log h(1, 1)] = [\phi]$.
\rom{(}Here uniqueness of $h$ means that
if $h'$ is another generalized metric with
the same property, then
$h = h' \ (\alev)$.\rom{)} Moreover,
$\omega(Z, [\phi])$
is $C^{\infty}$ around $x \in X(\CC)$ if and only if
$h$ is $C^{\infty}$ around $x$.
We denote this line bundle
$(\OO_X(Z), h)$ with the generalized metric $h$ by $\OO_Z((Z, [\phi]))$.
With this notation, 
for $(Z_1, [\phi_1]), (Z_2, [\phi_2]) \in \aDDiv(X)$
with $\phi_1, \phi_2 \in \LocInt(X(\CC))$, 
if $(Z_1, [\phi_1]) \sim (Z_2, [\phi_2])$, then
$\OO_X((Z_1, [\phi_1]))$ is isometric to
$\OO_X((Z_2, [\phi_2]))$ at every point around which
$\omega(Z_1, [\phi_1]) = \omega(Z_2, [\phi_2])$ is a
$C^{\infty}$ form.
\end{Proposition}

\Proof
First, let us see uniqueness.
Let $h$ and $h'$ be generalized metrics
of $\OO_X(Z)$ with
$[-\log h(1, 1)] = [-\log h'(1, 1)] = [\phi]$.
Take $a \in \LocInt(X(\CC))$
with $h' = e^a h$. Then,
by our assumption, $a = 0  \ (\alev)$.
Thus, $h = h'  \ (\alev)$.

Pick up an arbitrary point $x \in X(\CC)$.
Let $s$ be a local basis of $\OO_X(Z)$ around $x$.
Then, there is a non-zero 
rational rational function $f$ on $X(\CC)$
with $1 = f s$. Let us consider
\[
\exp(-\phi - \log |f|^2)
\]
around $x$. 
Let $s'$ be a another local basis of $\OO_X(Z)$
around $x$. We set $s' = us$ and $1 = f' s'$.
Then,
\[
\exp(-\phi - \log |f'|^2) =
\exp(-\phi - \log |f/u|^2) =
|u|^2 \exp(-\phi - \log |f|^2),
\]
which means that if we define the generalized
metric $h$ by
\[
 h(s, s) = \exp(-\phi - \log |f|^2),
\]
then $h$ is patched globally, and $h$ is a generalized metric
by Lemma~\ref{lem:criterion:gen:metric}.
Moreover,
\[
-\log h(1,1) = -\log h(fs, fs)
= -\log \left( |f|^2 h(s,s) \right) = \phi.
\]
Here, since $F_{\infty}^*(\phi) = \phi \ (\alev)$,
we can see $F_{\infty}^*(h) = \overline{h} \ (\alev)$.
Thus, we can construct our desired metric.

We set $b = \omega(Z, [\phi]) \in D^{1,1}(X(\CC))$.
Then, since $1 = fs$ around $x$,
we have $Z = (f)$ around $x$. Thus,
since $dd^c([\phi]) + \delta_{Z(\CC)} = b$ and
$dd^c(-[\log |f|^2]) + \delta_{(f)} = 0$,
\[
dd^c(-[\phi + \log |f|^2]) =
\delta_{Z(\CC)} - b - \delta_{(f)}
= - b
\]
around $x$. Therefore,
\begin{align*}
\text{$h$ is $C^{\infty}$ around $x$}
& \Longleftrightarrow
\text{$-\phi - \log |f|^2$ is $C^{\infty}$ around $x$} \\
& \Longleftrightarrow
\text{$dd^c(-[\phi + \log |f|^2])$ is $C^{\infty}$
around $x$} \qquad
\text{($\because$ \cite[Theorem~1.2.2]{GSArInt})}\\
& \Longleftrightarrow
\text{$b$ is $C^{\infty}$ around $x$}
\end{align*}

Finally, let us consider the last assertion.
By our assumption, there is a rational function $z$ on
$X$ such that
\[
(Z_1, [\phi_1]) = (Z_2, [\phi_2]) + \widehat{(z)}.
\]
Then,
$Z_1 = Z_2 + (z)$ and
$\phi_1 = \phi_2 - \log |z|^2$.
Let us consider the homomorphism
$\alpha : \OO_X(Z_1) \to \OO_X(Z_2)$
defined by $\alpha(s) = zs$.
Then, $\alpha$ is an isomorphism.
Let $1$ be the unit in the rational function field of $X$.
Then, $1$ gives rise to canonical rational sections of
$\OO_X(Z_1)$ and $\OO_X(Z_2)$.
Let $x$ be a point of $X(\CC)$ such that
$\omega(Z_1, [\phi_1])$ is $C^{\infty}$
around $x$,
and $s$ a local basis of
$\OO_X(Z_1)$ around $x$.
Then, $\alpha(s) = zs$
is a local basis of $\OO_X(Z_2)$ around
$x$. Choose a rational function $f$ with
$1 = fs$. Then, $1 = z^{-1}f\alpha(s)$.
Thus, if $h_1$ and $h_2$ are metrics of
$\OO_X((Z_1, [\phi_1]))$ and $\OO_X((Z_2, [\phi_2]))$
respectively, then
\[
h_1(s, s) = \exp(-\phi_1 -\log |f|^2) =
\exp( -\phi_2 - \log |z^{-1}f|^2) =
h_2(\alpha(s), \alpha(s))
\]
Hence, $\alpha$ is an isometry.
\QED

\subsection{Semi-ampleness and small sections}
\setcounter{Theorem}{0}
Let $X$ be an arithmetic variety and $S$ a subset of $X(\CC)$.
We set
\[
\aBChow^1(X;S)_{\QQ} = \{ \alpha \in \aBChow^1(X)_{\QQ} \mid
\text{$\omega(\alpha)$ is $C^{\infty}$ around $z$ for
all $z \in S$} \}.
\]
In the same way, we can define
$\aBChow^1(X;S)$, $\aBCycle^1(X;S)$,
$\aBCycle^1(X;S)_{\QQ}$, $\aBDiv(X;S)$,
$\aBDiv(X;S)_{\QQ}$, $\aBPic(X;S)$ and
$\aBPic(X;S)_{\QQ}$.
Let $x$ be a closed point of $X_{\QQ}$.
An element $\alpha$ of $\aBChow^1(X;S)_{\QQ}$
is said to be {\em semi-ample at $x$ with respect to $S$} if
there are a positive integer $n$ and $(E, g) \in \aBCycle^1(X;S)$
with the following properties:
\begin{enumerate}
\renewcommand{\labelenumi}{(\alph{enumi})}
\item
$E$ is effective and $x \not\in \Supp(E)$.

\item
$g(z) \geq 0$ for each $z \in S$.
(Note that $g(z)$ might be $\infty$.)

\item
$n \alpha$ coincides with $(E,g)$ in $\aBChow^1(X;S)_{\QQ}$.
\end{enumerate}
We remark that $\alpha \in \aBChow^1(X;S)_{\QQ}$ by the condition (c).
Moreover, it is easy to see that
if $\alpha_1$ and $\alpha_2$ are semi-ample at $x$ with respect to $S$,
so is $\lambda \alpha_1 + \mu \alpha_2$ for all 
non-negative rational numbers
$\lambda$ and $\mu$.

In terms of the natural action of $\Gal(\overline{\QQ}/\QQ)$
on $X(\overline{\QQ})$, we can consider the orbit
$O_{\Gal(\overline{\QQ}/\QQ)}(x)$ of $x$.
If we fix an embedding $\overline{\QQ} \to \CC$, we have
an injection $X(\overline{\QQ}) \to X(\CC)$.
It is easy to see that the image of $O_{\Gal(\overline{\QQ}/\QQ)}(x)$
does not depend on the choice of the embedding $\overline{\QQ} \to \CC$.
By abuse of notation, we denote this image by 
$O_{\Gal(\overline{\QQ}/\QQ)}(x)$.
Then, $O_{\Gal(\overline{\QQ}/\QQ)}(x)$ is one of the examples
of $S$.

\medskip
Let $U$ be a Zariski open set of $X$,
and $F$ a coherent $\OO_X$-module such that $F$ is locally free over $U$.
Let $h_F$ be a $C^{\infty}$ Hermitian metric of $F$ over $U(\CC)$.
We assume that $S \subseteq U(\CC)$.
For a closed point $x$ of $U_{\QQ}$, we say
$(F, h_F)$ is 
{\em generated by small sections at $x$ with respect to $S$}
if there are sections $s_1, \ldots, s_n \in H^0(X, F)$ such that
$F_x$ is generated by $s_1, \ldots, s_n$, and that
$h_F(s_i, s_i)(z) \leq 1$ for all $1 \leq i \leq n$ and
$z \in S$.

\begin{Proposition}
\label{prop:comparion:semiample:gen:small:sec}
We assume that $S \subseteq U(\CC)$.
For an element $(Z, g)$ of $\aBDiv(X;S)$,
$(Z,g)$ is semi-ample at $x$ with respect to $S$ if and only if
there is a positive integer $n$ such that
$\OO_X(n(Z,g))$ is generated by small sections at $x$
with respect to $S$.
\end{Proposition}

\Proof
First, we assume that $(Z,g)$ is semi-ample at $x$ with respect to $S$.
Then, there is $(E, f) \in \aBCycle^1(X;S)$ and a positive integer $n$
such that $n(Z, g) \sim (E, f)$,
$E$ is effective, $x \not\in \Supp(E)$, and
$f(z) \geq 0$ for each $z \in S$.
Note that $E$ is a Cartier divisor.
Then, by Proposition~\ref{prop:B:cycle:produce:hermitian:line:bundle},
$\OO_X(n(Z,g)) \simeq \OO_X((E, f))$.
Moreover, if $h$ is the metric of $\OO_X((E, f))$ and
$1$ is the canonical section of $\OO_X(E)$ with $\zero(1) = E$,
then $-\log(h(1,1)) = f$. 
Here $f(z) \geq 0$ for each $z \in S$.
Thus, $h(1,1)(z) \leq 1$ for each $z \in S$.
Therefore, $\OO_X((E, f))$ is generated by small sections at $x$
with respect to $S$.

Next we assume that $\OO_X(n(Z,g))$ is generated by small sections
at $x$ with respect to $S$ for some  positive integer $n$.
Then, there is a section $s$ of $\OO_X(nZ)$ such that
$h(s,s)(z) \leq 1$ for each $z \in S$.
Thus, if we set $E = \zero(s)$ and $f = -\log h(s,s)$,
then we can see $(Z,g)$ is semi-ample at $x$ with respect to $S$.
\QED

\begin{Proposition}
\label{prop:finite:sup:imply:globalsup}
Let $U$ be a Zariski open set of $X$,
and $L$ a line bundle on $X$.
Let $h$ be a $C^{\infty}$ Hermitian metric of $L$ over $U(\CC)$.
Fix a closed point $x$ of $U_{\QQ}$.
If $X$ is projective over $\ZZ$, then the followings are equivalent.
\begin{enumerate}
\renewcommand{\labelenumi}{(\arabic{enumi})}
\item
$(L, h)$ is generated by small sections at $x$ with respect to $U(\CC)$.

\item
$(L, h)$ is generated by small sections at $x$ with respect to any
finite subsets $S$ of $U(\CC)$.
\end{enumerate}
\end{Proposition}

\Proof
Clearly, (1) implies (2).
So we assume (2).
First of all, we can easily find $z_1, \ldots, z_n \in U(\CC)$ such that,
for any $s \in H^0(X(\CC), L_{\CC})$, if $s(z_1) = \cdots = s(z_n) = 0$,
then $s = 0$.
Thus, if we set
\[
\Vert s \Vert = \sqrt{h(s, s)(z_1)} + \cdots + \sqrt{h(s, s)(z_n)}
\]
for each $s \in H^0(X(\CC), L_{\CC})$, then $\Vert \ \Vert$ gives rise to
a norm on $H^0(X(\CC), L_{\CC})$.
Here we set
\[
B_z = \{ s \in H^0(X,L) \mid \text{$h(s,s)(z) \leq 1$} \}
\]
for each $z \in U(\CC)$.
Then, since $H^0(X, L)$ is a discrete subgroup of $H^0(X(\CC), L_{\CC})$,
$\bigcap_{i=1}^{n} B_{z_i}$ is a finite set.
Thus, adding finite points $z_{n+1}, \ldots, z_{N} \in U(\CC)$ 
to $z_1, \ldots, z_n$ if necessary, we have
\[
\bigcap_{z \in U(\CC)} B_z = \bigcap_{i=1}^N B_{z_i}.
\]
By our assumption, there is a section $s \in H^0(X, L)$ such that
$s(x) \not= 0$ and $h(s, s)(z_i) \leq 1$ for all $i=1, \ldots, N$.
Then, $s \in \bigcap_{i=1}^N B_{z_i} = \bigcap_{z \in U(\CC)} B_z$.
Thus, we get (2).
\QED

\subsection{Restriction to arithmetic curves}
\setcounter{Theorem}{0}
Let $X$ be an arithmetic variety,
$S$ a subset of $X(\CC)$,
$x$ a closed point of $X_{\QQ}$,
$K$ the residue field of $x$,
and $O_K$ the ring of integers in $K$.
We assume that 
the orbit of $x$ by $\Gal(\overline{\QQ}/\QQ)$ is contained in $S$,
namely, $O_{\Gal(\overline{\QQ}/\QQ)}(x) \subseteq S$, and
that the canonical morphism $\Spec(K) \to X$ 
induced by $x$ extends to
$\tilde{x} : \Spec(O_K) \to X$.

\begin{Proposition}
\label{prop:homo:ChowB:to:Chow:on:OK}
There is a natural homomorphism
\[
\tilde{x}^* : \aBPic(X;S)_{\QQ} \to \aChow^1(\Spec(O_K))_{\QQ}
\]
such that the restriction of $\tilde{x}^*$ to
$\aPic(X)_{\QQ}$ coincides with the usual pull-back
homomorphism.
\end{Proposition}

\Proof
Let $\alpha \in \aBPic(X;S)_{\QQ}$.
Choose $(Z, g) \in \aBDiv(X;S)$ and a positive integer $e$ such that
the class of $(1/e)(Z, g)$ in $\aBPic(X;S)_{\QQ}$
coincides with $\alpha$.
We would like to define $\tilde{x}^*(\alpha)$ by 
\[
(1/e) \acherncl_1\left( \tilde{x}^*(\OO_X((Z, g))) \right).
\]
For this purpose, we need to check that the above does not depend on
the choice $(Z, g)$ and $e$.
Let $(Z', g')$ and $e'$ be another $L^1$-cycle of codimension $1$ and
positive integer such that
the class of $(1/e')(Z', g')$ is $\alpha$.
Then, there is a positive integer $d$ such that
$de'(Z, g) \sim de(Z', g')$.
Thus, by Proposition~\ref{prop:B:cycle:produce:hermitian:line:bundle},
$\OO_Z(de'(Z, g))$ is isometric to $\OO_Z(de(Z', g'))$.
Hence,
\begin{align*}
de' \acherncl_1\left( \tilde{x}^*(\OO_X((Z, g))) \right) & =
\acherncl_1\left( \tilde{x}^*(\OO_X(de'(Z, g))) \right) \\
& = \acherncl_1\left( \tilde{x}^*(\OO_X(de(Z', g'))) \right) \\
& = de \acherncl_1\left( \tilde{x}^*(\OO_X((Z', g'))) \right).
\end{align*}
Therefore,
\[
(1/e) \acherncl_1\left( \tilde{x}^*(\OO_X((Z, g))) \right)
= (1/e') \acherncl_1\left( \tilde{x}^*(\OO_X((Z', g'))) \right).
\]
Thus, we can define $\tilde{x}^*$.
\QED

\subsection{Weak positivity of arithmetic $L^1$-divisors}
\label{subsec:wp:div}
\setcounter{Theorem}{0}
Let $X$ be an arithmetic variety,
$S$ a subset of $X(\CC)$, and $x$ a closed point of $X_{\QQ}$.
Let $\alpha \in \aBChow^1(X;S)_{\QQ}$ and
$\{ \alpha_n \}_{n=1}^{\infty}$ a sequence of
elements of $\aBChow^1(X;S)_{\QQ}$.
We say $\alpha$ is the limit of $\{ \alpha_n \}_{n=1}^{\infty}$
as $n \to \infty$, denoted by 
${\displaystyle \alpha = \lim_{n \to \infty} \alpha_n}$,
if there are 
(1) $Z_1, \ldots, Z_{l_1} \in \aBChow^1(X;S)_{\QQ}$,
(2) $g_1, \ldots, g_{l_2} \in \LocInt(X(\CC))$ with
$a(g_j) \in \aBChow^1(X;S)_{\QQ}$ for all $j$,
(3) sequences $\{ a_n^{1} \}_{n=1}^{\infty}, \ldots,
\{ a_n^{l_1} \}_{n=1}^{\infty}$ of rational numbers, and 
(4) sequences $\{ b_n^{1} \}_{n=1}^{\infty}, \ldots,
\{ b_n^{l_2} \}_{n=1}^{\infty}$ of real numbers
with the following properties:
\begin{enumerate}
\renewcommand{\labelenumi}{(\alph{enumi})}
\item
$l_1$ and $l_2$ does not depend on $n$.

\item
${\displaystyle \lim_{n \to \infty} a_n^{i} =
\lim_{n \to \infty} b_n^{j} = 0}$ for all $1 \leq i \leq l_1$ and
$1 \leq j \leq l_2$.

\item
${\displaystyle
\alpha = \alpha_n + 
\sum_{i=1}^{l_1} a_n^{i}Z_i + \sum_{j=1}^{l_2} a(b_n^{j} g_j)}$
in $\aBChow^1(X;S)_{\QQ}$ for all $n$.
\end{enumerate}
It is easy to see that if
${\displaystyle \alpha = \lim_{n \to \infty} \alpha_n}$ and
${\displaystyle \beta = \lim_{n \to \infty} \beta_n}$
in $\aBChow^1(X;S)_{\QQ}$, then
${\displaystyle \alpha + \beta = 
\lim_{n \to \infty} (\alpha_n + \beta_n)}$.

An element $\alpha$ of $\aBChow^1(X;S)_{\QQ}$
is said to be {\em weakly positive at $x$ with respect to $S$} if 
it is the limit of semi-ample $\QQ$-cycles at $x$ with respect to $S$,
i.e., there is a sequence $\{ \alpha_n \}_{n=1}^{\infty}$ of
elements of $\aBChow^1(X;S)_{\QQ}$ such that
$\alpha_n$'s are semi-ample at $x$ with respect to $S$ and
${\displaystyle \alpha = \lim_{n \to \infty} \alpha_n}$.

Let $K$ be the residue field of $x$ and
$O_K$ the ring of integers in $K$.
We assume that $O_{\Gal(\overline{\QQ}/\QQ)}(x) \subseteq S$, and
the canonical morphism $\Spec(K) \to X$ induced by $x$ extends to
$\tilde{x} : \Spec(O_K) \to X$.
Then, we have the following proposition.

\begin{Proposition}
\label{prop:non:negative:wp:div:via:beta}
If $X$ is regular and an element $\alpha$ of $\aBChow^1(X;S)_{\QQ}$ is
weakly positive at $x$ with respect to $S$, 
then $\adeg(\tilde{x}^*(\alpha)) \geq 0$.
\end{Proposition}

\Proof
Take $Z_1, \ldots, Z_{l_1}$,
$g_1, \ldots, g_{l_2}$,
$\{ a_n^{1} \}_{n=1}^{\infty}, \ldots,
\{ a_n^{l_1} \}_{n=1}^{\infty}$,
$\{ b_n^{1} \}_{n=1}^{\infty}, \ldots,
\{ b_n^{l_2} \}_{n=1}^{\infty}$, and
$\{ \alpha_n \}_{n=1}^{\infty}$ as in the previous
definition of weak positive arithmetic divisors.
Then,
\[
\adeg ( \tilde{x}^*(\alpha) ) = \adeg ( \tilde{x}^*(\alpha_n) ) + 
\sum_{i=1}^{l_1} a_n^{i} \adeg (\tilde{x}^*(Z_i)) + 
\sum_{j=1}^{l_2} b_n^{j} \adeg(\tilde{x}^*a(g_j)).
\]
Thus, since ${\displaystyle \lim_{n \to \infty} a_n^{i} =
\lim_{n \to \infty} b_n^{j} = 0}$ for all $1 \leq i \leq l_1$ and
$1 \leq j \leq l_2$ and $\adeg ( \tilde{x}^*(\alpha_n) ) \geq 0$
for all $n$, we have $\adeg ( \tilde{x}^*(\alpha) ) \geq 0$.
\QED

\subsection{Characterization of weak positivity}
\setcounter{Theorem}{0}
Let $X$ be a regular arithmetic variety, $S$ a subset of $X(\CC)$, and
$x$ a closed point of $X_{\QQ}$.
For an element $\alpha \in \aBChow^1(X)_{\QQ}$,
we say $\alpha$ is {\em ample at $x$ with respect to $S$}
if there are $(A, f) \in \aBCycle^1(X;S)$ and
a positive integer $n$ such that
$A$ is an effective and ample Cartier divisor on $X$, $x \not\in \Supp(A)$,
$f(z) > 0$ for all $z \in S$, and
$n \alpha$ is equal to $(A, f)$ in $\aBChow^1(X)_{\QQ}$.

First, let us consider the case where $X = \Spec(O_K)$.

\begin{Proposition}
\label{prop:wp:for:curve}
We assume that $X = \Spec(O_K)$, 
$x$ is the generic of $X$, and
$S = X(\CC)$.
For an element $\alpha \in \aChow^1(X;S)_{\QQ}$, we have the following.
\begin{enumerate}
\renewcommand{\labelenumi}{(\arabic{enumi})}
\item
$\alpha$ is ample at $x$ with respect to $S$ if and only if
$\adeg(\alpha) > 0$.

\item
$\alpha$ is weakly positive at $x$ with respect to $S$ if and only if
$\adeg(\alpha) \geq 0$.  
\end{enumerate}
\end{Proposition}

\Proof
(1)
Clearly, if $\alpha$ is ample at $x$ with respect to $S$, 
then $\adeg(\alpha) > 0$. Conversely,
we assume that $\adeg(\alpha) > 0$.
We take a positive integer $e$ and a Hermitian line bundle
$(L, h)$ on $X$ such that $\acherncl_1(L, h) = e \alpha$.
Then, $\adeg(L, h) > 0$. Thus, by virtue of Riemann-Roch formula and
Minkowski's theorem, there are a positive integer $n$ and
a non-zero element $s$ of $L^{\otimes n}$ 
with $(h^{\otimes n})(s, s)(z) < 1$ for all
$z \in S$. Thus, $\alpha$ is ample at $x$ with respect to $S$.

\medskip
(2)
First, we assume that $\alpha$ is weakly positive at $x$
with respect to $S$.
Then, by Proposition~\ref{prop:non:negative:wp:div:via:beta},
$\adeg(\alpha) \geq 0$.
Next, we assume that $\adeg(\alpha) \geq 0$. 
Let $\beta$ be an element of $\aChow^1(X;S)_{\QQ}$ such that
$\beta$ is ample at $x$ with respect to $S$.
Then, for any positive integer $n$,
$\adeg(\alpha + (1/n) \beta) > 0$. Thus, $\alpha + (1/n)\beta$
is ample at $x$ with respect to $S$ by (1). 
Hence, $\alpha$ is weakly positive at $x$ with respect to $S$.
\QED

Before starting a general case,
let us consider the following lemma.

\begin{Lemma}
\label{lem:make:semiample:by:ample}
We assume that $S$ is compact.
Let $\alpha$ be an element of $\aBChow^1(X)_{\QQ}$ such that
$\alpha$ is ample at $x$ with respect to $S$.
Then, we have the following.
\begin{enumerate}
\renewcommand{\labelenumi}{(\arabic{enumi})}
\item
Let $\beta$ be an element of $\aBChow^1(X;S)_{\QQ}$.
Then, there is a positive number $\epsilon_0$ such that
$\alpha + \epsilon \beta$ is semi-ample at $x$ with respect to $S$
for all rational numbers $\epsilon$ with $|\epsilon| \leq \epsilon_0$.

\item
Let $g$ be a locally integrable function on $X(\CC)$
with $a(g) \in \aBChow^1(X;S)_{\QQ}$. 
Then, there is a positive number
$\epsilon_0$ such that
$\alpha + a(\epsilon g)$ is semi-ample at $x$ with respect to $S$
for all real numbers $\epsilon$ with $|\epsilon| \leq \epsilon_0$.
\end{enumerate}
\end{Lemma}

\Proof
(1)
First, we claim that there is a positive number $t_0$ such that
$t \alpha + \beta$ is semi-ample at $x$ with respect to $S$ for
all rational numbers $t \geq t_0$.

Let us choose $(A, f) \in \aBCycle^1(X;S)$ and
a positive integer $n_0$ such that
$A$ is an effective and ample Cartier divisor on $X$, $x \not\in \Supp(A)$,
$f(z) > 0$ for all $z \in S$, and
$n_0 \alpha$ is equal to $(A, f)$ in $\aBChow^1(X)_{\QQ}$.
Moreover, we choose $(D, g) \in \aBCycle^1(X;S)$ and 
a positive integer $e$
such that $e \beta$ is equal to $(D, g)$ in
$\aBChow^1(X)_{\QQ}$.
Since $A$ is ample, there is a positive integer $n_1$
such that $\OO_X(n_1 A) \otimes \OO_X(D)$ is generated by global
sections at $x$. 
Thus, there are $(Z, h) \in \aBCycle^1(X;S)_{\QQ}$ such that
$Z$ is effective,
$x \not\in \Supp(Z)$ and $(Z, h) \sim n_1(A, f) + (D, g)$.

We  would like to find a positive integer $n_2$
with $n_2 f(z) + h(z) \geq 0$
for all $z \in S$.
Let $U$ be an open set of $X(\CC)$ such that $S \subseteq U$, and
$\omega(A, f)$ and $\omega(Z,h)$ are $C^{\infty}$ over $U$.
We set $\phi = \exp(-f)$ and $\psi = \exp(-h)$.
Then, $\phi$ and $\psi$ are continuous on $U$, and
$0 \leq \phi < 1$ on $S$.
Since $n_2 f + h = -\log(\phi^{n_2} \psi)$, it is sufficient to find
a positive integer $n_2$ with $\phi^{n_2} \psi \leq 1 $ on $S$.
If we set $a = \sup_{z \in S} \phi(z)$ and $b = \sup_{z \in S} \psi(z)$,
then $0 \leq a < 1$ and $0 \leq b$ because $S$ is compact.
Thus, there is a positive integer $n_2$ with $a^{n_2} b \leq 1$.
Therefore, $\phi^{n_1} \psi \leq 1$ on $S$.

Here we set 
$t_0 = (n_1+n_2)n_0e^{-1}$.
In order to see that $t\alpha + \beta$ is semi-ample at $x$
with respect to $S$ for $t \geq t_0$,
it is sufficient to show that
$(n_1+n_2)n_0 \alpha + e \beta$
is semi-ample at $x$ with respect to $S$ because $e t \geq (n_1+n_2)n_0$.
Here
\begin{align*}
(n_1 + n_2)n_0 \alpha + e \beta 
& \sim n_2(A, f) + \left( n_1(A, f) + (D, g) \right) \\
& \sim n_2(A, f) + (Z, h) \\
& = (n_2 A + Z, n_2 f + h),
\end{align*}
$x \not\in \Supp(n_2 A + Z)$, and
$(n_2 f + h)(z) \geq 0$ for all $z \in S$.
Thus, $(n_1 + n_2)n_0 \alpha + e \beta$
is semi-ample at $x$ with respect to $S$.
Hence, we get our claim.

In the same way, we can find a positive number $t_1$ such that
$t \alpha - \beta$ is semi-ample with respect to $S$
for all $t \geq t_1$.
Thus, if we set $\epsilon_0 = \min \{ 1/t_0, 1/t_1 \}$, then we have (1).

\medskip
(2)
In the same way as in the proof of (1),
we can find a positive number $\epsilon_0$ such that
$(f + \epsilon n_0 g)(z) \geq 0$ for all $z \in S$ and
all real number $\epsilon$ with $|\epsilon| \leq \epsilon_0$.
Thus we have (2)
because $n_0( \alpha + a(\epsilon g)) \sim (A, f + \epsilon n_0 g)$.
\QED

\begin{Proposition}
\label{prop:characterization:wp:div}
We assume that $S$ is compact.
Let $\beta$ be an element of $\aBChow^1(X;S)_{\QQ}$.
Then the following are equivalent.
\begin{enumerate}
\renewcommand{\labelenumi}{(\arabic{enumi})}
\item
$\beta$ is weakly positive at $x$ with respect to $S$.

\item
$\beta + \alpha$ is semi-ample at $x$ with respect to $S$
for any ample $\alpha \in \aBChow^1(X;S)_{\QQ}$ at $x$ 
with respect to $S$.
\end{enumerate}
\end{Proposition}

\Proof
(1) $\Longrightarrow$ (2):
Since $\beta$ is weakly positive at $x$ with respect to $S$,
there is a sequence of $\{ \beta_n \}$ such that
$\beta_n \in \aBChow^1(X;S)_{\QQ}$, 
$\beta_n$'s are semi-ample at $x$ with respect to $S$,
and $\lim_{n \to \infty} \beta_n = \beta$.
Take $Z_1, \ldots, Z_{l_1}$,
$g_1, \ldots, g_{l_2}$,
$\{ a_n^{1} \}_{n=1}^{\infty}, \ldots,
\{ a_n^{l_1} \}_{n=1}^{\infty}$, and
$\{ b_n^{1} \}_{n=1}^{\infty}, \ldots,
\{ b_n^{l_2} \}_{n=1}^{\infty}$
as in the definition of the limit in $\aBChow^1(X;S)_{\QQ}$.
Then, by Lemma~\ref{lem:make:semiample:by:ample},
there is a positive number $\epsilon_0$ such that
$\alpha + \epsilon Z_i$'s are semi-ample at $x$ with respect to $S$
for all rational numbers $\epsilon$ with $|\epsilon| \leq \epsilon_0$,
and $\alpha + a(\epsilon g_j)$'s are semi-ample at $x$ 
with respect to $S$
for all real numbers $\epsilon$ with $|\epsilon| \leq \epsilon_0$.
We choose $n$ such that $(l_1 + l_2)|a_n^{i}| \leq \epsilon_0$ and
$(l_1 + l_2)|b_n^{j}| \leq \epsilon_0$ for all $i$ and $j$.
Then,
\[
\beta + \alpha = \beta_n + 
\sum_{i=1}^{l_1} \frac{\alpha + (l_1 + l_2)a_n^{i} Z_i}{l_1 + l_2}  +
\sum_{j=1}^{l_2} \frac{\alpha + a((l_1+ l_2)b_n^{j} g_j)}{l_1 + l_2}.
\]
Here, $\alpha + (l_1 + l_2)a_n^{i} Z_i$ and
$\alpha + a((l_1+ l_2)b_n^{j} g_j)$ are semi-ample $x$ 
with respect to $S$.
Thus, we get the direction (1) $\Longrightarrow$ (2).

\medskip
(2) $\Longrightarrow$ (1):
Let $\alpha$ be an element of $\aBChow^1(X)_{\QQ}$ such that
$\alpha$ is ample at $x$ with respect to $S$.
We set $\beta_n = \beta + (1/n)\alpha$.
Then, by our assumption, $\beta_n$ is semi-ample at $x$
with respect to $S$.
Further, $\beta = \lim_{n \to \infty} \beta_n$.
\QED

\subsection{Small sections via generically finite morphisms}
\setcounter{Theorem}{0}
Let $g : V \to U$ be a proper and \'{e}tale morphism
of complex manifolds.
Let $(E, h)$ be a Hermitian vector bundle on $V$.
Then, a Hermitian metric $g_*(h)$ of $g_*(E)$
is defined by
\[
g_*(h)(s, t)(y) = \sum_{x \in g^{-1}(y)} h(s, t)(x)
\]
for any $y \in U$ and $s, t \in g_*(E)_y$.

\begin{Proposition}
\label{prop:find:small:section}
Let $X$ be a scheme such that 
every connected component of $X$ is a arithmetic variety.
Let $Y$ be a regular arithmetic variety, and
$g : X \to Y$ a proper and generically finite morphism such that
every connected component of $X$ maps surjectively to $Y$.
Let $U$ be a Zariski open set of $Y$ such that
$g$ is \'{e}tale over $U$. 
Let $S$ be a subset of $U(\CC)$ and
$y$ a closed point of $U_{\QQ}$.
Then, we have the following.
\begin{enumerate}
\renewcommand{\labelenumi}{(\arabic{enumi})}
\item
Let $\phi : E \to Q$ be a homomorphism of
coherent $\OO_X$-modules such that
$\phi$ is surjective over $g^{-1}(U)$, and
$E$ and $Q$ are locally free over $g^{-1}(U)$.
Let $h_E$ be a $C^{\infty}$ Hermitian metric
of $E$ over $g^{-1}(U)(\CC)$, and $h_Q$ the quotient metric
of $Q$ induced by $h_E$.
If $(g_*(E), g_*(h_E))$ is generated by small sections at $y$
with respect to $S$, then
$(g_*(Q), g_*(h_Q))$ is generated by small sections at $y$
with respect to $S$.

\item
Let $E_1$ and $E_2$ be coherent $\OO_X$-modules such that
$E_1$ and $E_2$ are locally free over $g^{-1}(U)$.
Let $h_1$ and $h_2$ be $C^{\infty}$ Hermitian metrics
of $E_1$ and $E_2$ over $g^{-1}(U)(\CC)$.
If $(g_*(E_1), g_*(h_1))$ and $(g_*(E_2), g_*(h_2))$ 
are generated by small sections
at $y$ with respect to $S$,
then so is $(g_*(E_1 \otimes E_2), g_*(h_1 \otimes h_2))$.

\item
Let $E$ be a coherent $\OO_X$-module such that
$E$ is locally free over $g^{-1}(U)$.
Let $h_E$ be a $C^{\infty}$ Hermitian metric of $E$ over
$g^{-1}(U)(\CC)$.
If $(g_*(E), g_*(h_E))$ 
is generated by small sections at $y$ with respect
to $S$, then
$(g_*(\Sym^n(E)), g_*(\Sym^n(h_E)))$ 
is generated by small sections at $y$
with respect to $S$.
\rom{(}For the definition of $\Sym^n(h_E)$, 
see \S\rom{\ref{subsec:formula:chern:sym:power}}.\rom{)}

\item
Let $F$ be a coherent $\OO_Y$-module such that
$F$ is locally free over $U$.
Let $h_F$ be a $C^{\infty}$ Hermitian metric of $F$ over $U(\CC)$.
Since $\rest{\det(F)}{U}$ is canonically isomorphic to
$\det(\rest{F}{U})$, $\det(h_F)$ gives rise to a $C^{\infty}$
Hermitian metric of $\det(F)$ over $U(\CC)$.
If $(F, h_F)$ is generated by small sections at $y$ with respect to $S$, 
then so is $(\det(F), \det(h_F))$. 
\end{enumerate}
\end{Proposition}

\Proof
(1) 
By our assumption,
$g_*(\phi) : g_*(E) \to g_*(Q)$ is surjective over $U$.
Let $s_1, \ldots, s_l \in H^0(Y, g_*(E)) = H^0(X, E)$
such that $g_*(E)_{y}$ is generated by
$s_1, \ldots, s_l$, and that $g_*(h_E)(s_i, s_i)(z) \leq 1$
for all $i$ and $z \in S$.
Then, $g_*(Q)_y$ is generated by 
$g_*(\phi)(s_1), \ldots, g_*(\phi)(s_l)$.
Moreover, 
by the definition of the quotient metric $h_Q$,
\[
g_*(h_Q)(g_*(\phi)(s_i), g_*(\phi)(s_i))(z) =
\sum_{x \in g^{-1}(z)}
h_Q(\phi(s_i), \phi(s_i))(x) \leq
\sum_{x \in g^{-1}(z)} 
h_E(s_i, s_i)(x) \leq 1
\]
for all $z \in S$.
Hence, $g_*(Q)$ is generated by small sections at $y$
with respect to $S$.

\medskip
(2)
Since $g$ is \'{e}tale over $U$,
$\alpha : g_*(E_1) \otimes g_*(E_2) \to g_*(E_1 \otimes E_2)$
is surjective over $U$.
By our assumption, there are
$s_1, \ldots, s_l \in H^0(Y, g_*(E_1))$ and
$t_1, \ldots, t_m \in H^0(Y, g_*(E_2))$ such that
$g_*(E_1)_y$ (resp. $g_*(E_2)_y$)
is generated by $s_1, \ldots, s_l$ (resp. $t_1, \ldots, t_m$),
and that $g_*(h_1)(s_i, s_i)(z) \leq 1$ and
$g_*(h_2)(t_j, t_j)(z) \leq 1$ for
all $i$, $j$ and $z \in S$.
Then, $g_*(E_1 \otimes E_2)_y$ is generated by
$\{ \alpha(s_i \otimes t_j) \}_{i,j}$.
Moreover,
{\allowdisplaybreaks
\begin{align*}
g_*(h_1 \otimes h_2)(\alpha(s_i \otimes t_j),
\alpha(s_i \otimes t_j))(z) & =
\sum_{x \in g^{-1}(z)} (h_1 \otimes h_2)(s_i \otimes t_j,
s_i \otimes t_j)(x) \\
& = \sum_{x \in g^{-1}(z)}
h_1(s_i, s_i)(x) h_2(t_j,t_j)(x) \\
& \leq 
\left( \sum_{x \in g^{-1}(z)} h_1(s_i, s_i)(x) \right)
\left( \sum_{x \in g^{-1}(z)} h_2(t_j, t_j)(x) \right) \\
& \leq 1
\end{align*}}
for all $z \in S$.
Thus, we get (2).

\medskip
(3) This is a consequence of (1) and (2).

\medskip
(4)
Let $r$ be the rank of $F$.
Since $F$ is generated by small sections at $y$
with respect to $S$,
there are $s_1, \ldots, s_r \in H^0(Y, F)$ such that
$F \otimes \kappa(y)$ is generated by $s_1, \ldots, s_r$ and
$h(s_i, s_i)(z) \leq 1$
for all $i$ and $z \in S$.
Let us consider an exact sequence:
\[
 0 \to F_{tor} \to F \to F/F_{tor} \to 0.
\]
Then, $\det(F) = \det(F/F_{tor}) \otimes \det (F_{tor})$.
Noting that $F_{tor} = 0$ on $U$,
let $g$ be a Hermitian metric of $\det(F/F_{tor})$ over $U(\CC)$
given by $\det(h_F)$.
Then, there is a Hermitian metric $k$ of $\det(F_{tor})$ over $U(\CC)$
such that 
$(\det(F), \det(h_F)) = (\det(F/F_{tor}), g) \otimes (\det (F_{tor}), k)$
over $U(\CC)$. If we identify $\det (F_{tor})$ with $\OO_{Y}$ over $U$,
$k$ is nothing more than the canonical metric of $\OO_{Y}$ over $U(\CC)$.

Let us fix a locally free sheaf $P$ on $Y$ and
a surjective homomorphism $P \to F_{tor}$.
Let $P'$ be the kernel of $P \to F_{tor}$.
Here $\left( \bigwedge^{\rank P'} P' \right)^{*}$
is an invertible sheaf on $Y$ because $Y$ is regular.
Thus we may identify $\det(F_{tor})$ with 
\[
\bigwedge^{\rank P} P \otimes 
\left( \bigwedge^{\rank P'} P' \right)^{*}.
\]
Further, a homomorphism
$\bigwedge^{\rank P'} P' \to \bigwedge^{\rank P} P$
induced by $P' \hookrightarrow P$ 
gives rise to
a non-zero section $t$ of $\det(F_{tor})$ because
\[
\Hom\left(\bigwedge^{\rank P'} P', \bigwedge^{\rank P} P\right) =
\Hom\left(\bigwedge^{\rank P'} P', \OO_Y\right) \otimes \bigwedge^{\rank P} P.
\]
Here $F_{tor} = 0$ on $U$. Thus,
$\det(F_{tor})$ is canonically isomorphic to $\OO_{Y}$ over $U$.
Since $P' = P$ over $U$, under the above isomorphism,
$t$ goes to the determinant of
$P' \overset{\operatorname{id}}{\longrightarrow} P$,
namely $1 \in \OO_{Y}$ over $U$.
Thus, $k(t,t)(z) = 1$ for each $z \in S$.

Let $\overline{s}_i$ be the image of $s_i$ in $F/F_{tor}$.
Then, $\overline{s}_1 \wedge \cdots \wedge \overline{s}_r$
gives rise to a section $s$ of $\det(F/F_{tor})$.
Thus, $s \otimes t$ is a section of $\det(F)$.
By our construction, $(s \otimes t)(y) \not= 0$.
Moreover, using Hadamard's inequality,
\begin{align*}
\det(h_F)(s \otimes t, s \otimes t)(z) & =
g(s, s)(z) \cdot k(t,t)(z) =
\det \left( h( s_i, s_j)(z) \right) \\
& \leq 
h(s_1, s_1)(z) \cdots h(s_r, s_r)(z) \leq 1
\end{align*}
for each $z \in S$.
Thus, we get (4).
\QED

\section{Arithmetic Riemann-Roch for generically finite morphisms}

\subsection{Quillen metric for generically finite morphisms}
\setcounter{Theorem}{0}
Before starting 
Proposition~\ref{prop:Quillen:metric:generalised:gen:finite:morph},
we recall the tensor product of two
matrices, which we will use in the proof.
For an 
$r\times r$ matrix $A=(a_{ij})$
and an
$n\times n$ matrix $B=(b_{kl})$, 
consider the following $rn\times rn$ matrix 
\begin{equation*}
\begin{pmatrix}
 a_{11}B & a_{12}B & \cdots & a_{1r}B \\
 a_{21}B & a_{22}B & \cdots & a_{2r}B \\
 \vdots  & \vdots  & \ddots & \vdots  \\
 a_{r1}B & a_{r2}B & \cdots & a_{rr}B \\
\end{pmatrix}.
\end{equation*}
This matrix, denoted by 
$A\otimes B$, 
is called {\em the tensor product of $A$ and $B$}. 
Then for 
$r \times r$ matrices $A,A'$ 
and
$n \times n$ matrices $B,B'$, 
we immediately see
\begin{equation*}
(A \otimes B) (A' \otimes B') = AA' \otimes BB', \cr
\det (A \otimes B) = (\det A)^n (\det B)^r.
\end{equation*}

Let $X$ be a smooth algebraic scheme over $\CC$, 
$Y$ a smooth algebraic variety over $\CC$, and
$f : X \to Y$ a proper and generically finite morphism.
We assume that every connected component of $X$ maps surjectively to
$Y$. Let $W$ be the maximal open set of $Y$ 
such that $f$ is \'{e}tale over there. 
Let $(E, h)$ be a Hermitian vector
bundle on $X$ such that on every connected component of $X$,
$E$ has the same rank $r$.

\begin{Proposition}
\label{prop:Quillen:metric:generalised:gen:finite:morph}
With notation and assumptions being as above, the Quillen metric 
$h_Q^{\overline E}$ on $\det Rf_*(E)$ over $W$ 
extends to a generalized metric on $\det Rf_*(E)$ over $Y$.
\end{Proposition}

\Proof
Let $n$ be the degree of $f$. 
Since $f$ is \'{e}tale over $W$, 
$f_*(E)$ is a locally free sheaf of rank $rn$ and 
$R^if_*(E) = 0$ for $i \ge 1$ over there. Thus 
\[
\det Rf_*(E) \vert {}_W = \bigwedge^{rn} f_*(E) \vert {}_W.
\]
If $y \in W$ is a complex point and 
$X_y=\{x_1,x_2,\cdots,x_n\}$ the fiber of $f$ over $y$,
then we have 
\[
\det Rf_*(E)_y = \det H^0(X_y,E).
\]
The Quillen metric on $\det Rf_*(E)$ over $W$ is defined as follows.
On $H{}^0(X_y,E)$ the $L^2$-metric is defined by the formula:
\[
h_{L^2}(s,t) = \sum_{\alpha=1}^n h(s,t)(x_\alpha),
\]
where $s,t \in H{}^0(X_y,E)$. 
This metric naturally induces the $L^2$-metric on 
$\det H{}^0(X_y,E)$. Since $X_y$ is zero-dimensional,
there is no need for zeta function regularization 
to obtain the Quillen metric. 
Thus the Quillen metric $h_Q^{\overline E}$ 
on $\det Rf_*(E) \vert {}_W$ is defined by 
the family of Hermitian line bundles 
$\{\det H^0(X_y,E)\}_{y \in W}$ 
with the induced $L^2$-metrics pointwisely.

To see that the Quillen metric over $W$
extends to a generalized metric over $Y$, let 
$s_1,s_2,\cdots,s_r$ be rational sections of $E$ such that 
at the generic point of every connected component of $X$,
they form a basis of $E$. Also let 
$\omega_1,\omega_2,\cdots,\omega_n$ be rational 
sections of $f_*(\OO_X)$ such that at the generic point 
they form a basis of 
$f_*(\OO_X)$. Since
\addtocounter{Claim}{1}
\begin{equation}
\label{eqn:2:prop:Quillen:metric:generalised:gen:finite:morph}
\det Rf_*(E) = \left(\bigwedge^{rn} (f_*(E))
               \right)^{**}.
\end{equation}
over $Y$, we can regard $\bigwedge_{ik}s_i\omega_k=s_1\omega_1\Wedge
s_1\omega_2\Wedge\cdots\Wedge
s_1\omega_n\Wedge\cdots\Wedge s_r\omega_n$ 
as a non-zero rational section of $\det Rf_*(E)$. 
Shrinking $W$, we can find a non-empty Zariski open set $W_0$ of $W$
such that $s_i$'s and $\omega_j$'s has no poles or zeros over $f^{-1}(W_0)$.

To proceed with our argument, we need the following lemma.

\begin{Lemma}
\label{lem:formula:for:Quillen:metric}
Let $L$ be the total quotient field of $X$, and
$K$ the function field of $Y$.
Then,
\[
\log h_Q^{\overline E} \left( \bigwedge_{ik}s_i\omega_k,
                       \bigwedge_{ik}s_i\omega_k \right) =
r \log \left| \det (\Tr_{L/K}(\omega_i \cdot \omega_j)) \right|
+ f_* \log \det (h(s_i, s_j))
\]
over $W_0$.
\end{Lemma}

\Proof
Let $y \in W_0$ be a complex point, and 
$\{x_1,x_2,\ldots,x_n\}$ the fiber of $f^{-1}(y)$ over $y$. Then, 
{\allowdisplaybreaks
\newcommand{\bigzerol}{\smash{\hbox{\Large 0}}}
\newcommand{\bigzerou}{\smash{\lower1.7ex\hbox{\Large 0}}}
\begin{align*}
 \log h_Q^{\overline E}
  \left( \bigwedge_{ik}s_i\omega_k,\bigwedge_{ik}s_i\omega_k \right)(y)
 & = \log \det 
      \left(
        \sum_{\alpha =1}^n
        h(s_i\omega_k,s_j\omega_l)(x_{\alpha})
      \right)_{ij,kl} \\
 & = \log \det 
      \left(
        \sum_{\alpha =1}^n
        \omega_k(x_\alpha) h(s_i,s_j)(x_{\alpha}) 
        \overline{\omega_l(x_\alpha)}
      \right)_{ij,kl} \\
 & = \log \det \left( 
      (I_r \otimes \Omega) 
      \begin{pmatrix}
        H(x_1)    &        &  \bigzerou  \\
                  & \ddots &             \\
        \bigzerol &        &  H(x_n)     \\
      \end{pmatrix}
      \overline{{}^t(I_r \otimes \Omega)}
               \right) \\
 & = \log \det \left\{
      \vert \det (\Omega) \vert ^{2r}
      \prod_{\alpha=1}^n 
      \det \left(h(s_i,s_j)(x_{\alpha})\right)_{ij}
                \right\} \\
 & = r \log \det \vert \det (\Omega) \vert ^{2} +
     \sum_{\alpha=1}^n 
     \log \det \left(h(s_i,s_j)\right)(x_{\alpha}), 
\end{align*}}
where 
$\Omega = (\omega_k(x_\alpha))_{k \alpha}$ 
and 
$H(x_\alpha) = (h(s_i,s_j)(x_{\alpha}))_{ij}$. 
On the other hand, we have  
\[
\sum_{\alpha=1}^n 
     \log \det \left(h(s_i,s_j)\right) (x_{\alpha})
= \left(f_* \log \det \left(h(s_i,s_j)\right)\right)(y).
\] 
Moreover, using the 
following Lemma~\ref{lemma:trace:global:local}, we have 
\[
\vert \det (\Omega) \vert ^{2} =
\vert \det (\Omega {}^t \Omega) \vert =
\left\vert \det \left(\sum_{\alpha=1}^n \omega_k(x_\alpha) 
                \omega_l(x_\alpha) \right)_{kl} \right\vert
= \vert \det \left(\Tr_{L/K}(\omega_k\cdot \omega_l)\right)_{kl} \vert. 
\]
Thus we get the lemma. 
\QED

\begin{Lemma}
\label{lemma:trace:global:local}
Let $f : \Spec(B) \to \Spec(A)$ 
be a finite \'{e}tale morphism of regular affine schemes. 
Let ${\mathfrak m}$ be the maximal ideal of $A$ and 
${\mathfrak n}_1,{\mathfrak n}_2,\cdots,{\mathfrak n}_n$ the prime ideals lying over 
${\mathfrak m}$. Assume that $\kappa ({\mathfrak m})$ is algebraically closed 
and hence $\kappa ({\mathfrak n}_i)$ is \rom{(}naturally\rom{)}
isomorphic to  $\kappa ({\mathfrak m})$ for each $1 \leq i \leq n$. 
Let $b$ be an element of $B$ and $b({\mathfrak n}_i)$ the value of
$b$ in $\kappa ({\mathfrak n}_i) \cong \kappa ({\mathfrak m})$.
Then 
\[
\Tr_{B/A}(b)({\mathfrak m}) = \sum_{i=1}^n b({\mathfrak n}_i)
\]
in $\kappa ({\mathfrak m})$,
where $\Tr_{B/A}(b)({\mathfrak m})$ is 
the value of $\Tr_{B/A}(b)$ in $\kappa ({\mathfrak m})$.
\end{Lemma}

\Proof
It is easy to see that every ${\mathfrak n}_i$ is the maximal ideal and that 
${\mathfrak m}B =  {\mathfrak n}_1{\mathfrak n}_2 \cdots {\mathfrak n}_n$. 
Let $\widehat{A}$ be the completion of $A$ with respect to ${\mathfrak m}$,  
$\widehat{B}$ the completion of $B$ with respect to ${\mathfrak m}B$, and 
$\widehat{B_i}$ the completion of $B$ 
with respect to ${\mathfrak n}_i$ for each $1 \leq i \leq n$. 
Then by Chinese remainder theorem, 
$\widehat{B} = \prod_{i=1}^n \widehat{B_i}$ as an $\widehat{A}$-algebra. 
Note that $\widehat{A} /{\mathfrak m} \widehat{A} =  \kappa ({\mathfrak m})$ 
and $\widehat{B_i} /{\mathfrak n}_i \widehat{B_i} =  \kappa ({\mathfrak n}_i)$. 
Since $\widehat{A} \to \widehat{B_i}$ is \'{e}tale and $\kappa ({\mathfrak m}) 
\cong \kappa ({\mathfrak n}_i)$, we have 
$\widehat{A} \cong \widehat{B_i}$.
Let $e_1=(1,0,\cdots,0),e_2=(0,1,\cdots,0),\cdots,e_n=(0,0,\cdots,1) 
\in \prod_{i=1}^n \widehat{B_i}
=\widehat{B}$ be a free basis of $\widehat{B}$ over 
$\widehat{A}$.  We put $b e_i = b_i e_i$ with $b_i \in
\widehat{B_i} \cong \widehat{A}$  for each $1 \leq i \leq n$. 
Then $b_i \equiv b({\mathfrak n}_i) \pmod{{\mathfrak n}_i}$. 
Now the lemma follows from 
\[
\Tr_{B/A}(b) = \Tr_{\widehat{B} / \widehat{A}}(b)
= \sum_{i=1}^n b_i
\]
in $\widehat{A}$.
\QED

Let us go back to the proof of
Proposition~\ref{prop:Quillen:metric:generalised:gen:finite:morph}.
Since $\rest{\det (\Tr_{L/K}(\omega_i \cdot \omega_j))}{W_0}$ extends to a rational
function $\det (\Tr_{L/K}(\omega_i \cdot \omega_j))$ on $Y$, 
\[
\log \left |
\det (\Tr_{L/K}(\omega_i \cdot \omega_j)) \right|
\in \LocInt(Y).
\]
Moreover, 
by Proposition~\ref{prop:push:forward:B:pq},
$f_* \log \det (h(s_i, s_j)) \in \LocInt(Y)$.
Thus, by Lemma~\ref{lem:formula:for:Quillen:metric},
\[
\rest{\log h_Q^{\overline E} 
\left( \bigwedge_{ik}s_i\omega_k, \bigwedge_{ik}s_i\omega_k \right)}{W_0}
\]
extends to a locally integrable function on $Y$.
Hence by Lemma~\ref{lem:criterion:gen:metric} 
the Quillen metric over $W$ 
extends to a generalized metric over $Y$.
\QED

\begin{Remark}
\label{rmk:Quillen:metric:continuous:finite:morphism}
In the above situation, Let $W'$ be a open set of $Y$ such that 
$f$ is flat and finite over there. 
Then the Quillen metric extends to a continuous function over $W'$
by the same formula as in 
(\ref{lem:formula:for:Quillen:metric})
\end{Remark}

\subsection{Riemann-Roch for generically finite morphisms}
\setcounter{Theorem}{0}
In this subsection, we formulate 
the arithmetic Riemann-Roch theorem 
for generically finite morphisms.  

\begin{Theorem}
\label{thm:arith:Riemann:Roch:gen:finite:morphism}
Let $X$ be a scheme such that every connected component of $X$ is an arithmetic variety.
Let $Y$ be a regular arithmetic variety, and
$f : X \to Y$ a proper and generically finite morphism such that
every connected component of $X$ maps surjectively to $Y$.
Let $(E, h)$ a Hermitian vector bundle on $X$ such that
on each connected component of $X$, $E$ has the same rank $r$. Then,
\[
\acherncl_1 \left( \det Rf_*(E), h_Q^{\overline{E}} \right) -
r \acherncl_1 \left( \det Rf_*(\OO_X), h_Q^{\overline{\OO}_X}
\right) \in 
\aBChow^1(Y)
\]
and
\[
\acherncl_1 \left( \det Rf_*(E), h_Q^{\overline{E}} \right) -
r \acherncl_1 \left( \det Rf_*(\OO_X), h_Q^{\overline{\OO}_X}
\right) = f_* \left( \acherncl_1 (E, h) \right)
\]
in $\aBChow^1(Y)_{\QQ}$,
where $h_Q^{\overline{E}}$ and $h_Q^{\overline{\OO}_X}$ are
the Quillen metric of $\det Rf_*(E)$ and $\det Rf_*(\OO_X)$ respectively.
\end{Theorem}

\Proof
Let $X = \coprod_{\alpha \in A} X_{\alpha}$ be the decomposition into connected
components of $X$. Since $f$ is proper, $A$ is a finite set.
We set $f_{\alpha} = \rest{f}{X_{\alpha}}$ and $(E_{\alpha}, h_{\alpha}) = \rest{(E, h)}{X_{\alpha}}$.
Then
\[
Rf_*(E) = \bigoplus_{\alpha \in A} R (f_{\alpha})_*(E_{\alpha}),\quad
Rf_*(\OO_X) = \bigoplus_{\alpha \in A} R (f_{\alpha})_*(\OO_{X_{\alpha}}),\quad
\acherncl_1 (E, h) = \sum_{\alpha \in A} \acherncl_1 (E_{\alpha}, h_{\alpha}).
\]
Hence we have the following:
\[
\begin{cases}
{\displaystyle \acherncl_1 \left( \det Rf_*(E), h_Q^{\overline{E}} \right) =
\sum_{\alpha \in A} \acherncl_1 \left( \det R(f_{\alpha})_*(E_{\alpha}), 
h_Q^{\overline{E}_{\alpha}} \right)}, \\
{\displaystyle \acherncl_1 \left( \det Rf_*(\OO_X), h_Q^{\overline{\OO}_X} \right) =
\sum_{\alpha \in A} \acherncl_1 \left( \det R(f_{\alpha})_*(\OO_{X_{\alpha}}), 
h_Q^{\overline{\OO}_{X_{\alpha}}} \right)}, \\
{\displaystyle f_* \left( \acherncl_1 (E, h) \right) =
\sum_{\alpha \in A} f_* \left( \acherncl_1 (E_{\alpha}, h_{\alpha}) \right)}.
\end{cases}
\]
Thus, we may assume that $X$ is connected, i.e.,
$X$ is an arithmetic variety.

\medskip
Let $K=K(Y)$ and $L=K(X)$ 
be the function fields of 
$Y$ and $X$ respectively. 
Let $n$ be the degree of $f$ 
and
$\omega_1,\omega_2,\cdots,\omega_n$ 
rational functions on $X$ 
such that at the generic point they form a basis of $K$-vector space $L$.
Further, let $s_1,s_2,\ldots,s_r$ 
be rational sections of $E$
such that at the generic point they form a basis of
$L$-vector space $E_L$. 
Then $s_1\omega_1\Wedge s_1\omega_2\Wedge\cdots\Wedge
s_1\omega_n\Wedge\cdots\Wedge s_r\omega_n$,
$s_1\Wedge\cdots\Wedge s_r$ 
and
$\omega_1\Wedge\cdots\Wedge\omega_n$ 
are non-zero rational sections of 
$\det f_*(E)$, 
$\det (E)$ 
and
$\det f_*(\OO_X)$ respectively.
Here we shall prove the following equality in
$\aDCycle^1(Y)$:
\addtocounter{Claim}{1} 
\begin{multline}
\label{eqn:desired:thm:arith:Riemann:Roch:gen:finite:morphism}
 \left(\zero \left(\bigwedge_{ik}s_i\omega_k\right),
   \left[-\log h_Q^{\overline E}
   \left(\bigwedge_{ik}s_i\omega_k,\bigwedge_{ik}s_i\omega_k\right)\right] 
 \right) \\
 - r
 \left(\zero\left(\bigwedge_k\omega_k\right),
  \left[-\log h_Q^{\overline \OO_X}
  \left(\bigwedge_{k}\omega_k,\bigwedge_{k}\omega_k\right)\right]
 \right) \\
 = f_*\left(\zero\left(\bigwedge_i s_i\right),
       \left[-\log \det h\left(\bigwedge_i s_i,\bigwedge_i s_i\right)\right]
      \right),
\end{multline}
where 
$\bigwedge_{ik}s_i\omega_k=s_1\omega_1\Wedge
s_1\omega_2\Wedge\cdots\Wedge
s_1\omega_n\Wedge\cdots\Wedge s_r\omega_n$,
$\bigwedge_k\omega_k
=\omega_1\Wedge\cdots\Wedge\omega_n$ and 
$\bigwedge_i s_i=s_1\Wedge\cdots\Wedge s_r$. 

First we shall show the equality of divisors.
Let $Y_0$ be the maximal Zariski open set
of $X$ such that $f$ is flat over $Y_0$.
Then, $\codim_Y(Y \setminus Y_0)\ge 2$ 
by \cite[III,Proposition~9.7]{Hartshorne}. 
Since $f$ is generically finite, 
$f$ is in fact finite over $Y_0$. 
Then $\Cycle^1(Y)=\Cycle^1(Y_0)$ and thus  
it suffices to prove the equality of divisors over $Y_0$. 
Since it suffices to prove it locally, let
$U=\Spec(A)$ 
be an affine open set of $Y_0$ 
and $f^{-1}(U)=\Spec(B)$ 
the open set of $X_0=f^{-1}(Y_0)$.
Shrinking $U$ if necessary,
we may assume that 
$B$ is a free $A$-module of rank $n$ 
and that 
$E$ is a free $B$-module of rank $r$. 
Let $d_1,d_2,\cdots,d_n$
be a basis of $B$ over $A$, 
and
$e_1,e_2,\cdots,e_r$ 
be a basis of $E$ over $B$.
Note that $K$ and $L$ are the
quotient fields of $A$ and $B$ respectively. 
In the following we freely identify a rational function (or
section) by the corresponding element at the generic point. 
In this sense, we set 
\begin{align*}
\omega_k = \sum_{l=1}^{n}a^{kl}d_l \quad (k=1,2,\cdots,n) \\
s_i = \sum_{j=1}^{r}\sigma_{ij}e_j \quad (i=1,2,\cdots,r),
\end{align*} 
where 
$a^{kl} \in K\,(1\le k,l\le n)$ 
and 
$\sigma_{ij}\in L\,(1\le i,j\le r)$.

For each 
$\sigma_{ij}(1\le i,j\le r)$, 
let
$T_{\sigma_{ij}} : L \to L$ 
be multiplication by
$\sigma_{ij}$. 
With respect to a basis
$\omega_1,\omega_2,\cdots,\omega_n$ 
of $L$ over $K$, 
$T_{\sigma_{ij}}$ gives rise to
the matrix
$(c_{ij}^{kl})_{1\le k,l\le n}\in M_n(K)$
defined by
\[
  \sigma_{ij}\omega_k
 =\sum_{l=1}^n c_{ij}^{kl}\omega_l 
 \quad (k=1,2,\cdots,n).
\]
We also denote this matrix by 
$T_{\sigma_{ij}}$.
Then,
{\allowdisplaybreaks
\begin{align*}
 \bigwedge_{ik}s_i\omega_k
 & = \bigwedge_{ik}\left(\sum_{j=1}^{r}
     \sigma_{ij}e_j\right)\omega_k 
   = \bigwedge_{ik}\left(\sum_{j=1}^{r}\sum_{l=1}^{n}
     c_{ij}^{kl}\right)e_j\omega_l \\
 & = \det (c_{ij}^{kl})_{ik,jl}
     \bigwedge_{ik}e_i\omega_k
   = \det (c_{ij}^{kl})_{ik,jl}
     \bigwedge_{ik}e_i\left(\sum_{l=1}^{n}a^{kl}d_l\right) \\
 & = \det (c_{ij}^{kl})_{ik,jl}
     \bigwedge_{ik}\left(\sum_{j=1}^{r}
     \delta_{ij}a^{kl}\right)e_j \omega_l
   = \det (c_{ij}^{kl})_{ik,jl}
     \det (\delta_{ij}a^{kl})_{ik,jl}
     \bigwedge_{ik}e_i d_l.
\end{align*}}
On the other hand, since the matrices 
$T_{\sigma_{ij}}$ and $T_{\sigma_{i'j'}}$ 
commute with each other, we have
{\allowdisplaybreaks
\begin{align*}
  \det (c_{ij}^{kl})_{ik,jl}
 & = \det 
  \begin{pmatrix}
   T_{\sigma_{11}} & T_{\sigma_{12}} & \cdots &
   T_{\sigma_{1r}} \\
   T_{\sigma_{21}} & T_{\sigma_{22}} & \cdots &
   T_{\sigma_{2r}} \\
   \vdots & \vdots & \ddots & \vdots \\
   T_{\sigma_{r1}} & T_{\sigma_{r2}} & \cdots &
   T_{\sigma_{rr}} \\
  \end{pmatrix}\\
 & = \det \left(
      \sum_{\tau\in \Perm_r}
       \sign(\tau)
      T_{\sigma_{1\tau(1)}}\cdot\cdots\cdot
      T_{\sigma_{r\tau(r)}}
          \right) \\
 & = \det(T_{\det(\sigma_{ij})_{ij}}) \\
 & = \Norm_{L/K}(\det(\sigma_{ij})_{ij}).
\end{align*}}
Moreover, we have
\begin{align*}
     \det (\delta_{ij}a^{kl})_{ik,jl}
 & = \det (I_r \otimes (a^{kl})_{kl}) \\
 & = (\det  (a^{kl})_{kl})^r.
\end{align*}
>From the above three equalities, 
$\zero\left(\bigwedge_{ik}s_i\omega_k\right)$
is given by the rational function
\[
       \Norm_{L/K}(\det(\sigma_{ij})_{ij})
       (\det  (a^{kl})_{kl})^r.
\]
Further
\[
\bigwedge_{i}s_i =  (\det  (\sigma_{ij})_{ij}) 
                      \bigwedge_{i} e_k
\quad\text{and}\quad
\bigwedge_{k}\omega_k =  (\det  (a^{kl})_{kl})
                      \bigwedge_{k} d_k.
\]
Hence we have 
\[
  \zero\left(\bigwedge_{ik}s_i\omega_k\right)
 -r\left(\zero\left(\bigwedge_k\omega_k\right)\right)
 = f_*\left(\zero\left(\bigwedge_i s_i\right)\right).
\]

Next we shall show the equality of currents. 
Since all the currents in the equality come from 
locally integrable functions by
Proposition~\ref{prop:push:forward:B:pq} and 
Proposition~\ref{prop:Quillen:metric:generalised:gen:finite:morph}, 
it suffices to show the equality over a non-empty Zariski open set of 
every connected component of $Y(\CC)$. 
So let $W_0$ 
be a non-empty Zariski open set of a connected component of $Y(\CC)$ 
such that $f_{\CC}$ is \'{e}tale 
and that 
$s_i\, (1 \le i \le r)$ or $\omega_k \,(1 \le k \le n)$ 
have no poles or zeroes over there.  
Then over $W_0$ 
all these currents are defined by $C^{\infty}$ functions. 
Let
$y \in Y(\CC)$ 
be a complex point 
and 
$x_1,x_2,\cdots,x_n$ 
be the fiber $f_{\CC}^{-1}(y)$ over $y$. 
>From the proof of Lemma~\ref{lem:formula:for:Quillen:metric}, 
as  $C^{\infty}$ functions around $y$,
{\allowdisplaybreaks
\begin{align*}
 -\log h_Q^{\overline E}
  \left(\bigwedge_{ik}s_i\omega_k,\bigwedge_{ik}s_i\omega_k\right)(y)
 & = -\log \det \left\{
      \vert \det (\Omega) \vert ^{2r}
      \prod_{\alpha=1}^n 
      \det \left(h(s_i,s_j)(x_{\alpha})\right)_{ij}
                \right\}, 
\end{align*}}
where 
$\Omega = (\omega_k(x_\alpha))_{k \alpha}$ 
and 
$H(x_\alpha) = (h(s_i,s_j)(x_{\alpha}))_{ij}$. Also, 
{\allowdisplaybreaks
\begin{align*}
 -\log h_Q^{\overline \OO_X}
  \left(\bigwedge_{k}\omega_k,\bigwedge_{k}\omega_k\right)(y)
 & = -\log \det 
      \vert \det (\Omega) \vert ^{2}.
\end{align*}}
On the other hand, by the definition of the push-forward $f_*$,
\begin{align*}
   f_* \left(
     -\log \det h\left(\bigwedge_i s_i,\bigwedge_i s_i\right)
     \right)(y)
 & = \sum_{\alpha =1}^n 
     -\log \det h\left(\bigwedge_i s_i,\bigwedge_i s_i\right)(x_\alpha) \\
 & = \sum_{\alpha =1}^n 
     -\log \det \left(
                h(s_i,s_j)(x_\alpha)
                \right)_{ij}.   
\end{align*}
Hence we have the desired equality of
currents by the above three equalities.

Thus we have showed the equality
(\ref{eqn:desired:thm:arith:Riemann:Roch:gen:finite:morphism}).
Since the right hand side belongs in fact to 
$\aBCycle^1(Y)$, 
the left hand side must also belong to 
$\aBCycle^1(Y)$, 
and thus we have the equality in 
$\aBCycle^1(Y)$. \QED

\section{Arithmetic Riemann-Roch for stable curves}
\subsection{Bismut-Bost formula}
\setcounter{Theorem}{0}
Let $X$ be a smooth algebraic variety over $\CC$,
$L$ a line bundle on $X$,
and $h$ a generalized metric of $L$ over $X$.
Let $s$ be a rational section of $L$. Then, by the definition of
the generalized metric $h$, $-\log h(s, s)$ gives rise to a current
$-[\log h(s,s)]$. Moreover, it is easy to see that a current
\[
dd^c(-[\log h(s,s)]) + \delta_{\operatorname{div}(s)}
\]
does not depend on the choice of $s$. Thus, we define $c_1(L, h)$ to be
\[
c_1(L, h) = dd^c(-[\log h(s,s)]) + \delta_{\operatorname{div}(s)}.
\]

Let $f : X \to Y$ be a proper morphism of smooth algebraic varieties
$\CC$ such that every fiber of $f$ is a reduced and connected curve
with only ordinary double singularities.
We set $\Sigma = \{ x \in X \mid \text{$f$ is not smooth at $x$.} \}$ and
$\Delta = f_*(\Sigma)$.
Let $|\Delta|$ be the support of $\Delta$.
We fix a Hermitian metric of $\omega_{X/Y}$. Then, in \cite{BBQm},
Bismut and Bost proved the following.

\begin{Theorem}
\label{thm:Quillen:metric:stable:curve}
Let $\overline{E} = (E, h)$ be a Hermitian vector bundle on $X$.
Then, the Quillen metric $h_Q^{\overline{E}}$ of $\det Rf_*(E)$ on 
$Y \setminus |\Delta|$ gives
rise to a generalized metric of $\det Rf_*(E)$ on $Y$.
Moreover,
\[
c_1 \left( \det Rf_*(E), h_Q^{\overline{E}} \right) = - f_* \left[
\todd({\overline{\omega}_{X/Y}}^{-1}) \chernch(\overline{E}) \right]^{(2,2)} 
- \frac{\rank E}{12} \delta_{\Delta}.
\]
\end{Theorem}

\subsection{Riemann-Roch for stable curves}
\setcounter{Theorem}{0}
In this subsection, we prove the arithmetic Riemann-Roch theorem 
for stable curves.
 
\begin{Theorem}
\label{thm:arith:Riemann:Roch:stable:curves}
Let $f : X \to Y$ be a projective morphism of regular arithmetic varieties
such that every fiber of $f_{\CC} : X(\CC) \to Y(\CC)$ is a reduced
and connected curve with only ordinary double singularities.
We fix a Hermitian metric of the dualizing sheaf $\omega_{X/Y}$.
Let $\overline{E} = (E, h)$ be a Hermitian vector bundle on $X$.
Then, 
\[
\acherncl_1 \left( \det Rf_*(E), h_Q^{\overline{E}} \right) -
\rank (E) \acherncl_1 \left( \det Rf_*(\OO_X), h_Q^{\overline{\OO}_X}
\right) \in 
\aBChow^1(Y)
\]
and
\begin{multline*}
\acherncl_1 \left( \det Rf_*(E), h_Q^{\overline{E}} \right) -
\rank (E) \acherncl_1 \left( \det Rf_*(\OO_X), h_Q^{\overline{\OO}_X}
\right) \\ = f_* \left( \frac{1}{2} \left(
\acherncl_1 (\overline{E})^2  -
\acherncl_1 (\overline{E}) \cdot \acherncl_1 (\overline{\omega}_{X/Y}) \right)
- \acherncl_2 (\overline{E})
\right)
\end{multline*}
in $\aBChow^1(Y)_{\QQ}$,
where $h_Q^{\overline{E}}$ and $h_Q^{\overline{\OO}_X}$ are
the Quillen metric of $\det Rf_*(E)$ and $\det Rf_*(\OO_X)$ respectively.
\end{Theorem}

\Proof We prove the theorem in two steps.

{\bf Step 1.}\quad
First, we assume that $f_{\QQ} : X_{\QQ} \to Y_{\QQ}$ is
smooth. In this case, by \cite{GSRR},
\[
\acherncl_1 \left( \det Rf_*(E), h_Q^{\overline{E}} \right) =
f_* \left( \achernch(E, h)\atodd(Tf, h_f) -
a(\chernch(E_{\CC})\todd(Tf_{\CC})R(Tf_{\CC}))
\right)^{(1)}.
\]
in $\aChow^1(Y)_{\QQ}$.
Since
\[
\achernch(\overline{E}) = \rank(E) + \acherncl_1(\overline{E}) +
\left( \frac{1}{2} \acherncl_1(\overline{E})^2 - 
\acherncl_2(\overline{E}) \right) + \text{(higher terms)}
\]
and
\[
\atodd(Tf, h_f) = 1 - \frac{1}{2} \acherncl_1(\overline{\omega}_{X/Y}) +
\atodd_2(Tf, h_f) + \text{(higher terms)},
\]
we have
\[
\left( \achernch(E, h)\atodd(Tf, h_f) \right)^{(2)} =
\frac{1}{2} \left(
\acherncl_1 (\overline{E})^2 -
\acherncl_1 (\overline{E}) \cdot \acherncl_1 (\overline{\omega}_{X/Y}) \right)
- \acherncl_2 (\overline{E}) + \rank(E) \atodd_2(Tf, h_f).
\]
On the other hand,
since the power series $R(x)$ has no constant term, the $(1,1)$ part of
\[
\chernch(E_{\CC})\todd(Tf_{\CC})R(Tf_{\CC})
\]
is $\rank(E) R_1(Tf_{\CC})$, where $R_1(Tf_{\CC})$ is the $(1,1)$ part of
$R(Tf_{\CC})$. Therefore, we obtain
\addtocounter{Claim}{1}
\begin{multline}
\label{eqn:1:thm:arith:Riemann:Roch:stable:curves}
\acherncl_1 \left( \det Rf_*(E), h_Q^{\overline{E}} \right) =
f_* \left( \frac{1}{2} \left(
\acherncl_1 (\overline{E})^2 -
\acherncl_1 (\overline{E}) \cdot \acherncl_1 (\overline{\omega}_{X/Y}) \right)
- \acherncl_2 (\overline{E})
\right) \\
+ \rank(E) f_* \left( \atodd_2(Tf, h_f) - a(R_1(Tf_{\CC})) \right).
\end{multline}
Applying (\ref{eqn:1:thm:arith:Riemann:Roch:stable:curves})
to the case $(E, h) = (\OO_X, h_{can})$,
we have
\addtocounter{Claim}{1}
\begin{equation}
\label{eqn:2:thm:arith:Riemann:Roch:stable:curves}
\acherncl_1 \left( \det Rf_*(\OO_X), h_Q^{\overline{\OO}_X} \right) =
f_* \left( \atodd_2(Tf, h_f) - a(R_1(Tf_{\CC})) \right).
\end{equation}
Thus, combining (\ref{eqn:1:thm:arith:Riemann:Roch:stable:curves})
and (\ref{eqn:2:thm:arith:Riemann:Roch:stable:curves}),
we have our formula in the case where
$f_{\QQ} : X_{\QQ} \to Y_{\QQ}$ is smooth.

{\bf Step 2.}\quad
Next, we consider the general case. 
The first assertion is a consequence of Theorem~\ref{thm:Quillen:metric:stable:curve} because
using Theorem~\ref{thm:Quillen:metric:stable:curve},
\begin{multline*}
\cherncl_1 \left( \det Rf_*(E), h_Q^{\overline{E}} \right) -
\rank (E) \cherncl_1 \left( \det Rf_*(\OO_X), h_Q^{\overline{\OO}_X}
\right) \\
 = 
- f_* \left[
\todd({\overline{\omega}_{X/Y}}^{-1}) \chernch(\overline{E})
 \right]^{(2,2)} 
+ \rank (E) f_* \left[
\todd({\overline{\omega}_{X/Y}}^{-1}) \chernch(\overline{\OO}_X)
 \right]^{(2,2)}
\end{multline*}
belongs to $\LocInt(\Omega_{Y(\CC)}^{1,1})$ 
by Proposition~\ref{prop:push:forward:B:pq}.
The second assertion is a consequence of the useful Lemma~\ref{lem:criterion:linear:equiv:B:cycle}. 
In fact, both sides of the second assertion are arithmetic $L^1$-cycles 
on Y by the first assertion and the Proposition~\ref{prop:push:forward:arith:cycle}:
If we take 
$\Delta=\{y \in Y_{\QQ} \mid \text{$f_{\QQ}$ is not smooth over $y$} \}$ 
and define $\overline {\Delta}$ to be the closure of $\Delta$ in $Y$, then 
the compliment $U=Y \setminus \overline {\Delta}$ 
contains no irreducible
components  of fibers of $Y \to \Spec(\ZZ)$ 
and $f_{\CC}$ is smooth over $U(\CC)$: 
The arithmetical linear equivalence of both sides restricted to $U$ is
a consequence of Step~1. 
Thus by Lemma~\ref{lem:criterion:linear:equiv:B:cycle}, we also have our
formula in the general case.
\QED

\section{Asymptotic behavior of analytic torsion}

\renewcommand{\theTheorem}{\arabic{section}.\arabic{Theorem}}
\renewcommand{\theequation}{\arabic{section}.\arabic{Theorem}}

Let $M$ be a compact K\"{a}hler manifold of dimension $d$, 
$\overline{E} = (E, h_E)$ a flat vector bundle of rank $r$ on $M$ 
with a flat metric $h_E$, and
$\overline{A} = (A, h_A)$ a Hermitian vector bundle on $M$. 
For $0 \leq q \leq d$, let $\Delta_{q,n}$ be the Laplacian on 
$A^{0,q}\left( \Sym^n(\overline{E}) \otimes \overline{A}
\right)$ and 
$\Delta_{q,n}'$ 
the restriction of $\Delta_{q,n}$ 
to $\Image \partial \oplus \Image \overline{\partial}$. Let 
$\sigma(\Delta_{q,n}') = \{ 0<\lambda_1 \leq \lambda_2 \leq \cdots \}$ 
be the sequence of eigenvalues of 
$\Delta_{q,n}'$.
Here we count each eigenvalue up to its multiplicity.
Then, the associated zeta function $\zeta_{q,n}(s)$ is given by
\[
\zeta_{q,n}(s) = \Tr\left[ (\Delta_{q,n}')^{-s} \right] =
\sum_{i=1}^{\infty} \lambda_i^{-s}.
\] 
It is well known that
$\zeta_{q,n}(s)$ converges absolutely for $\Re(s)>d$ 
and that it has a meromorphic continuation to the whole complex plane 
without pole at $s=0$. The analytic torsion 
$T \left( \Sym^n(\overline{E}) \otimes \overline{A} \right)$ 
is defined by
\[
T \left( \Sym^n(\overline{E}) \otimes \overline{A} \right)
= \sum_{q=0}^{d} (-1)^q q \zeta_{q,n}'(0).
\]
In the following we closely follow \cite[\S 2]{Vojta}. 

The Theta function associated with $\sigma(\Delta_{q,n}')$ is
defined by 
\[
\Theta_{q,n}(t) = \Tr \left[ \exp (-t \Delta_{q,n}') \right]
                = \sum_{i=1}^{\infty} e^{- \lambda_i t}.
\]
By taking Mellin transformation, we have, for $\Re(s)>d$,
\[
\zeta_{q,n}(s) = 
\frac{1}{\Gamma(s)} 
\int_0^{\infty} \Theta_{q,n}(t) t^s \frac{dt}{t}.
\]
We put 
\[
\tilde \zeta_{q,n}(s) = 
\frac{1}{\rank (\Sym^n(E))} n^{-d} \frac{1}{\Gamma(s)} 
\int_0^{\infty} \Theta_{q,n}\left(\frac{t}{n}\right) t^s \frac{dt}{t}.
\] 
Then we have 
\[
\frac{1}{\rank (\Sym^n(E))} n^{-d} \zeta_{q,n}(s) 
= n^{-s} \tilde \zeta_{q,n}(s)
\]
and thus
\addtocounter{Theorem}{1}
\begin{equation}
\label{eqn:1:relation:zeta:zero}
\frac{1}{\rank (\Sym^n(E))} n^{-d} \zeta_{q,n}'(0) 
= -(\log n) \tilde \zeta_{q,n}(0) + \zeta_{q,n}'(0)
\end{equation}

Bismut and Vasserot~\cite[(14),(19)]{BVAT} showed that 
$\Theta_{q,n}(t)$ has the following properties
(note that these 
parts of \cite{BVAT} do not depend on the assumption of
positivity of a line bundle, as indicated in Vojta
\cite[Proposition~2.7.3]{Vojta}):
\begin{enumerate}
\renewcommand{\labelenumi}{(\alph{enumi})}
 \item For every $k \in \NN$, $0 \leq q \leq d$ and $n \in \NN$, 
   there are real numbers $a_{q,n}^j\;(-d \leq j \leq k)$ 
   such that 
   \[
   \frac{1}{\rank (\Sym^n(E))} n^{-d} \Theta_{q,n}\left(\frac{t}{n}\right)
   = \sum_{j=-d}^{k} a_{q,n}^j t^j + o(t^k)
   \]
   as $t \downarrow 0$, 
   with $o(t^k)$ uniform with respect to $n \in \NN$. 
 \item For every $0
   \leq q \leq d$ and $j \geq -d$,  there are real numbers
   $a_{q}^j$  such that 
   \[
   a_{q,n}^j = a_{q}^j + O\left(\frac{1}{\sqrt n}\right)
   \]
   as $n \to  \infty$.
\end{enumerate}

Also by (b), 
we can replace the $o(t^k)$ in
(a) by $O(t^{k+1})$ and still have the uniformity statement. 
Thus we can write, for every $k \in \NN$,
\[
\frac{1}{\rank (\Sym^n(E))} n^{-d} \Theta_{q,n}\left(\frac{t}{n}\right)
   = \sum_{j=-d}^{k} a_{q,n}^k t^j + \rho_{q,n}^k(t)
\]
with $\rho_{q,n}^k(t) =  o(t^{k+1})$. Then
{\allowdisplaybreaks
\begin{align*}
\tilde \zeta_{q,n}(s)  & = 
\frac{1}{\rank (\Sym^n(E))} n^{-d} \frac{1}{\Gamma(s)} 
\int_1^{\infty} \Theta_{q,n}\left(\frac{t}{n}\right) t^s \frac{dt}{t} \\ 
 & \qquad + \frac{a_{q,n}^j}{\Gamma(s)} 
\int_0^1 t^{j+s-1} dt 
 +
\sum_{j=-d}^{k} \frac{1}{\Gamma(s)} 
\int_0^1 \rho_{q,n}^k(t) dt \\
 &  =
\frac{1}{\rank (\Sym^n(E))} n^{-d} \frac{1}{\Gamma(s)} 
\int_1^{\infty} \Theta_{q,n}\left(\frac{t}{n}\right) t^s \frac{dt}{t} \\
 & \qquad + \sum_{j=-d}^{k} \frac{a_{q,n}^j}{\Gamma(s) (j+s)} 
+ \frac{1}{\Gamma(s)} 
\int_0^1 \rho_{q,n}^k(t)t^s \frac{dt}{t}.
\end{align*}
}
In the last expression, the first integral is holomorphic 
for all $s \in \CC$, while the second integral is holomorphic for 
$\Re (s) > -k-1$; the middle term is a meromorphic function in the 
whole complex plane. 

Putting $k=0$ and $s=0$ in the above equation, we have
\addtocounter{Theorem}{1}
\begin{equation}
\label{eqn:3:zeta:zero}
\tilde \zeta_{q,n}(0) = a_{q,n}^0.
\end{equation}
Moreover, by differentiating the above equation  
when $k=0$, we have
\addtocounter{Theorem}{1}
\begin{multline}
\label{eqn:4:diff:zeta:zero}
\tilde \zeta_{q,n}'(0) = 
\frac{1}{\rank (\Sym^n(E))} n^{-d} 
\int_1^{\infty} \Theta_{q,n}\left(\frac{t}{n}\right) \frac{dt}{t} \\
+ \sum_{j=-d}^{-1} \frac{a_{q,n}^j}{j} - a_{q,n}^0 \Gamma'(1)
+ \frac{1}{\Gamma(s)} 
\int_0^1 \rho_{q,n}^0(t)  \frac{dt}{t}.
\end{multline}

We have now the following Proposition.

\begin{Proposition}
\label{prop:lower:bound:zeta:function}
There exists a constant c such that for all $n \in \NN$,
\[
\zeta_{q,n}'(0) \geq -c n^{d+r-1} \log n
\] 
\end{Proposition}

\Proof
By (\ref{eqn:1:relation:zeta:zero}),
(\ref{eqn:3:zeta:zero}) and (\ref{eqn:4:diff:zeta:zero}), we have
\begin{multline*}
\zeta_{q,n}'(0)
=
- \rank (\Sym^n(E)) n^d (\log n) a_{q,n}^0 \\
+
\rank (\Sym^n(E)) n^d  \left(
\frac{1}{\rank (\Sym^n(E))} n^{-d} 
\int_1^{\infty} \Theta_{q,n}\left(\frac{t}{n}\right) \frac{dt}{t} 
                        \right. \\
\left.
+ \sum_{j=-d}^{-1} \frac{a_{q,n}^j}{j} - a_{q,n}^0 \Gamma'(1)
+ \frac{1}{\Gamma(s)} 
\int_0^1 \rho_{q,n}^0(t)  \frac{dt}{t}
\right)
\end{multline*}
In the first term of the right hand side, 
$a_{q,n}^0$ is bounded with respect to $n$ by (b). 
In the second term of the right hand side, the first integral 
is non-negative; 
the sum of $a_{q,n}^j$'s 
is bounded with respect to $n$ by (b); 
the term $- a_{q,n}^0 \Gamma'(1)$ 
is also bounded with respect to  $n$ by (b); 
the second integral is also bounded with respect to  $n$, 
for $\rho_{q,n}^0(t)=O(t)$ uniformly with respect to  $n$.
Moreover,  
\[
\rank (\Sym^n(E)) = \binom{n+r-1}{r-1}=O(n^{r-1})
\]  
as $n \to \infty$. Thus, there is a constant $c$ 
such that for all $n \in \NN$, 
\[
\zeta_{q,n}'(0) \geq -c n^{d+r-1} \log n.
\]
\QED

\medskip
In the following sections, we only need the case of $d=1$, namely 
where $M$ is a compact Riemann surface. In this case, 
the above Proposition~\ref{prop:lower:bound:zeta:function} 
gives an asymptotic upper bound of analytic torsion.  

\begin{Corollary}
\label{cor:asymp:analytic:torsion}
Let $C$ be a compact Riemann surface,
$\overline{E} = (E, h_E)$ a flat vector bundle of rank $r$ on $C$ 
with a flat metric $h$, and
$\overline{A} = (A, h_A)$ a Hermitian vector bundle on $C$.
Then, there is a constant $c$ such that for all $n \in \NN$, 
\[
T \left( \Sym^n(\overline{E}) \otimes \overline{A} \right)
\leq c n^r \log n.
\]
\end{Corollary}

\Proof
Since $\dim C = 1$
\[
T \left( \Sym^n(\overline{E}) \otimes \overline{A} \right)
= -\zeta_{1,n}'(0).
\]
Now the corollary follows from 
Proposition~\ref{prop:lower:bound:zeta:function}.
\QED

\renewcommand{\theTheorem}{\arabic{section}.\arabic{subsection}.\arabic{Theorem}}
\renewcommand{\theequation}{\arabic{section}.\arabic{subsection}.\arabic{Theorem}.\arabic{Claim}}

\section{Formulae for arithmetic Chern classes}

\subsection{Arithmetic Chern classes of symmetric powers}
\label{subsec:formula:chern:sym:power}
\setcounter{Theorem}{0}
Let $M$ be a complex manifold and $(E, h)$
a Hermitian vector bundle on $M$.
Since $E^{\otimes n}$ has the natural Hermitian metric 
$h^{\otimes n}$,
we can define a Hermitian metric $\Sym^n(h)$ of $\Sym^n(E)$
to be the quotient metric of $E^{\otimes n}$ in terms of
the natural surjective homomorphism $E^{\otimes n} \to \Sym^n(E)$.
We denote $(\Sym^n(E), \Sym^n(h))$ by $\Sym^n(E, h)$.
If $x \in M$ and $\{ e_1, \ldots, e_r \}$ is an
orthonormal basis of $E_x$ with respect to $h_x$,
then it is easy to see that
\[
(\Sym^n(h))_x \left( e_1^{\alpha_1} \cdots e_r^{\alpha_r},
e_1^{\beta_1} \cdots e_r^{\beta_r} \right) =
\begin{cases}
{\displaystyle \frac{\alpha_1 ! \cdots \alpha_r !}{n !}} & 
\text{if $(\alpha_1, \ldots, \alpha_r) = (\beta_1, \ldots, \beta_r)$}, \\
0 & \text{otherwise}.
\end{cases}
\]
Then we have the following proposition.

\begin{Proposition}
\label{prop:chern:class:sym:power}
Let $X$ be an arithmetic variety and 
$\overline{E} = (E, h)$ a Hermitian vector
bundle of rank $r$ on $X$. Then, we have the following.
\begin{enumerate}
\renewcommand{\labelenumi}{(\arabic{enumi})}
\item
${\displaystyle \acherncl_1 \left( \Sym^n(\overline{E}) \right) =
\frac{n}{r} \binom{n+r-1}{r-1} \acherncl_1(\overline{E}) +
a \left( \sum_{\substack{\alpha_1 + \cdots + \alpha_r = n, \\ 
\alpha_1 \geq 0, \ldots, \alpha_r \geq 0}}
\log \left( \frac{n!}{\alpha_1 ! \cdots \alpha_r !} \right) \right)}$.

\item
If $X$ is regular, then
\begin{multline*}
\achernch_2 \left( \Sym^n(\overline{E}) \right) =
\binom{n+r}{r+1} \achernch_2(\overline{E}) + 
\frac{1}{2} \binom{n+r-1}{r+1} \acherncl_1(\overline{E})^2 \\
+ a \left( 
\frac{n}{r}
\sum_{\substack{\alpha_1 + \cdots + \alpha_r = n, \\ 
\alpha_1 \geq 0, \ldots, \alpha_r \geq 0}}
\log \left( \frac{n!}{\alpha_1 ! \cdots \alpha_r !} \right)
c_1(\overline{E}) \right).
\end{multline*}
\end{enumerate}
\end{Proposition}

\Proof
In \cite{SoVan}, C. Soul\'{e} gives similar formulae in implicit forms.
We follow his idea to calculate them.

\medskip
(1)
First of all, we fix notation.
We set
\[
S_{r, n} = \{
(\alpha_1, \ldots, \alpha_r) \in
(\ZZ_{+})^r \mid \alpha_1 + \cdots + \alpha_r = n \},
\]
where $\ZZ_{+} = \{ x \in \ZZ | x \geq 0 \}$.
For $I = (\alpha_1, \ldots, \alpha_r) \in S_{r, n}$ and
rational sections $s_1, \ldots, s_r$ of $E$,
we denote
$s_1^{\alpha_1} \cdots s_r^{\alpha_r}$ by $s^I$ and
$\alpha_1 ! \cdots \alpha_r !$ by $I !$.

Let $s_1, \ldots, s_r$ be independent
rational sections of $E$.
Then, $\{ s^I \}_{I \in S_{r,n}}$
forms independent rational sections of $\Sym^n(E)$.
First, let us see that
\addtocounter{Claim}{1}
\begin{equation}
\label{eqn:1:prop:chern:class:sym:power}
\zero \left( \bigwedge_{I \in S_{r,n}} s^I \right)
= \frac{n}{r} \binom{n+r-1}{r-1}
\zero ( s_1 \wedge \cdots \wedge s_r ).
\end{equation}
This is a local question.
So let $x \in X$ and $\{ \omega_1, \ldots, \omega_r \}$
be a local basis of $E$ around $x$.
We set $s_i = \sum_{j=1}^r a_{ij}\omega_j$.
Then, 
$s_1 \wedge \cdots \wedge s_r = \det(a_{ij})
\omega_1 \wedge \cdots \wedge \omega_r$.
Let $K$ be a rational function field of $X$.
Since the characteristic of $K$ is zero,
any $1$-dimensional representation of
$\operatorname{GL}_r(K)$ is a power of the determinant. 
Thus, there is an integer $N$ with
\[
\bigwedge_{I \in S_{r,n}} s^I =
\det(a_{ij})^N \bigwedge_{I \in S_{r,n}} \omega^I.
\]
Here, by an easy calculation, we can see that
\[
N = \frac{n}{r} \binom{n+r-1}{r-1}.
\]
Thus, we get
(\ref{eqn:1:prop:chern:class:sym:power}).

Next, let us see that
\addtocounter{Claim}{1}
\begin{multline}
\label{eqn:2:prop:chern:class:sym:power}
- \log \det \left(
\Sym^n(h)(s^I, s^J)
\right)_{I,J \in S_{r,n}} = \\
- \frac{n}{r} \binom{n+r-1}{r-1}
\log \det (h(s_i, s_j))_{i,j} +
\sum_{I \in S_{r,n}}
\log \left( \frac{n !}{I !} \right).
\end{multline}
Let $x \in X(\CC)$ and $\{ e_1, \ldots, e_r \}$ an
orthonormal basis of $E \otimes \kappa(x)$.
We set $s_i = \sum_{i=1}^r b_{ij}e_j$.
Moreover, we set
$s^I = \sum_{J \in S_{r,n}} b_{IJ}e^J$.
Then, in the same way as before,
$\det (b_{IJ}) = \det (b_{ij})^N$.
Further, since
\[
\Sym^n(h)(s^I, s^J) = \sum_{I', J' \in S_{r,n}}
b_{II'}\Sym^n(h)(e^{I'}, e^{J'}) \overline{b_{J'J}},
\]
we have
\begin{align*}
\det \left( \Sym^n(h)(s^I, s^J) \right)_{I,J \in S_{r,n}} & =
|\det(b_{IJ})|^2
\det \left( \Sym^n(h)(e^I, e^J) \right)_{I,J \in S_{r,n}} \\
& =
|\det(b_{ij})|^{2N}  \prod_{I \in S_{r,n}} \frac{I !}{n !}.
\end{align*}
Thus, we get (\ref{eqn:2:prop:chern:class:sym:power}).
Therefore, combining (\ref{eqn:1:prop:chern:class:sym:power})
and (\ref{eqn:2:prop:chern:class:sym:power}),
we obtain (1).

\medskip
(2)
First, we recall an elementary fact.
Let $\Phi \in \RR[X_1, \ldots, X_r]$ be a symmetric homogeneous
polynomial, and $M_r(\CC)$ the algebra of complex $r \times r$
matrices.
Then, there is a unique polynomial map
$\underline{\Phi} : M_r(\CC) \to \CC$ such that
$\underline{\Phi}$ is invariant under conjugation by $\GL_r(\CC)$
and its value on a diagonal matrix
$\diag(\lambda_1, \ldots, \lambda_r)$ is equal to
$\Phi(\lambda_1, \ldots, \lambda_r)$.

Let us consider the natural homomorphism
\[
\rho_{r,n} : \Aut_{\CC}(\CC^r) \to \Aut_{\CC}(\Sym^n(\CC^r))
\]
as complex Lie groups, which induces a homomorphism
\[
\gamma_{r, n} = d(\rho_{r,n})_{\operatorname{id}}
: \HEnd_{\CC}(\CC^r) \to \HEnd_{\CC}(\Sym^n(\CC^r))
\]
as complex Lie algebras.
Let $\{ e_1, \ldots, e_r \}$ be the standard basis of $\CC^r$.
Then, $\{ e_I \}_{I \in S_{r,n}}$ forms a basis of
$\Sym^n(\CC^r)$, where $e_I = e_{1}^{\alpha_1} \cdots e_r^{\alpha_r}$
for $I = (\alpha_1, \ldots, \alpha_r)$.
Let us consider the symmetric polynomial
\[
\chernch_2^{r,n} = \frac{1}{2}
\sum_{I \in S_{r,n}} X_I^2
\]
in $\RR[ X_I ]_{I \in S_{r,n}}$.
Then, by the previous remark, using the basis $\{ e_I \}_{I \in S_{r,n}}$,
we have a polynomial map
\[
\underline{\chernch_2^{r,n}} : \HEnd_{\CC}(\Sym^n(\CC^r))
\to \CC
\]
such that
$\underline{\chernch_2^{r,n}}$ is
invariant under conjugation by $\Aut_{\CC}(\Sym^n(\CC^r))$
and 
\[
\underline{\chernch_2^{r,n}}
\left( \diag(\lambda_I)_{I \in S_{r, n}} \right) =
\chernch_2^{r,n}(\ldots, \lambda_I, \ldots).
\]
Here we consider a polynomial map given by
\[
\begin{CD}
\theta_{r, n} :
\HEnd_{\CC}(\CC^r) @>{\gamma_{r,n}}>>
\HEnd_{\CC}(\Sym^n(\CC^r))
@>{\underline{\chernch_2^{r,n}}}>> \CC.
\end{CD}
\]
Since $\gamma_{r,n}(P A P^{-1}) = 
\rho_{r,n}(P)\gamma_{r,n}(A)\rho_{r,n}(P)^{-1}$ for
all $A \in \HEnd_{\CC}(\CC^r)$ and $P \in \Aut_{\CC}(\CC^r)$,
$\theta_{r,n}$ is invariant under conjugation by
$\Aut_{\CC}(\CC^r)$.
Let us calculate
\[
\theta_{r,n}(\diag(\lambda_1, \ldots, \lambda_r)).
\]
First of all,
\[
\gamma_{r,n}(\diag(\lambda_1, \ldots, \lambda_r)) =
\diag\left( \ldots, 
\left( \alpha_1 \lambda_1 + \cdots + \alpha_r \lambda_r \right),
\ldots \right)_{(\alpha_1, \ldots,  \alpha_r) \in S_{r,n}}.
\]
Thus,
\[
\theta_{r,n}(\diag(\lambda_1, \ldots, \lambda_r)) =
\frac{1}{2}
\sum_{(\alpha_1, \ldots, \alpha_r) \in S_{r,n}}
\left( \alpha_1 \lambda_1 + \cdots + \alpha_r \lambda_r \right)^2.
\]
On the other hand, by easy calculations,
we can see that
\[
\sum_{(\alpha_1, \ldots, \alpha_r) \in S_{r,n}}
\left( \alpha_1 \lambda_1 + \cdots + \alpha_r \lambda_r \right)^2 =
\binom{n+r}{r+1}\left( \lambda_1^2 + \cdots + \lambda_r^2 \right) +
\binom{n+r-1}{r+1}\left( \lambda_1 + \cdots + \lambda_r \right)^2.
\]
Therefore, we get
\[
\theta_{r,n}(\diag(\lambda_1, \ldots, \lambda_r)) =
\frac{1}{2} \binom{n+r}{r+1}\left( \lambda_1^2 + \cdots + \lambda_r^2 \right) +
\frac{1}{2} \binom{n+r-1}{r+1}\left( \lambda_1 + \cdots + \lambda_r \right)^2.
\]
Hence,
\addtocounter{Claim}{1}
\begin{equation}
\label{eqn:1:(2):prop:chern:class:sym:power}
\theta_{r,n} = \binom{n+r}{r+1} \underline{\chernch_2} +
\frac{1}{2} \binom{n+r-1}{r+1} \underline{(c_1)^2},
\end{equation}
where ${\displaystyle \chernch_2(X_1, \ldots, X_r) = 
\frac{1}{2}(X_1^2 + \cdots + X_r^2)}$ and
$c_1(X_1, \ldots, X_r) = X_1 + \cdots + X_r$.

Let $M$ be a complex manifold and $\overline{F} = (F, h_F)$ 
a Hermitian vector bundle of rank $r$ on $M$.
Let $K_{\overline{F}}$ be the curvature form of $\overline{F}$, and
$K_{\Sym^n(\overline{F})}$ the curvature form of $\Sym^n(\overline{F})$.
Then,
\[
K_{\Sym^n(\overline{F})} = 
\left( \gamma_{r, n} \otimes \operatorname{id}_{A^{1,1}(M)} \right)
(K_{\overline{F}}).
\]
Thus, by (\ref{eqn:1:(2):prop:chern:class:sym:power}),
\addtocounter{Claim}{1}
\begin{equation}
\label{eqn:2:(2):prop:chern:class:sym:power}
\chernch_2 \left( \Sym^n(F,h_F) \right) =
\binom{n+r}{r+1} \chernch_2(F,h_F) +
\frac{1}{2} \binom{n+r-1}{r+1} \cherncl_1(F,h_F)^2.
\end{equation}

Now let $\overline{E} = (E, h)$ be a Hermitian vector bundle on 
a regular arithmetic variety $X$.
Let $h'$ be another Hermitian metric of $E$. Then, using the definition
of Bott-Chern secondary characteristic classes and 
(\ref{eqn:2:(2):prop:chern:class:sym:power}),
\begin{multline*}
\achernch_2 \left( \Sym^n(E, h) \right) -
\achernch_2 \left( \Sym^n(E, h') \right) = \\
a \left( 
\binom{n+r}{r+1} \widetilde{\chernch_2}(E, h, h')  +
\frac{1}{2} \binom{n+r-1}{r+1} \widetilde{\cherncl_1^2}(E, h, h') 
\right).
\end{multline*}
Thus,
\[
\achernch_2 \left( \Sym^n(E, h) \right) -
\binom{n+r}{r+1} \achernch_2(E, h) -
\frac{1}{2} \binom{n+r-1}{r+1} \acherncl_1(E, h)^2
\]
does not depend on the choice of the metric $h$.
Therefore, in order to show (2),
by using splitting principle \cite[3.3.2]{GSCh},
we may assume that
\[
(E, h) = \overline{L}_1 \oplus \cdots \oplus \overline{L}_r,
\]
where $\overline{L}_i = (L_i, h_i)$'s 
are Hermitian line bundles.
Then,
\[
\Sym^n(\overline{E}) =
\bigoplus_{\substack{\alpha_1 + \cdots + \alpha_r = n, \\
\alpha_1 \geq 0, \ldots, \alpha_r \geq 0}}
\overline{L}_1^{\otimes \alpha_1} \otimes \cdots \otimes 
\overline{L}_r^{\otimes \alpha_r} \otimes 
\left( \OO_X, \frac{\alpha_1 ! \cdots \alpha_r !}{n!} h_{can} \right).
\]
Therefore, $\achernch_2 \left( \Sym^n(\overline{E}) \right)$
is equal to
\[
\sum_{\substack{\alpha_1 + \cdots + \alpha_r = n, \\
\alpha_1 \geq 0, \ldots, \alpha_r \geq 0}}
\left\{
\achernch_2 \left( \overline{L}_1^{\otimes \alpha_1} \otimes \cdots \otimes 
\overline{L}_r^{\otimes \alpha_r} \right)
-\log \left( \frac{\alpha_1 ! \cdots \alpha_r !}{n!} \right)
a\left( \cherncl_1 
\left( \overline{L}_1^{\otimes \alpha_1} \otimes \cdots \otimes 
\overline{L}_r^{\otimes \alpha_r} \right) \right)
\right\}.
\]
On the other hand, since
\begin{multline*}
\sum_{(\alpha_1, \cdots, \alpha_r) \in S_{r,n}}
\log \left( \frac{n!}{\alpha_1 ! \cdots \alpha_r !} \right)
\left( \alpha_1 X_1 + \cdots + \alpha_r X_r \right) \\
= \left(
\frac{n}{r}
\sum_{(\alpha_1, \cdots, \alpha_r) \in S_{r,n}}
\log \left( \frac{n!}{\alpha_1 ! \cdots \alpha_r !} \right)
\right)(X_1 + \cdots + X_r),
\end{multline*}
we have
{\allowdisplaybreaks
\begin{align*}
\achernch_2 \left( \Sym^n(\overline{E}) \right)
= & \binom{n+r}{r+1} \achernch_2(\overline{E}) +
\frac{1}{2} \binom{n+r-1}{r+1} \acherncl_1(\overline{E})^2  \\
& + \sum_{(\alpha_1, \cdots, \alpha_r) \in S_{r,n}}
\log \left( \frac{n!}{\alpha_1 ! \cdots \alpha_r !} \right)
a \left( \alpha_1 \cherncl_1(\overline{L_1}) + \cdots + 
\alpha_r\cherncl_1(\overline{L_r}) \right) \\
= & \binom{n+r}{r+1} \achernch_2(\overline{E}) +
\frac{1}{2} \binom{n+r-1}{r+1} \acherncl_1(\overline{E})^2 \\
& \qquad 
+ \left( 
\frac{n}{r}
\sum_{(\alpha_1, \cdots, \alpha_r) \in S_{r,n}}
\log \left( \frac{n!}{\alpha_1 ! \cdots \alpha_r !} \right)
\right)
a ( \cherncl_1(\overline{E})).
\end{align*}
}
Thus, we get (2).
\QED

\subsection{Arithmetic Chern classes of
$\overline{E} \otimes \overline{E}^{\vee}$}
\setcounter{Theorem}{0}
Here, let us consider arithmetic Chern classes of
$\overline{E} \otimes \overline{E}^{\vee}$.

\begin{Proposition}
\label{prop:ch2:end}
Let $X$ be a regular arithmetic variety and
$(E, h)$ a Hermitian vector bundle of rank $r$ on $X$.
Then,
\[
\achernch_2(E \otimes E^{\vee}, h \otimes h^{\vee})
= 2r \achernch_2(E, h) - \acherncl_1(E, h)^2
= (r-1) \acherncl_1(E, h)^2 - 2r \acherncl_2(E, h).
\]
\end{Proposition}

\Proof
Since $\achernch_i(E^{\vee}, h^{\vee}) = (-1)^i
\achernch_i(E, h)$ and
$\achernch(E \otimes E^{\vee}, h \otimes h^{\vee})
= \achernch(E, h) \cdot \achernch(E^{\vee}, h^{\vee})$,
we have
\begin{align*}
\achernch_2(E \otimes E^{\vee}, h \otimes h^{\vee}) & =
r \achernch_2(E, h) + \acherncl_1(E, h) \cdot 
\acherncl_1(E^{\vee}, h^{\vee}) + r \achernch_2(E^{\vee}, h^{\vee}) \\
& = 2r \achernch_2(E, h) - \acherncl_1(E, h)^2.
\end{align*}
The last assertion is derived from
the fact
\[
\achernch_2(E, h) = \frac{1}{2} \acherncl_1(E,h)^2 
- \acherncl_2(E, h).
\]
\QED

\section{The proof of the relative Bogomolov's inequality in the arithmetic case}
\renewcommand{\theTheorem}{\arabic{section}.\arabic{subsection}}

The purpose of this section is to give 
the proof of the following theorem.

\addtocounter{subsection}{1}
\begin{Theorem}[Relative Bogomolov's inequality in the arithmetic case]
\label{thm:relative:Bogomolov:inequality:arithmetic:case}
Let $f : X \to Y$ be a projective morphism of regular arithmetic varieties
such that every fiber of $f_{\CC} : X(\CC) \to Y(\CC)$ is a reduced
and connected curve with only ordinary double singularities.
Let $(E, h)$ be a Hermitian vector bundle of rank $r$ on $X$, 
and $y$ a closed point of $Y_{\QQ}$.
If $f$ is smooth over $y$ and $\rest{E}{X_{\bar{y}}}$ is
semi-stable, then
\[
\adis_{X/Y}(E, h) = f_* \left( 2r \acherncl_2(E, h) - (r-1)\acherncl_1(E, h)^2 \right)
\]
is weakly positive at $y$ with respect to any subsets
$S$ of $Y(\CC)$ with the following
properties: \rom{(1)} $S$ is finite, and
\rom{(2)} $f_{\CC}^{-1}(z)$ is smooth and $\rest{E_{\CC}}{f_{\CC}^{-1}(z)}$ is
poly-stable for all $z \in S$.
\end{Theorem}

\renewcommand{\theTheorem}{\arabic{section}.\arabic{subsection}.\arabic{Theorem}}

\subsection{Sketch of the proof of the relative Bogomolov's inequality}

The proof of the relative Bogomolov's inequality is very long, so that
for reader's convenience, 
we would like to give a rough sketch of
the proof of it.

\bigskip
{\bf Step 1.}\quad
Using the Donaldson's Lagrangian, we reduce to the case where
the Hermitian metric $h$ of $E$ along $f_{\CC}^{-1}(z)$ is
Einstein-Hermitian for each $z \in S$.

\medskip
{\bf Step 2.}\quad
We set
\[
\overline{F}_n = \Sym^n \left(
\End(\overline{E}) \otimes f^*(\overline{H}) \right) \otimes
\overline{A} \otimes f^*(\overline{H}),
\]
where $\overline{A}$ is a Hermitian line bundle
on $X$ and $\overline{H}$ is a Hermitian line bundle on $Y$.
Later we will specify these $\overline{A}$ and $\overline{H}$.
By virtue of the arithmetic Riemann-Roch for stable curves
(cf. Theorem~\ref{thm:arith:Riemann:Roch:stable:curves}) and
formulae of arithmetic Chern classes for symmetric powers
(cf. \S\ref{subsec:formula:chern:sym:power}), we can see that
\[
\frac{1}{(r^2 + 1)!} \adis_{X/Y}(\overline{E}) =
- \lim_{n \to \infty} \frac{\acherncl_1(\det Rf_*(F_n), h_n)}{n^{r^2+1}},
\]
where $h_n$ is a generalized metric of $\det Rf_*(F_n)$
such that $\acherncl_1(\det Rf_*(F_n), h_n)
\in \aBChow^1(Y)$ and
$h_n$ coincides with the Quillen metric $h_Q^{\overline{F}_n}$
at each $z \in S$.

\medskip
{\bf Step 3.}\quad
We assume that $A$ is very ample and $A \otimes \omega_{X/Y}^{-1}$
is ample.
We choose an arithmetic variety $B \subset X$
such that $B \in |A^{\otimes 2}|$, $B \to Y$ is \'{e}tale over $y$,
and $B(\CC) \to Y(\CC)$ is \'{e}tale over each $z \in S$.
(Exactly speaking, $B$ is not realized as an element of $|A^{\otimes 2}|$.
For simplicity, we assume it.)
We set $\overline{G}_n = \rest{\overline{F}_n}{B}$ and
$g = \rest{f}{B}$. 
Here we suppose that 
$g_*(\rest{\End(\overline{E})}{B}) \otimes \overline{H}$ and
$g_*(\rest{\overline{A}}{B}) \otimes \overline{H}$
are generated by small sections at $y$
with respect to $S$.

Applying the Riemann-Roch formula for generically finite morphisms
(cf. Theorem~\ref{thm:arith:Riemann:Roch:gen:finite:morphism}),
we can find a generalized metric $g_n$ of
$\det g_*(G_n)$ such that
$g_n$ is equal to the Quillen metric 
of $\overline{G}_n$ at each $z \in S$,
$\acherncl_1(\det g_*(G_n), g_n) \in \aBChow^1(Y)$, and
\[
\lim_{n \to \infty} \frac{\acherncl_1(\det g_*(G_n), g_n)}{n^{r^2+1}}
= 0.
\]

Let us consider the exact sequence:
\[
0 \to f_*(F_n) \to g_*(G_n) \to R^1f_*(F_n \otimes A^{\otimes -2})
\]
induced by $0 \to F_n \otimes A^{\otimes -2} \to F_n \to G_n \to 0$.
Let $Q_n$ be the image of
\[
g_*(G_n) \to R^1f_*(F_n \otimes A^{\otimes -2}).
\]
The natural $L^2$-metric of $g_*(G_n)$ 
around $z$ induces the quotient metric $\tilde{q}_n$
of $Q_n$ around $z$ for each $z \in S$.
Thus, we can find a $C^{\infty}$ metric $q_n$ of $\det Q_n$
such that $q_n$ is equal to $\det \tilde{q}_n$ at each $z \in S$.

Since
\[
\det Rf_*(F_n) = \det g_*(G_n) \otimes (\det Q_n)^{\otimes -1} \otimes 
\left( \det R^1f_*(F_n) \right)^{\otimes -1},
\]
we have the generalized metric $t_n$ of $\det R^1f_*(F_n)$ such that
\[
(\det Rf_*(F_n), h_n) = 
(\det g_*(G_n), g_n) \otimes (\det Q_n, q_n)^{\otimes -1} \otimes 
(\det R^1f_*(F_n), t_n)^{\otimes -1}.
\]

{\bf Step 4.}\quad
We set $
a_n = \max_{z \in S} \{ \log t_n(s_n,s_n)(z) \}$,
where $s_n$ is the canonical section of
$\det R^1f_*(F_n)$.
In this step, we will show that
$\acherncl_1(\det Q_n, q_n)$ is semi-ample at $y$ with respect to $S$ and
$a_n \leq O(n^{r^2}\log(n))$.
The semi-ampleness of $\acherncl_1(\det Q_n, q_n)$
at $y$ is derived from
Proposition~\ref{prop:find:small:section} and the fact
that $g_*(\rest{\End(\overline{E})}{B}) \otimes \overline{H}$ and
$g_*(\rest{\overline{A}}{B}) \otimes \overline{H}$
are generated by small sections at $y$ with respect to $S$.
The estimation of $a_n$ involves asymptotic behavior of analytic torsion
(cf. Corollary~\ref{cor:asymp:analytic:torsion}) 
and a comparison of sup-norm with $L^2$-norm
(cf. Lemma~\ref{lem:comparison:sup:L2}).

{\bf Step 5.}\quad
Thus, using the last equation in Step 3,
we can get a decomposition
\[
-\frac{\acherncl_1(\det Rf_*(F_n), h_n)}{n^{r^2+1}} =
\alpha_n + \beta_n
\]
such that $\alpha_n$ is semi-ample at $y$ with respect to $S$ and
$\lim_{n \to \infty} \beta_n = 0$.

\subsection{Preliminaries}
\setcounter{Theorem}{0}
First of all, we will prepare three lemmas for the proof
of the relative Bogomolov's inequality.

\begin{Lemma}
\label{lem:comparison:sup:L2}
Let $M$ be a $d$-dimensional compact K\"{a}hler manifold,
$\overline{E} = (E, h)$ a flat Hermitian vector bundle
of rank $r$ on $M$, and $\overline{V} = (V, k)$
a Hermitian line bundle.
Then, there is a constant $c$ such that,
for any $n > 0$ and any $s \in H^0(M, \Sym^n(E) \otimes V)$,
\[
  \Vert s \Vert_{\sup} \leq c n^{d+r-1} \Vert s \Vert_{L^2}.
\]
\end{Lemma}

\Proof
Let $f : P = \Proj \left(\bigoplus_{i \geq 0} \Sym^i(E) \right) \to M$ 
be the projective bundle of $E$, and
$L = \OO_P(1)$ the tautological line bundle of $E$ on $P$.
Let $h_L$ be the quotient metric of $L$ induced by
the surjective homomorphism $f^*(E) \to L$ and 
the Hermitian metric $f^*(h)$ of $f^*(E)$.
Let $\Omega_M$ be a K\"{a}hler form of $M$.
Since $\overline{E}$ is flat,
$c_1(L, h_L)$ is positive semi-definite of rank $r-1$.
Thus, $f^*(\Omega_M) + c_1(L, h_L)$ gives rise to
a fundamental $2$-form $\Omega_P$ on $P$.
Moreover, by virtue of the flatness of $\overline{E}$,
we have $c_1(L, h_L)^r = 0$.
Thus,
\[
\Omega_P^{d+r-1} = \binom{d+r-1}{d} f^*(\Omega_M^d) \wedge c_1(L, h_L)^{r-1}.
\]

By \cite[Lemma~30]{GSRR}, there is a constant $c$ such that
\[
\Vert s' \Vert_{\sup} \leq c n^{d+r-1} \Vert s' \Vert_{L^2}
\]
for any $n > 0$ and any $s' \in H^0(P, L^{\otimes n} \otimes f^*(V))$,
where ${\displaystyle \Vert s' \Vert_{L^2} = \int_{P} |s'|^2 \Omega_P^{d+r-1}}$.
We denote a homomorphism
\[
f^*(\Sym^n(E)) \otimes f^*(V) \to L^{\otimes n} \otimes f^*(V)
\]
by $\alpha_n$.
As in the proof of \cite[(44)]{GSRR}, we can see that,
for any $s \in H^0(M, \Sym^n(E) \otimes V)$,
\[
|s|^2 = \binom{n+r-1}{r-1}
\int_{P \to M} | \alpha_n(s)  |^2 c_1(L, h_L)^{r-1}.
\]
Thus,
\[
|s|^2 \leq \binom{n+r-1}{r-1}
\int_{P \to M} \Vert \alpha(s) \Vert_{\sup}^2 c_1(L, h_L)^{r-1}
= \binom{n+r-1}{r-1} \Vert \alpha(s) \Vert_{\sup}^2.
\]
Therefore, we get 
\[
\Vert s \Vert_{\sup}^2 \leq \binom{n+r-1}{r-1} \Vert \alpha_n(s) \Vert_{\sup}^2
\]
for all $s \in H^0(M, \Sym^n(E) \otimes V)$.

On the other hand,
{\allowdisplaybreaks
\begin{align*}
\Vert \alpha_n(s) \Vert_{L^2}^2 & =
\int_{P} |\alpha_n(s)|^2 \Omega_P^{r} \\
& =
\binom{d+r-1}{d} \int_M \int_{P \to M} |\alpha_n(s)|^2 f^*(\Omega_M^d) \wedge c_1(L, h_L)^{r-1} \\
& =
\binom{d+r-1}{d} \int_M \Omega_M^d \int_{P \to M} |\alpha_n(s)|^2 c_1(L, h_L)^{r-1} \\
& =
\binom{d+r-1}{d} \binom{n+r-1}{r-1}^{-1} \int_M |s|^2 \Omega_M^d  \\
& =
\binom{d+r-1}{d} \binom{n+r-1}{r-1}^{-1} \Vert s \Vert_{L^2}^2.
\end{align*}
}
Therefore,
\begin{align*}
\Vert s \Vert_{\sup}^2 
& \leq \binom{n+r-1}{r-1} \Vert \alpha_n(s) \Vert_{\sup}^2 \\
& \leq \binom{n+r-1}{r-1} c^2 n^{2(d+r-1)} \Vert \alpha_n(s) \Vert_{L^2}^2 \\
& = \binom{d+r-1}{d} c^2 n^{2(d+r-1)} \Vert s \Vert_{L^2}^2.
\end{align*}
Thus, we get our lemma.
\QED

Here we recall Einstein-Hermitian metrics of vector bundles.
Let $M$ be a $d$-dimensional 
compact K\"{a}hler manifold with a K\"{a}hler form $\Omega_M$,
and $E$ a vector bundle on $M$.
We say $E$ is {\em stable} (resp. {\em semistable}) 
{\em with respect to $\Omega_M$}
if, for any subsheaves $F$ of $E$
with $0 \subsetneq F \subsetneq E$,
\[
\frac{1}{\rank F} \int_M c_1(F) \wedge \Omega_M^{d-1} 
< 
\frac{1}{\rank E} \int_M c_1(E) \wedge \Omega_M^{d-1}.
\]
\[
\left( \text{resp.}\quad
\frac{1}{\rank F} \int_M c_1(F) \wedge \Omega_M^{d-1} 
\leq
\frac{1}{\rank E} \int_M c_1(E) \wedge \Omega_M^{d-1}. \right)
\]
Moreover, $E$ is said to be {\em poly-stable with respect to $\Omega_M$}
if $E$ is semistable with respect to $\Omega_M$ and
$E$ has a decomposition $E = E_1 \oplus \cdots \oplus E_s$ of vector bundles
such that each $E_i$ is stable with respect to $\Omega_M$.
Let $h$ be a Hermitian metric of $E$.
We say $h$ is {\em Einstein-Hermitian with respect to $\Omega_M$} 
if there is a constant $\rho$
such that
$K(E, h) \wedge \Omega_M^{d-1} = \rho \Omega_M^d \otimes \operatorname{id}_E$,
where $K(E, h)$ is the curvature form given by $(E, h)$ and
$\operatorname{id}_E$ is the identity map in $\Hom(E, E)$.
The Kobayashi-Hitchin correspondence tells us that
$E$ has an Einstein-Hermitian metric with respect to $\Omega_M$ if and only if
$E$ is poly-stable with respect to $\Omega_M$.

\begin{Lemma}
\label{lem:sum:EH:metric}
Let $M$ be a compact K\"{a}hler manifold with
a K\"{a}hler form $\Omega_M$, and $E$ a poly-stable vector
bundle with respect to $\Omega_M$ on $M$.
If $h$ and $h'$ are Einstein-Hermitian metrics of $E$
with respect to $\Omega_M$, then so is $h+h'$.
\end{Lemma}

\Proof
Let $E = E_1 \oplus \cdots \oplus E_s$ be a decomposition into
stable vector bundles.
If we set $h_i = \rest{h}{E_i}$ and
$h'_i = \rest{h'}{E_i}$ for each $i$, then $h_i$ and $h'_i$
are Einstein-Hermitian metrics of $E_i$ and we have the following
orthogonal decompositions:
\[
(E, h) = \bigoplus_{i=1}^s (E_i, h_i)
\quad\text{and}\quad
(E, h') = \bigoplus_{i=1}^s (E_i, h'_i)
\]
(cf. \cite[Chater~IV, \S~3]{Ko}).
Thus, we may assume that $E$ is stable.
In this case, by virtue of the uniqueness of Einstein-Hermitian metric,
there is a positive constant $c$ with $h' = ch$.
Thus, $h + h' = (1+c)h$. Hence $h+h'$ is Einstein-Hermitian.
\QED

\begin{Lemma}
\label{lem:polystable:complex:conjugation}
Let $C$ be a compact Riemann surface.
Considering $C$ as a projective variety over $\CC$,
let $\overline{C} = C \otimes_{\CC} \CC$ be the tensor
product via the complex conjugation.
Let $E$ be a vector bundle on $C$, and
$\overline{E} = E \otimes_{\CC} \CC$ on $\overline{C}$.
Then, $E$ is poly-stable on $C$ if and only if $\overline{E}$
is poly-stable on $\overline{C}$.
\end{Lemma}

\Proof
This is an easy consequence of the fact that
if $F$ is a vector bundle on $C$, then
$\deg(F) = \deg(\overline{F})$.
\QED

\subsection{Complete proof of the relative Bogomolov's inequality}
\setcounter{Theorem}{0}
\renewcommand{\theClaim}{\arabic{section}.\arabic{subsection}.\arabic{Claim}}
\renewcommand{\theequation}{\arabic{section}.\arabic{subsection}.\arabic{Claim}}

Let us start the complete proof of the relative Bogomolov's inequality.

Considering $S \cup F_{\infty}(S)$ instead of $S$,
we may assume that $F_{\infty}(S) = S$ 
by virtue of Lemma~\ref{lem:polystable:complex:conjugation}.
For each $z \in S$, let $\Omega_z$ be the K\"{a}hler form induced
by the metric of $\overline{\omega}_{X/Y}$ along $f_{\CC}^{-1}(z)$.
Since $\rest{E_{\CC}}{f_{\CC}^{-1}(z)}$ is poly-stable for all $z \in S$,
there is a $C^{\infty}$ Hermitian metric $h'$ of $E_{\CC}$ such that
$\rest{h'}{f_{\CC}^{-1}(z)}$ is Einstein-Hermitian with respect to $\Omega_z$
for all $z \in S$.
It is easy to see that $\rest{\overline{F_{\infty}^*(h')}}{f_{\CC}^{-1}(z)}$
is Einstein-Hermitian with respect to $\Omega_z$
for all $z \in S$. Thus,
if $h'$ is not invariant under $F_{\infty}$,
then, considering $h' + \overline{F_{\infty}^*(h')}$,
we may assume that $h'$ is invariant under $F_{\infty}$.
For, by Lemma~\ref{lem:sum:EH:metric},
$h' + \overline{F_{\infty}^*(h')}$
is Einstein-Hermitian with respect to $\Omega_z$ on $f_{\CC}^{-1}(z)$
for each $z \in S$.

Here we claim:

\begin{Claim}
\label{claim:assume:Einstein:Hermitian}
There is a $\gamma \in \LocInt(Y(\CC))$ such that
$a(\gamma) \in \aBChow^1(Y;S)$ and
$\gamma(z) \geq 0$ for each $z \in S$, and
\[
\adis_{X/Y}(E, h) = \adis_{X/Y}(E, h') + a(\gamma).
\]
\end{Claim}

\Proof
We set ${\displaystyle \phi = \sqrt[r]{\det(h')/\det(h)}}$.
Then, it is easy to see that
$\adis_{X/Y}(E, \phi h) = \adis_{X/Y}(E, h)$.
Thus, we may assume that $\det(h) = \det(h')$.
Then, we have
\[
\adis_{X/Y}(E, h) - \adis_{X/Y}(E, h') = 
a \left( - f_* (2r \widetilde{\chernch_2}(E, h, h')) \right).
\]
Hence if we set $\gamma = - f_* (2r \widetilde{\chernch_2}(E, h, h'))$,
then $a(\gamma) \in \aBChow^1(Y;S)$.
On the other hand, by \cite[(ii) of Corollary~1.30]{BGSAT},
$- f_*(\widetilde{\chernch_2}(E, h, h'))(z)$ is
nothing more than Donaldson's Lagrangian (for details,
see \cite[\S6]{MoBG}).
Thus, we get $\gamma(z) \geq 0$ for each $z \in S$.
\QED

\medskip
By the above claim, we may assume that
$\rest{h}{f_{\CC}^{-1}(z)}$ is Einstein-Hermitian for each $z \in S$.
Let $\overline{A} = (A, h_A)$ be a Hermitian line bundle
on $X$ such that $A$ is very ample, and 
$A \otimes \omega_{X/Y}^{\otimes -1}$ is ample.
If we take a general member $M'$ of $|A_{\QQ}^{\otimes 2}|$,
then, by Bertini's theorem (cf. \cite[Theorem~6.10]{JB}), 
$M'$ is smooth over $\QQ$,
and $M' \to Y_{\QQ}$ is \'{e}tale over $y$. 
Note that if $Z$ is an algebraic set of $\PP^N_{\CC}$,
$U$ is a non-empty Zariski open set of $\PP^N_{\QQ}$, and
$U(\QQ) \subseteq Z(\CC)$, then $Z = \PP^N_{\CC}$.
Hence, we may assume that $M'(\CC) \to Y(\CC)$ is
\'{e}tale over $z$ for all $z \in S$.
Let $M' = M'_1 + \cdots + M'_{l_1} + M'_{l_1 + 1} + \cdots + M'_{l_2}$
be the decomposition of $M'$ into irreducible components
(actually, the decomposition into connected components 
because $M'$ is smooth over $\QQ$) such that $f_{\QQ}(M'_i) = Y_{\QQ}$
for $1 \leq i \leq l_1$ and $f_{\QQ}(M'_j) \subsetneq Y_{\QQ}$
for $l_1 + 1 \leq j \leq l_2$.
Let $M_i$ ($i=1, \ldots, l_1$) be the closure of $M'_i$ in $X$.
We set $M = M_1 + \cdots + M_{l_1}$ and
$B = M_1 \coprod \cdots \coprod M_{l_1}$ (disjoint union).
Then, there is a line bundle $L$ on $X$ with $M \in | A^{\otimes 2} \otimes L |$.
Note that $\rest{L}{X_y} \simeq \OO_{X_y}$ and $
\rest{L_{\CC}}{f_{\CC}^{-1}(z)} \simeq \OO_{f_{\CC}^{-1}(z)}$ for all $z \in S$ 
because
$y \not\in \bigcup_{j=l_1+1}^{l_2} f_{\QQ}(M'_j)$ and
$z \not\in \bigcup_{j=l_1+1}^{l_2} f_{\CC}(M'_j(\CC))$.
We denote the morphism $B \to M \to X$ by $\iota$, and
the morphism 
$B \overset{\iota}{\longrightarrow} X \overset{f}{\longrightarrow} Y$
by $g$.
We remark that the morphism $B \to M$ is an isomorphism over $\QQ$.
Further, we set
\[
\overline{F} = \End(E, h) = (E \otimes E^{\vee}, h \otimes h^{\vee}).
\]
Then, $h \otimes h^{\vee}$ is a flat metric along $f_{\CC}^{-1}(z)$
for each $z \in S$
because $h \otimes h^{\vee}$ is Einstein-Hermitian and
$\deg \left( E \otimes E^{\vee} \right) = 0$ along $f_{\CC}^{-1}(z)$.
We choose a Hermitian line bundle $\overline{H} = (H, h_H)$ on $Y$
such that $g_*(\iota^*(A)) \otimes H$ and
$g_*(\iota^*(F)) \otimes H$ are generated by small sections at $y$
with respect to $S$.
Moreover, we set
\[
\overline{F}_n = \Sym^n \left(
\overline{F} \otimes f^*(\overline{H}) \right) \otimes
\overline{A} \otimes f^*(\overline{H})
= \left( \Sym^n \left( F \otimes f^*(H) \right) \otimes A \otimes f^*(H),
k_n \right).
\]

\begin{Claim}
\label{claim:terms:right:R:R:formla}
There are $Z_0, \ldots, Z_{r^2} \in \aBChow^1(Y;S)_{\QQ}$ and
$\beta \in \LocInt(Y(\CC))$ such that
$a(\beta) \in \aBChow^1(Y;S)$,
and
\[
f_* \left( \achernch_2(\overline{F}_n) - 
\frac{1}{2}
\acherncl_1(\overline{F}_n) \cdot \acherncl_1(\overline{\omega}_{X/Y})
\right) =
\frac{n^{r^2 + 1}}{(r^2 + 1) !}
f_* (\achernch_2(\overline{F})) +
\sum_{i=0}^{r^2} Z_i n^i + a(b_n \beta),
\]
where
${\displaystyle
b_n = \sum_{\substack{\alpha_1 + \cdots + \alpha_{r^2} = n, \\ 
\alpha_1 \geq 0, \ldots, \alpha_{r^2} \geq 0}}
\log \left( \frac{n!}{\alpha_1 ! \cdots \alpha_{r^2} !} \right)}$.
\end{Claim}

\Proof
Since $\Sym^n(\overline{F} \otimes f^*(\overline{H})) \otimes \overline{A}
\otimes f^*(\overline{H})$
is isometric to
$\Sym^n(\overline{F}) \otimes f^*(\overline{H})^{\otimes (n+1)} \otimes
\overline{A}$,
\begin{multline*}
\achernch_2(\overline{F}_n) =
\achernch_2(\Sym^n(\overline{F})) +
\acherncl_1(\Sym^n(\overline{F})) \cdot 
\acherncl_1(f^*(\overline{H})^{\otimes (n+1)} \otimes \overline{A}) \\
+\binom{n+r^2-1}{r^2-1} 
\achernch_2(f^*(\overline{H})^{\otimes (n+1)} \otimes \overline{A}).
\end{multline*}
Here since $\det(\overline{F}) = \overline{\OO}_X$,
by Proposition~\ref{prop:chern:class:sym:power},
\[
\acherncl_1(\Sym^n(\overline{F})) = a(b_n)
\quad\text{and}\quad
\achernch_2(\Sym^n(\overline{F})) = \binom{n+r^2}{r^2+1}
\achernch_2(\overline{F}).
\]
Thus, by Proposition~\ref{prop:projection:formula:line:bundle},
\begin{align*}
f_* \left( \acherncl_1(\Sym^n(\overline{F})) \cdot 
\acherncl_1(f^*(\overline{H})^{\otimes (n+1)} \otimes \overline{A}) \right)
& = f_* \left(
b_n a \left((n+1) f^*(c_1(\overline{H})) + c_1(\overline{A}) \right) \right) \\
& = a \left( b_n f_*(c_1(\overline{A})) \right).
\end{align*}
On the other hand, using the projection formula 
(cf. Proposition~\ref{prop:projection:formula:line:bundle}),
\begin{align*}
f_* \left( 
\achernch_2(f^*(\overline{H})^{\otimes (n+1)} \otimes \overline{A}) \right)
& = \frac{1}{2} f_* \left[ \left(
(n+1) \acherncl_1(f^*(\overline{H})) + \acherncl_1(\overline{A}) 
\right)^2 \right] \\
& = \frac{1}{2} f_* \left[
(n+1)^2 \acherncl_1(f^*(\overline{H}))^2 + 
2(n+1) \acherncl_1(f^*(\overline{H})) \cdot \acherncl_1(\overline{A}) +
\acherncl_1(\overline{A})^2 \right] \\
& = (n+1) \deg_f(A) \acherncl_1(\overline{H}) + 
\frac{1}{2} f_* \left( \acherncl_1(\overline{A})^2 \right),
\end{align*}
where $\deg_f(A)$ is the degree of $A$ on the generic fiber of $f$.
Therefore, we have
\begin{multline*}
f_* \achernch_2(\overline{F}_n) = 
\binom{n+r^2}{r^2+1} f_* \achernch_2(\overline{F}) + \\
\binom{n+r^2-1}{r^2-1} \left(
(n+1) \deg_f(A) \acherncl_1(\overline{H}) + 
\frac{1}{2} f_* \left( \acherncl_1(\overline{A})^2 \right) 
\right) +
a \left( b_n f_*(c_1(\overline{A})) \right).
\end{multline*}
Thus, there are $Z'_0, \ldots, Z'_{r^2} \in \aChow^1(Y;S)_{\QQ}$
such that
\addtocounter{Claim}{1}
\begin{equation}
\label{eqn:ch2:Fn}
f_* \achernch_2(\overline{F}_n) = 
\frac{n^{r^2+1}}{(r^2+1)!} f_* \achernch_2(\overline{F}) +
\sum_{i=0}^{r^2} Z'_i n^i +
a \left( b_n f_*(c_1(\overline{A})) \right).
\end{equation}
Further, since 
$\acherncl_1(\overline{F}_n) \cdot \acherncl_1(\overline{\omega}_{X/Y})$
is equal to
\[
\left(
\acherncl_1(\Sym^n(\overline{F})) + \binom{n+r^2-1}{r^2-1}
((n+1) \acherncl_1(f^*(\overline{H})) + \acherncl_1(\overline{A}) ) 
\right) \cdot
\acherncl_1(\overline{\omega}_{X/Y}),
\]
we have
\begin{align*}
f_* \left(
\acherncl_1(\overline{F}_n) \cdot \acherncl_1(\overline{\omega}_{X/Y})
\right) 
& = a \left( b_n f_*(c_1(\overline{\omega}_{X/Y})) \right) + \\
& \qquad
\binom{n+r^2-1}{r^2-1} \left(
(n+1) (2g-2) \acherncl_1(\overline{H}) + 
f_* \left( \acherncl_1(\overline{A}) \cdot 
\acherncl_1(\overline{\omega}_{X/Y}) \right)
\right).
\end{align*}
Thus, there are $Z''_0, \ldots, Z''_{r^2} \in \aChow^1(Y;S)_{\QQ}$
such that
\addtocounter{Claim}{1}
\begin{equation}
\label{eqn:c1:Fn:c1:w}
f_* \left(
\acherncl_1(\overline{F}_n) \cdot \acherncl_1(\overline{\omega}_{X/Y})
\right) =
\sum_{i=0}^{r^2} Z''_i n^i + 
a \left( b_n f_*(c_1(\overline{\omega}_{X/Y})) \right).
\end{equation}
Thus, combining (\ref{eqn:ch2:Fn}) and
(\ref{eqn:c1:Fn:c1:w}), we get our claim.
\QED

Let $h_{X/Y}$ be a $C^{\infty}$ Hermitian metric of 
$\det Rf_* \OO_X$ over $Y(\CC)$
such that $h_{X/Y}$ is invariant under $F_{\infty}$.
Then, since the Quillen metric $h^{\overline{\OO_X}}_Q$
of $\det Rf_* \OO_X$ is a generalized metric,
there is a real valued $\phi \in \LocInt(Y(\CC))$ such that
$h^{\overline{\OO_X}}_Q = e^{\phi} h_{X/Y}$ and
$F_{\infty}^*(\phi) = \phi \ (\alev)$.
Adding a suitable real valued $C^{\infty}$ function $\phi'$ with
$F_{\infty}^*(\phi') = \phi'$ to $\phi$
(replace $h_{X/Y}$ by $e^{-\phi'}h_{X/Y}$ accordingly),
we may assume that $\phi(z) = 0$ for
all $z \in S$.
Here, we set
${\displaystyle h_n = \exp \left( - \binom{n+r^2-1}{r^2-1} \phi \right) 
h^{\overline{F}_n}_Q}$.
Then, $h_n$ is a generalized metric of $\det Rf_* F_n$ with
$F_{\infty}^*(h_n) = \overline{h}_n \ (\alev)$.
Moreover,
\begin{multline*}
\acherncl_1 \left( \det Rf_* F_n, h_n \right) -
\binom{n+r^2-1}{r^2-1} \acherncl_1 \left( \det Rf_* \OO_X, h_{X/Y} \right) \\
= \acherncl_1 \left( \det Rf_* F_n, h^{\overline{F}_n}_Q \right) -
\binom{n+r^2-1}{r^2-1} \acherncl_1 
\left( \det Rf_* \OO_X, h^{\overline{\OO_X}}_Q \right).
\end{multline*}
Here, since
\[
\acherncl_1 \left( \det Rf_* F_n, h^{\overline{F}_n}_Q \right) -
\binom{n+r^2-1}{r^2-1} \acherncl_1 
\left( \det Rf_* \OO_X, h^{\overline{\OO_X}}_Q \right)
\in \aBChow^1(Y;S)_{\QQ}
\]
by Theorem~\ref{thm:arith:Riemann:Roch:stable:curves}
and $\acherncl_1(\det Rf_* \OO_X, h_{X/Y}) \in \aChow^1(Y;S)$,
we have
\[
\acherncl_1 \left( \det Rf_* F_n, h_n \right) \in \aBChow^1(Y;S)_{\QQ}.
\]
Further, by the arithmetic Riemann-Roch theorem for
stable curves (cf. Theorem~\ref{thm:arith:Riemann:Roch:stable:curves}),
\begin{multline*}
\acherncl_1 \left( \det Rf_*(F_n), h_n \right) -
\binom{n+r^2-1}{r^2-1} \acherncl_1 
\left( \det Rf_*(\OO_X), h_{X/Y} \right) \\
= f_* \left( \achernch_2(\overline{F}_n) -
\frac{1}{2}
\acherncl_1 (\overline{F}_n) \cdot \acherncl_1 (\overline{\omega}_{X/Y}) 
\right).
\end{multline*}
Therefore,
by Claim~\ref{claim:terms:right:R:R:formla},
there are $W_0, \ldots, W_{r^2} \in \aBChow^1(Y;S)_{\QQ}$ and
$\beta \in \LocInt(Y(\CC))$ such that
$a(\beta) \in \aBChow^1(Y;S)$, and
\addtocounter{Claim}{1}
\begin{equation}
\label{eqn:1:proof:arith:BG:inq}
\acherncl_1 \left( \det Rf_*(F_n), h_n \right) =
\frac{n^{r^2 + 1}}{(r^2 + 1) !}
f_* (\achernch_2(\overline{F})) +
\sum_{i=0}^{r^2} W_i n^i + a(b_n \beta).
\end{equation}

\begin{Claim}
\label{claim:dis:lim:c1:Fn}
${\displaystyle
\frac{1}{(r^2+1)!} \adis_{X/Y}(\overline{E}) =
-\lim_{n \to \infty} 
\frac{\acherncl_1 \left( \det Rf_*(F_n), h_n \right)}{n^{r^2+1}}
}$ in $\aBChow^1(Y;S)_{\QQ}$.
\end{Claim}

\Proof
By virtue of Proposition~\ref{prop:ch2:end},
$f_*(\achernch_2(\overline{F})) = -\adis_{X/Y}(\overline{E})$.
Thus, by (\ref{eqn:1:proof:arith:BG:inq}),
it is sufficient to show that $0 \leq b_n \leq O(n^{r^2})$.

It is well known that
\[
\frac{\log(\theta_1) + \cdots + \log(\theta_N)}{N} \leq 
\log \left( \frac{\theta_1 + \cdots + \theta_N}{N} \right)
\]
for positive numbers $\theta_1, \ldots, \theta_N$.
Thus, noting
${\displaystyle 
\sum_{\substack{\alpha_1 + \cdots + \alpha_{r^2} = n, \\ 
\alpha_1 \geq 0, \ldots, \alpha_{r^2} \geq 0}}
\frac{n!}{\alpha_1 ! \cdots \alpha_{r^2} !} = (r^2)^n}$,
we have
\[
0 \leq
\sum_{\substack{\alpha_1 + \cdots + \alpha_{r^2} = n, \\ 
\alpha_1 \geq 0, \ldots, \alpha_{r^2} \geq 0}}
\log \left( \frac{n!}{\alpha_1 ! \cdots \alpha_{r^2} !} \right) 
\leq
\binom{n+r^2-1}{r^2-1} \log \left(
\frac{ (r^2)^n }{\binom{n+r^2-1}{r^2-1}}
\right) 
\leq O(n^{r^2}).
\]
\QED

We set $\overline{G}_n = \iota^*(\overline{F}_n)$.
Then, by Theorem~\ref{thm:arith:Riemann:Roch:gen:finite:morphism},
\[
\acherncl_1 \left( \det Rg_*(G_n), h_Q^{\overline{G}_n} \right) -
\binom{n+r^2-1}{r^2 -1}
\acherncl_1 \left( \det Rg_*(\OO_{B}), 
h_Q^{\overline{\OO}_{B}} \right) \in \aBChow^1(Y;S)
\]
and
\[
\acherncl_1 \left( \det Rg_*(G_n), h_Q^{\overline{G}_n} \right) -
\binom{n+r^2-1}{r^2 -1}
\acherncl_1 \left( \det Rg_*(\OO_{B}), 
h_Q^{\overline{\OO}_{B}}
\right)
= g_* \left( \acherncl_1 (\overline{G}_n) \right).
\]
As before, we can take a $C^{\infty}$ Hermitian metric $h_{B/Y}$
of $\det Rg_*(\OO_{B})$ over $Y(\CC)$ and
a real valued $\varphi \in \LocInt(Y(\CC))$ such that
$h_Q^{\overline{\OO}_{B}} = e^{\varphi} h_{B/Y}$,
$F_{\infty}^*(h_{B/Y}) = \overline{h}_{B/Y}$, 
$F_{\infty}^*(\varphi) = \varphi \ (\alev)$,
and $\varphi(z) = 0$ for all $z \in S$. We set
\[
g_n = \exp \left( - \binom{n+r^2-1}{r^2-1} \varphi \right) 
h^{\overline{G}_n}_Q.
\]
Then,
\[
\acherncl_1 \left( \det Rg_*(G_n), g_n \right) -
\binom{n+r^2-1}{r^2 -1}
\acherncl_1 \left( \det Rg_*(\OO_{B}), h_{B/Y} \right)
= g_* \left( \acherncl_1 (\overline{G}_n) \right)
\]
and $\acherncl_1 \left( \det Rg_*(G_n), g_n \right)
\in \aBChow^1(Y;S)$.
Moreover, in the same as in Claim~\ref{claim:terms:right:R:R:formla},
we can see that
\[
g_* \left( \acherncl_1 (\overline{G}_n) \right) =
a(\deg(g)b_n) + \binom{n+r^2-1}{r^2-1}
\left( (n+1) g_* \acherncl_1(g^*(\overline{H})) +
g_* \acherncl_1(\iota^*(\overline{A})) \right).
\]
Thus, there are $W'_0, \ldots, W'_{r^2} \in
\aBChow^1(Y;S)_{\QQ}$ such that
\[
\acherncl_1 \left( \det Rg_*(G_n), g_n \right)
= \sum_{i=0}^{r^2} W'_i n^i + a(b_n \deg(g)).
\]
Therefore, we have
\addtocounter{Claim}{1}
\begin{equation}
\label{eqn:2:proof:arith:BG:inq}
\lim_{n \to \infty}
\frac{\acherncl_1 \left( \det Rg_*(G_n), g_n \right)}{n^{r^2+1}}
= 0
\end{equation}
in $\aBChow^1(Y;S)_{\QQ}$.

Let us consider an exact sequence:
\[
0 \to F_n \otimes A^{\otimes -2} \otimes L^{\otimes -1} \to F_n \to \rest{F_n}{M} \to 0.
\]
Since $F$ is semi-stable and of degree $0$
along $X_y$ and $\rest{L}{X_y} = \OO_{X_y}$,
we have 
\[
f_*(F_n \otimes A^{\otimes -2}\otimes L^{\otimes -1}) = 0
\]
on $Y$. Thus, the above exact sequence gives rise to
\[
0 \to f_*(F_n) \to (\rest{f}{M})_*(\rest{F_n}{M})
\to R^1f_*(F_n \otimes A^{\otimes -2}\otimes L^{\otimes -1})
\to R^1 f_*(F_n).
\]
Let $Q_n$ be the cokernel of
\[
f_*(F_n) \to (\rest{f}{M})_*(\rest{F_n}{M}) \to g_*(G_n).
\]
Let $U$ be the maximal Zariski open set of $Y$ such that
$f$ is smooth over $U$ and $g$ is \'{e}tale over $U$.
Moreover, let $U_n$ be the maximal Zariski open set of $Y$
such that 
\[
\begin{cases}
\text{(a) $U_n \subset U$, } \\
\text{(b) $(\rest{f}{M})_*(\rest{F_n}{M})$ coincides with $g_*(G_n)$ over $U_n$,} \\
\text{(c) $R^1f_*(F_n) = 0$ over $U_n$, and} \\
\text{(d) $f_*(F_n)$, $g_*(G_n)$ and $Q_n$ are
locally free over $U_n$.}
\end{cases}
\]
Then, $y \in (U_n)_{\QQ}$ and $S \subseteq U_n(\CC)$.
For, since $A \otimes \omega_{X/Y}^{-1}$ is ample on $X_y$ and
$E$ is semi-stable on $X_y$,
we can see that $R^1 f_*(F_n) = 0$ around $y$, which implies that
$f_*(F_n)$ is locally free around $y$.
Further, since $f_*(F_n)$ and $(\rest{f}{M})_*(\rest{F_n}{M})$ are free at $y$,
$R^1 f_*(F_n) = 0$ around $y$, and
$(\rest{f}{M})_*(\rest{F_n}{M})$ coincides with $g_*(G_n)$ around $y$,
we can easily check that $Q_n$ is free at $y$. Thus, 
$y \in (U_n)_{\QQ}$.
In the same way, we can see that $S \subseteq U_n(\CC)$.

Next let us consider a metric of $\det Q_n$.
$g_*(G_n)$ has the Hermitian metric
$(\rest{f}{M})_*\left( \rest{k_n}{M} \right)$ over $U_n(\CC)$, 
where $k_n$ is the Hermitian metric of $\overline{F}_n$.
Let $\tilde{q}_n$ be the quotient metric of $Q_n$
over $U_n(\CC)$ induced by
$(\rest{f}{M})_*\left( \rest{k_n}{M} \right)$. Let $q_n$ be a $C^{\infty}$
Hermitian metrics of
$\det Q_n$ over $Y(\CC)$ such that $F_{\infty}^*(q_n) = q_n$ and
$q_n(z) = \det \tilde{q}_n(z)$ for all $z \in S$.
(If $q_n$ is not invariant under $F_{\infty}$, then consider
$(1/2)\left(q_n + \overline{F_{\infty}^*(q_n)}\right)$.)

Here since
$\det Rf_*(F_n) \simeq \det f_*(F_n) \otimes \left(\det R^1 f_*(F_n)\right)^{-1}$ 
and
$\det f_*(F_n) \simeq \det g_*(G_n) \otimes (\det Q_n)^{-1}$,
we have
\[
\det Rf_*(F_n) \simeq \det g_*(G_n) \otimes (\det Q_n)^{-1} \otimes 
\left(\det R^1 f_*(F_n)\right)^{-1}.
\]
Further, we have generalized metrics
$h_n$, $g_n$ and $q_n$ of
$\det Rf_*(F_n)$, $\det g_*(G_n)$ and $\det Q_n$.
Thus, there is a generalized metric $t_n$
of $\det R^1 f_*(F_n)$ such that
the above is an isometry.

As in the proof of Proposition~\ref{prop:find:small:section},
let us construct a section of $\det R^1f_*(F_n)$.
First, we fix a locally free sheaf $P_n$ on $Y$ and
a surjective homomorphism $P_n \to R^1 f_*(F_n)$.
Let $P'_n$ be the kernel of $P_n \to R^1 f_*(F_n)$.
Then, $P'_n$ is a torsion free sheaf and has the same rank
as $P_n$ because $R^1 f_*(F_n)$ is a torsion sheaf. 
Noting that $\left( \bigwedge^{\rank P'_n} P'_n \right)^{*}$
is an invertible sheaf on $Y$,
we can identify $\det R^1 f_*(F_n)$ with 
\[
\bigwedge^{\rank P_n} P_n \otimes 
\left( \bigwedge^{\rank P'_n} P'_n \right)^{*}.
\]
Moreover, the homomorphism
$\bigwedge^{\rank P'_n} P'_n \to \bigwedge^{\rank P_n} P_n$
induced by $P'_n \hookrightarrow P_n$ 
gives rise to
a non-zero section $s_n$ of $\det R^1 f_*(F_n)$.
Note that $s_n(y) \not= 0$ and $s_n(z) \not= 0$ 
for all $z \in S$ because $R^1 f_*(F_n) = 0$ at $y$ and $z$.

Here we set
\[
a_n = \max_{z \in S} \{ \log t_n(s_n,s_n)(z) \}.
\]
By our construction, we have
\[
\acherncl_1(\det R^1 f_*(F_n), e^{-a_n}t_n) \in \aBChow^1(Y;S).
\]
and an isometry
\addtocounter{Claim}{1}
\begin{multline}
\label{eqn:isometry:det:bundle}
(\det Rf_*(F_n), h_n) \simeq \\
(\det g_*(G_n), g_n) \otimes (\det Q_n, q_n)^{-1} \otimes
(\det R^1 f_*(F_n), e^{-a_n}t_n)^{-1} \otimes (\OO_Y, e^{-a_n}h_{can}).
\end{multline}

Here we claim:

\begin{Claim}
\label{claim:gen:small:sec:Q:n}
$(\det Q_n, q_n)$ is generated by small sections at $y$ with respect to $S$.
\end{Claim}

\Proof
First of all,
\[
g_*\left( \iota^*(F) \otimes g^*(H) \right) = 
g_*(\iota^*(F)) \otimes H
\quad\text{and}\quad
g_*\left( \iota^*(A) \otimes g^*(H) \right) = 
g_*(\iota^*(A)) \otimes H
\]
are generated by small section at $y$ with respect to $S$.
Thus, by (2) and (3) of Proposition~\ref{prop:find:small:section},
\[
g_*(G_n) = g_* \left(
\Sym^n (\iota^*(F) \otimes g^*(H)) \otimes \iota^*(A) \otimes g^*(H)
\right)
\]
is generated by small sections at $y$ with respect to $S$.
Thus, by (1) of Proposition~\ref{prop:find:small:section},
$(Q_n, \tilde{q}_n)$ is generated by small sections at $y$ with
respect to $S$.
Hence, by (4) of Proposition~\ref{prop:find:small:section},
$(\det Q_n, q_n)$ is generated by small sections at $y$
with respect to $S$ because
$q_n(z) = \det \tilde{q}_n (z)$ for all $z \in S$.
\QED

Next we claim:

\begin{Claim}
\label{claim:estimate:sequences:an:bn}
$a_n \leq O(n^{r^2} \log(n))$.
\end{Claim}

\Proof
It is sufficient to show that
$\log t_n(s_n,s_n)(z) \leq O(n^{r^2}\log(n))$
for each $z \in S$.
Let $\{ e_1, \ldots, e_{l_n} \}$
be an orthonormal basis of 
$g_*(G_n) \otimes \kappa(z)$ with respect to
$g_*(\rest{k_n}{B})(z)$ such that
$\{ e_1, \ldots, e_{m_n} \}$ forms a basis
of $f_*(F_n) \otimes \kappa(z)$.
Then, $e_1 \wedge \cdots \wedge e_{m_n}$,
$e_1 \wedge \cdots \wedge e_{l_n}$ and
$\bar{e}_{m_n + 1} \wedge \cdots \wedge \bar{e}_{l_n}$
form bases of
$\det (f_*(F_n)) \otimes \kappa(z)$,
$\det (g_*(G_n)) \otimes \kappa(z)$, and
$\det (Q_n) \otimes \kappa(z)$ respectively, and
$(e_1 \wedge \cdots \wedge e_{m_n}) \otimes
(\bar{e}_{m_n + 1} \wedge \cdots \wedge \bar{e}_{l_n}) =
e_1 \wedge \cdots \wedge e_{l_n}$,
where $\bar{e}_{m_n + 1}, \ldots, \bar{e}_{l_n}$
are images of $e_{m_n + 1}, \ldots, e_{l_n}$ in
$Q_n \otimes \kappa(z)$.
Then,
\[
\left| (e_1 \wedge \cdots \wedge e_{m_n}) \otimes s_n^{\otimes -1} 
\right|_{h_n}^2(z) =
\frac{|e_1 \wedge \cdots \wedge e_{l_n}|_{g_n}^2(z)}{
|\bar{e}_{m_n + 1} \wedge \cdots \wedge \bar{e}_{l_n}|_{q_n}^2(z)
|s_n|_{t_n}^2(z)} = |s_n|_{t_n}^{-2}(z),
\]
where $|a|_{\lambda} = \sqrt{\lambda(a,a)}$ for
$\lambda = h_n, g_n, q_n, t_n$.
Moreover, let $\Omega_{z}$ be the K\"{a}hler form
induced by the metric of $\overline{\omega}_{X/Y}$ along $f_{\CC}^{-1}(z)$.
Then, there is a Hermitian metric $v_n$ of $H^0(f_{\CC}^{-1}(z), F_n)$
defined by
\[
v_n(s, s') 
= \int_{f_{\CC}^{-1}(z)} k_n(s, s') \Omega_{z}.
\]
Here $R^1 f_*(F_n) = 0$ at $z$. Thus,
$(\det R^1 f_*(F_n))_z$ is canonically isomorphic to $\OO_{Y(\CC), z}$.
Since $(P'_n)_z = (P_n)_z$, under the above isomorphism,
$s_n$ goes to the determinant of
$(P_n)_z \overset{\operatorname{id}}{\longrightarrow} (P_n)_z$,
namely $1 \in \OO_{Y(\CC), z}$.
Hence, by the definition of Quillen metric,
\[
\left| (e_1 \wedge \cdots \wedge e_{m_n}) \otimes s_n^{\otimes -1} 
\right|_{h_n}^2(z)
= \det( v_n(e_i, e_j) )
\exp \left( -T \left( \rest{\overline{F}_n}{f_{\CC}^{-1}(z)} 
\right) \right).
\]
Therefore,
\[
\log |s_n|_{t_n}^2(z) = 
T \left( \rest{\overline{F}_n}{f_{\CC}^{-1}(z)} \right) - 
\log \det( v_n(e_i, e_j) ).
\]
By Corollary~\ref{cor:asymp:analytic:torsion},
\[
T \left( \rest{\overline{F}_n}{f_{\CC}^{-1}(z)} \right)
\leq O(n^{r^2} \log(n)).
\]
Thus, in order to get our claim,
it is sufficient to show that
\[
- \log \det( v_n(e_i, e_j) ) \leq  O(n^{r^2-1}\log(n)).
\]
Let $s$ be an arbitrary section of $H^0(f_{\CC}^{-1}(z), F_n)$.
Then, by Lemma~\ref{lem:comparison:sup:L2},
\[
g_*\left( \rest{k_n}{B} \right)(s, s) = 
\sum_{x \in g_{\CC}^{-1}(z)} |s|_{k_n}^2(x) \leq
\deg(g) \sup_{x \in f_{\CC}^{-1}(z)} \{ |s|_{k_n}^2(x) \} \leq
\deg(g) c^2 n^{2r^2} \Vert s \Vert^2_{L^2}
\]
for some constant $c$ independent of $n$.
Thus, by \cite[Lemma~3.4]{MoBG} and our choice of $e_i$'s,
\[
1 = \det \left( g_*\left(\rest{k_n}{B}\right)(e_i, e_j) \right)
\leq
\left( \deg(g) c^2 n^{2r^2} \right)^{\dim_{\CC} H^0(f_{\CC}^{-1}(z), F_n)}
\det \left( v_n(e_i, e_j) \right).
\]
Using Riemann-Roch theorem, we can easily see that
\[
\dim_{\CC} H^0(f_{\CC}^{-1}(z), F_n) \leq O(n^{r^2-1}).
\]
Thus, we have
\[
-\log \det( v_n(e_i, e_j) ) \leq  O(n^{r^2-1}\log(n)).
\]
Hence, we obtain our claim.
\QED

Let us go back to the proof of our theorem.
By the isometry (\ref{eqn:isometry:det:bundle}), we get
\begin{align*}
-\acherncl_1(\det Rf_*(F_n), h_n) & =
-\acherncl_1( \det g_*(G_n), g_n) + \acherncl_1(\det Q_n, q_n)
+ \acherncl_1(\det R^1 f_*(F_n), e^{-a_n}t_n) - a(a_n) \\
& = \left[
\acherncl_1(\det Q_n, q_n) + \acherncl_1(\det R^1 f_*(F_n),
e^{-a_n}t_n) + a\left( \max \{ -a_n, 0 \} \right) \right]
\\ & \qquad\qquad
+ \left[ -\acherncl_1( \det g_*(G_n), g_n) +
a (\min \{ -a_n, 0 \}) \right].
\end{align*}
Here we set
\[
\begin{cases}
{\displaystyle
\alpha_n = \frac{(r^2+1)!}{n^{r^2+1}} \left[
\acherncl_1(\det Q_n, q_n) + \acherncl_1(\det R^1 f_*(F_n), e^{-a_n}t_n)
+ a\left( \max \{ -a_n, 0 \} \right)
\right],} \\
{} \\
{\displaystyle
\beta_n = \frac{(r^2+1)!}{n^{r^2+1}} \left[
-\acherncl_1( \det g_*(G_n), g_n) +
a (\min \{ -a_n, 0 \})
\right].}
\end{cases}
\]
Then,
\[
\frac{-(r^2+1)!\acherncl_1(\det Rf_*(F_n), h_n)}{n^{r^2+1}} =
\alpha_n + \beta_n.
\]
By (\ref{eqn:2:proof:arith:BG:inq})
and Claim~\ref{claim:estimate:sequences:an:bn},
${\displaystyle \lim_{n \to \infty} \beta_n = 0}$
in $\aBChow^1(Y;S)_{\QQ}$.
Therefore, by Claim~\ref{claim:dis:lim:c1:Fn},
\[
\adis_{X/Y}(\overline{E}) =
\lim_{n \to \infty} 
\frac{-(r^2+1)! \acherncl_1 \left( \det Rf_*(F_n), h_n \right)}
{n^{r^2+1}} = \lim_{n\to\infty} (\alpha_n + \beta_n) =
\lim_{n \to \infty} \alpha_n
\]
in $\aBChow^1(Y;S)_{\QQ}$.
On the other hand, it is obvious that
\[
\acherncl_1(\det R^1 f_*(F_n), e^{-a_n}t_n)
\quad\text{and}\quad
a\left( \max \{ -a_n, 0 \} \right)
\]
is semi-ample at $y$ with respect to $S$.
By Claim~\ref{claim:gen:small:sec:Q:n},
$\acherncl_1(\det Q_n, q_n)$ is semi-ample at $y$ with respect to $S$.
Thus, $\alpha_n$ is semi-ample at $y$ with respect to $S$.
Hence we get our theorem.
\QED
\renewcommand{\theClaim}{\arabic{section}.\arabic{subsection}.\arabic{Theorem}.\arabic{Claim}}
\renewcommand{\theequation}{\arabic{section}.\arabic{subsection}.\arabic{Theorem}.\arabic{Claim}}


\section{Preliminaries for Cornalba-Harris-Bost's inequality}
\label{section:Bost:type:inequality:preparation}

This section is a preparatory one for the next section, 
where we will prove the relative Cornalba-Harris-Bost's inequality 
(cf. Theorem~\ref{thm:semistability:imply:average:semi-ampleness}).  
Moreover, in the next section, we will see how 
the relative Bogomolov's inequality 
(Theorem~\ref{thm:relative:Bogomolov:inequality:arithmetic:case})  
and 
the relative Cornalba-Harris-Bost's inequality 
(Theorem~\ref{thm:semistability:imply:average:semi-ampleness}) 
are related 
(cf. Proposition~\ref{prop:Bogomolov:to:Bost}). 

\subsection{Normalized Green forms}
\label{subsec:normalized:Green:form}
\setcounter{Theorem}{0}

Let $Y$ be a smooth quasi-projective variety over $\CC$, 
$\overline{E}=(E,h)$ a Hermitian vector bundle of rank $r$ on $Y$. 
Let $\pi : \PP(E) \to Y$ be the canonical morphism, 
where $\PP(E) = \Proj (\bigoplus_{i \ge 0} \Sym^{i}(E^{\lor}))$. 
We equip the canonical quotient bundle $\OO_{E}(1)$ on $\PP(E)$ 
with the quotient metric via  
$\pi ^* (E^{\vee}) \to \OO_{E}(1)$. 
We will denote this Hermitian line bundle by $\overline{\OO_{E}(1)}$. 
Furthermore, let $\Omega = \cherncl_1(\overline{\OO_{E}(1)})$ be 
the first Chern form. 

The purpose of this subsection is that, for  
every cycle $X \subset \PP(E)$ 
whose all irreducible components map surjectively to $Y$, 
we give a Green form $g_X$ 
such that on a general fiber,
it is an $\Omega$-normalized 
Green current in the sense of \cite[2.3.2]{BGS}.

Let $X$ be a cycle  of codimension $p$ on $\PP(E)$ such that 
every irreducible component of $X$ maps surjectively to $Y$.
An $L^1$-form $g_X$ on $\PP(E)$ 
satisfying the following conditions is called
an {\em $\Omega$-normalized Green form}, 
(or simply a {\em normalized Green form} when 
no confusion is likely).
\begin{enumerate}
\renewcommand{\labelenumi}{(\roman{enumi})}
\item
There are $d$-closed $L^1$-forms $\gamma_i$ of type $(p-i, p-i)$
on $Y$ ($i=0, \ldots, p$) with
\[
dd^c([g_X]) + \delta_X 
= \sum_{i=0}^{p} \left[ \pi^*(\gamma_i) \wedge \Omega^i \right].
\]

\item
$\pi_*(g_X \wedge \Omega^{r-p}) = 0$.
\end{enumerate}
Note that $\gamma_p$ is the degree of $X$ along a general fiber of $\pi$.

Let $X = \sum_i a_i X_i$ be the irreducible decomposition
of $X$ as cycles.
Let $\tilde{X}_i \to X_i$ be a desingularization of $X_i$, and
$\tilde{f}_i : \tilde{X}_i \to Y$ the induced morphism.
The main result of this subsection is the following.

\begin{Proposition}
\label{prop:normalized:Green:form}
With notation as above,
there exists an $\Omega$-normalized Green form $g_X$ on $\PP(E)$
satisfying the following property.
If $y \in Y$ and
$\tilde{f}_i$ is smooth over $y$ for every $i$,
then there is an open set $U$ containing $y$ such that
$\gamma_0, \ldots, \gamma_p$ are $C^{\infty}$ on $U$ and that
$\rest{g_X}{\pi^{-1}(U)}$ is a Green form of logarithmic type for $X_U$,
where $\gamma_0, \ldots, \gamma_p$ are $L^1$-forms
in the definition of $\Omega$-normalized Green form. 
\end{Proposition}

To prove the above proposition,
let us begin with the following two lemmas.
 
\begin{Lemma}
\label{lemma:auxiliary:green:form}
There exist a Green form $g$ of logarithmic type along $X$, and
$d$-closed $C^{\infty}$ forms $\beta_i$ of type $(p-i, p-i)$
on $Y$ \rom{(}$i=0, \ldots, p$\rom{)} such that 
\[
dd^c([g]) + \delta_X 
= \sum_{i=0}^{p} \left[ \pi^*(\beta_i) \wedge \Omega^i \right].
\]
\end{Lemma}

\Proof
We divide the proof into three steps.

{\bf Step 1.} : The case where $Y$ is projective.

Let $g_1$ be a Green form of logarithmic type along $X$ such that 
\[
dd^c([g_1]) + \delta_X 
= [\omega]
\]
where $\omega$ is a smooth form on $\PP(E)$.
Then, we can find a smooth form $\eta$ on $\PP(E)$ of the form 
\[
\eta = \sum_{i=0}^{p} \pi^*(\beta_i) \wedge \Omega^i
\]
which represents the same cohomology class as $\omega$,
where $\beta_i$ is a $d$-closed $C^{\infty}$-form of type
$(p-i, p-i)$ on $Y$. 
Since $\omega - \eta$ is $d$-exact $(p,p)$-form, 
by the $dd^c$-lemma, there is a smooth $(p-1,p-1)$-form $\phi$ with 
$\omega - \eta = dd^c(\phi)$. 
Thus, if we set $g = g_1 - \phi$, 
then $g$ is of logarithmic type along $X$ and 
\[
dd^c([g]) + \delta_X = dd^c([g_1]) - dd^c(\phi) + \delta_X = [\eta].
\]

\medskip
{\bf Step 2.} : Let $h'$ be another Hermitian metric of $E$, and
$\Omega'$ the Chern form of $\OO_{E}(1)$ arising from $h'$.
In this step, we will prove that
if the lemma holds for $h'$, then so does it for $h$.

By our assumption,
there exist a Green form $g'$ of logarithmic type along $X$, and
$d$-closed $C^{\infty}$ forms $\beta'_i$ ($i=0, \ldots, p$)
of type $(p-i, p-i)$
on $Y$ such that 
\[
dd^c([g']) + \delta_X 
= \sum_{i=0}^{p} \left[ \pi^*(\beta'_i) \wedge {\Omega'}^i \right].
\]
On the other hand, there is a real $C^{\infty}$-function $a$ on
$\PP(E)$ with $\Omega' - \Omega = dd^c(a)$.
Here note that if $v$ is a $\partial$ and $\overline{\partial}$-closed 
form on $\PP(E)$, then
$dd^c(v \wedge a) = v \wedge dd^c(a)$.
Thus, it is easy to see that there is a
$C^{\infty}$ form $\theta$ on $\PP(E)$ such that
\[
\sum_{i=1}^p \pi^*(\beta'_i) \wedge {\Omega'}^i = 
dd^c(\theta) + \sum_{i=1}^p \pi^*(\beta'_i) \wedge {\Omega}^i.
\]
Therefore, if we set $g = g' - \theta$ and $\beta_i = \beta'_i$,
then we have our assertion for $h$.

\medskip
{\bf Step 3.} : General case.

Using Hironaka's resolution \cite{Hiro}, 
there is a smooth projective variety
$Y'$ over $\CC$ such that $Y$ is an open set of $Y'$.
Moreover, using \cite[Exercise~5.15 in Chapter~II]{Hartshorne},
there is a coherent sheaf $E'$ on $Y'$ with
$\rest{E'}{Y} = E$.
Further, taking a birational modification along $Y' \setminus Y$
if necessary, we may assume that $E'$ is locally free.
Let $h'$ be a Hermitian metric of $E'$ over $Y'$.
Since $\PP(E)$ is an Zariski open set of $\PP(E')$,
let $X'$ be the closure of $X$ in $\PP(E')$.
Then, by Step~1, our assertion holds for 
$(E', h')$ and $X'$.
Thus, so does it for $(E, \rest{h'}{Y})$ and $X$.
Therefore, by Step~2, we can conclude our lemma.
\QED

\begin{Lemma}
\label{lemma:push:g:is:L1}
Let $g$ be a Green form of logarithmic type along $X$ 
and $\omega$ a $C^{\infty}$-form with $dd^c([g]) + \delta_X = [\omega]$. 
If we set $\varsigma = \pi_* (g \wedge \Omega^{r-p})$, then 
$\varsigma \in L^1_{loc}(Y)$ and 
$dd^c([\varsigma]) \in L^1_{loc}(\Omega_Y^{1,1})$.
Moreover, if $y \in Y$ and $\tilde{f}_i$ is smooth over $y$
for every $i$,
then $\varsigma$ is $C^{\infty}$ around $y$.
\end{Lemma}

\Proof
By Proposition~\ref{prop:push:forward:B:pq}, 
$\varsigma$ is an $L^1$-function on $Y$ and  
\begin{align*}
dd^c([\varsigma]) 
 & = dd^c( \pi_* ([g \wedge \Omega^{r-p}]))
   = \pi_* dd^c([g \wedge \Omega^{r-p}]) \\
 & = \pi_* dd^c([g]) \wedge \Omega^{r-p} 
 = \pi_* ([\omega] \wedge \Omega^{r-p})
     - \pi_* (\delta_X \wedge \Omega^{r-p}) \\
 & = \pi_* [\omega \wedge \Omega^{r-p}] 
     - \sum_i a_i \pi_* (\delta_{X_i} \wedge \Omega^{r-p}) \\
 & = \pi_* [\omega \wedge \Omega^{r-p}] 
     - \sum_i a_i (\tilde{f}_i)_* [\tilde{f}_i^* (\Omega^{r-p})]. 
\end{align*}
Thus, $dd^c([\varsigma]) \in L^1_{loc}(\Omega_Y^{1,1})$.
Moreover,  if $y \in Y$ and $\tilde{f}_i$ is smooth over $y$
for every $i$, then, by the above formula,
$dd^c([\varsigma])$ is $C^{\infty}$ around $y$.
Thus, by virtue of
\cite[(i) of Theorem~1.2.2]{GSArInt},
$\varsigma$ is $C^{\infty}$ around $y$.
\QED

Let us start the proof of Proposition~\ref{prop:normalized:Green:form}.
Let $g$ be a Green form constructed in 
Lemma~\ref{lemma:auxiliary:green:form}. 
Then, there are $d$-closed $\beta_i$'s with
$\beta_i \in A^{p-i,p-i} (Y)$ and
\[
dd^c([g]) + \delta_X 
= \sum_{i=0}^{p} \left[ \pi^*(\beta_i) \wedge \Omega^i \right].
\]
If we set 
$\varsigma = \pi_* (g \wedge \Omega^{r-p})$, 
then by Lemma~\ref{lemma:push:g:is:L1}, 
$\varsigma$ is locally an $L^1$-form. We put 
\[
g_X = g - \pi^*(\varsigma) \Omega^{p-1}, 
\] 
which is clearly locally an $L^1$-form on $\PP(E)$. 
We will show that $g_X$ satisfies the conditions (i) and (ii). 
Using $\int_{\PP(E) \to Y} \Omega^{r-1} = 1$, (ii) can be readily checked. 
Moreover,
\begin{align*}
dd^c([g_X]) + \delta_X 
& = \sum_{i=0}^{p} \left[ \pi^*(\beta_i) \wedge \Omega^i \right] 
- dd^c [\pi^*(\varsigma) \Omega^{p-1}] \\
& = \beta_p \Omega^p + 
\pi^*([\beta_{p-1}] - dd^c([\varsigma])) \wedge \Omega^{p-1} +
\sum_{i=0}^{p-2} \left[ \pi^*(\beta_i) \wedge \Omega^i \right].
\end{align*}
The remaining assertion is easily derived from
Lemma~\ref{lemma:push:g:is:L1}.
\QED

\begin{Remark}
\label{rem:norm:Green:general:fiber}
Let $y$ be a point of $Y$ such that
$\tilde{f}_i$ is smooth over $y$ for every $i$.
Then, by Proposition~\ref{prop:normalized:Green:form},
on the fiber $\pi^{-1}(y)$, 
$\rest{g_X}{\pi^{-1}(y)}$ is a Green form of logarithmic type along
$X_{y}$. Moreover,
\[
dd^c ([\rest{g_X}{\pi^{-1}(y)}]) + \delta_{X_y} = \deg (X_y) 
[ \rest{\Omega^p}{\pi^{-1}(y)} ]
\]
and
\[
\int_{\pi^{-1}(y)} \left( \rest{g_X}{\pi^{-1}(y)} \right)
\left( \rest{\Omega^{r-p}}{\pi^{-1}(y)} \right) = 0.
\]
Thus, $\rest{g_X}{\pi^{-1}(y)}$ is a
$\Omega$-normalized Green form on $\pi^{-1}(y)$, and it is also
a $\Omega$-normalized Green current in the sense of \cite[2.3.2]{BGS}.
\end{Remark}

\subsection{Associated Hermitian vector bundles}
\label{subsec:associated:herm:vb}
\renewcommand{\theClaim}{\arabic{section}.\arabic{subsection}.\arabic{Theorem}}
\renewcommand{\theequation}{\arabic{section}.\arabic{subsection}.\arabic{Theorem}}
\setcounter{Theorem}{0}

Let $\bfGL_r = \Spec \ZZ [X_{11},X_{12},\cdots,X_{rr}]_{\det(X_{ij})}$ 
be the general linear group of rank $r$ and 
$\bfSL_r = \Spec \ZZ [X_{11},X_{12},\cdots,X_{rr}]/(\det(X_{ij})-1)$ 
be the special linear group of rank $r$. 

Let $\rho : \bfGL_r \to \bfGL_R$ be a morphism of group schemes. 
First, we note that 
\[
\rho(\CC) (\overline{A}) = \overline{\rho(\CC) (A)},
\]
where $\rho(\CC) : \bfGL_r(\CC) \to \bfGL_R(\CC)$ is the induced morphism 
and $A \in \bfGL_r(\CC)$. 
Indeed, the above equality is nothing but the associativity of the map 
\[
\Spec \CC \overset{-}{\longrightarrow} \Spec \CC 
\overset{A}{\longrightarrow}    \bfGL_r 
\overset{\rho}{\longrightarrow} \bfGL_R.
\]
Next, we consider the following condition for $\rho$;
\addtocounter{Theorem}{1}
\begin{equation}
\label{eqn:condition}
\rho({}^t A) = {}^t \rho (A) \qquad \text{for any $A \in \bfGL_r$}. 
\end{equation}
In the group scheme language, 
this condition means $\rho$ commutes with the transposed morphism. 

Let $\Unitary_r(\CC) 
= \{A \in \bfGL_r(\CC) \mid {}^t A \cdot \overline{A} = I_r \}$ 
be the unitary group of rank $r$. 
If a group morphism $\rho : \bfGL_r \to \bfGL_R$ 
commutes with the transposed morphism, then 
\begin{align*}
 I_R  & = \rho(\CC) (I_r) = \rho(\CC) ({}^t A \cdot \overline{A}) \\
      & = \rho(\CC)({}^t A) \cdot \rho(\CC)(\overline{A}) 
      = {}^t \rho(\CC)(A) \cdot \overline{\rho(\CC)(A)},
\end{align*}
namely, $\rho(\CC)$ maps  $\Unitary_r(\CC)$ into $\Unitary_R(\CC)$.

Let $k$ be an integer. 
A morphism $\rho : \bfGL_r \to \bfGL_R$ of group schemes is said to be 
{\em of degree $k$} if 
\[
\rho(t I_r) = t^k I_R \qquad \text{for any $t$}.
\]
In the group scheme language, 
this means that the diagram 
\begin{equation*}
 \begin{CD}
  \bfGL_1         @>{\lambda_r}>>  \bfGL_r        \\
  @V{\alpha}VV                  @VV{\rho}V   \\
  \bfGL_1         @>{\lambda_R}>>  \bfGL_R
 \end{CD}
\end{equation*}
commutes, where $\lambda_r$ and $\lambda_R$ are given by 
$t \mapsto \diag(t,t,\cdots,t)$ and $\alpha$ is given by $t \mapsto t^k$. 

\medskip
Let $Y$ be an arithmetic variety, 
$\overline{E} = (E,h)$ a Hermitian vector bundle of rank $r$ on $Y$ and 
$\rho : \bfGL_r \to \bfGL_R$ be 
a morphism of group schemes satisfying 
commutativity with the transposed morphism. 
In the following, we will show that we can naturally construct 
a Hermitian vector bundle 
$\overline{E}^{\rho} = (E^{\rho},h^{\rho})$, 
which we will call the 
{\em associated Hermitian vector bundle} with respect to 
$\overline{E}$ and $\rho$. 

First, we construct $E^{\rho}$. 
Let $\{ Y_{\alpha} \}$ be an affine open covering such that 
$\phi_{\alpha}: E \vert_{Y_{\alpha}}
\overset{\sim}{\longrightarrow} \OO_{Y_{\alpha}}^{\oplus r}$ 
gives a local trivialization. 
On $Y_{\alpha} \cap Y_{\beta}$, we set the transition function 
$g_{\alpha\beta} = \phi_{\alpha} \cdot \phi_{\beta}^{-1}$, 
which can be seen as an element of 
$\bfGL_r(\Gamma(\OO_{Y_{\alpha} \cap Y_{\beta}}))$. 
Then we define the associated vector bundle $E^{\rho}$ 
as the vector bundle of rank $R$ on $Y$ 
with the transition functions 
$\rho(\Gamma(\OO_{Y_{\alpha} \cap Y_{\beta}})) (g_{\alpha\beta})$; 
\[
E^{\rho} = \coprod_{\alpha} \OO_{Y_{\alpha}}^{\oplus R} / \sim.
\]

Next, we define metric on $E^{\rho}$. 
Let $h^{\alpha}$ be the Hermitian metric 
on $\OO_{Y_{\alpha}}^{\oplus r}$ over $Y_{\alpha}$ 
such that 
$\phi_{\alpha}: E \vert_{Y_{\alpha}} 
\overset{\sim}{\longrightarrow} \OO_{Y_{\alpha}}^{\oplus r}$ 
becomes isometry over $Y_{\alpha}(\CC)$. Let 
\begin{gather*}
e_1^{\alpha} = {}^t (1,0,\cdots,0), \cdots, e_r^{\alpha} = {}^t (0,\cdots,0,1) 
\in \Gamma(\OO_{Y_{\alpha}}^{\oplus r}), \\
f_1^{\alpha} = {}^t (1,0,\cdots,0), \cdots, f_R^{\alpha} = {}^t (0,\cdots,0,1) 
\in \Gamma(\OO_{Y_{\alpha}}^{\oplus R})
\end{gather*}
be the standard local frames of 
$\OO_{Y_{\alpha}}^{\oplus r}$ and $\OO_{Y_{\alpha}}^{\oplus R}$. 
We set 
\[
H_{\alpha} = (h^{\alpha}(e_{i}^{\alpha},e_{j}^{\alpha}))_{1 \leq i,j \leq r}.
\]
Then $H_{\alpha}$ is a $C^{\infty}$-map over $Y_{\alpha}(\CC)$ 
and, for each point $y$ in $Y_{\alpha}(\CC)$, 
$H_{\alpha}(y)$ is a positive definite Hermitian matrix. 
Let $\rho(C^{\infty}(Y_{\alpha}(\CC))) : 
\bfGL_r(C^{\infty}(Y_{\alpha}(\CC))) \to \bfGL_R(C^{\infty}(Y_{\alpha}(\CC)))$ 
be the induced map. 

\addtocounter{Theorem}{1}
\begin{Claim}
$\rho(C^{\infty}(Y_{\alpha}(\CC))) (H_{\alpha})$ 
is a $C^{\infty}$-map over $Y_{\alpha}(\CC)$ and, 
for each point $y$ in $Y_{\alpha}(\CC)$, 
$\rho(C^{\infty}(Y_{\alpha}(\CC)))(H_{\alpha})(y)$ 
is a positive definite Hermitian matrix. 
\end{Claim}

\Proof
The first assertion is obvious. For the second one, we note that 
there is a matrix $A \in \bfGL_r(\CC)$ 
such that ${}^t A \cdot \overline{A} = H_{\alpha}(y)$. 
Then it is easy to see that $\rho(C^{\infty}(Y_{\alpha}(\CC)))(H_{\alpha})(y)$ 
is a positive definite Hermitian matrix by using \eqref{eqn:condition}.
\QED

Now we define a metric $h^{\rho_{\alpha}}$ 
on $\OO_{Y_{\alpha}}^{\oplus R}$ over $Y_{\alpha}$ by 
\[
h^{\rho_{\alpha}}(f_{k}^{\alpha},f_{l}^{\alpha})
= \rho(C^{\infty}(Y_{\alpha}(\CC))) (H_{\alpha})_{kl}
\]
for $1 \leq k,l \leq R$. 

\addtocounter{Theorem}{1}
\begin{Claim}
$\{h^{\rho_{\alpha}}\}_{\alpha}$ glue together 
to form a Hermitian metric on $E^{\rho}$.
\end{Claim}

\Proof
Let $s_{\alpha} = {}^t (s_1^{\alpha},\cdots,s_R^{\alpha}) 
\in \Gamma(\OO_{Y_{\alpha}}^{\oplus R} \vert_{Y_{\alpha} \cap Y_{\beta}})$ 
and $s_{\beta} = {}^t (s_1^{\beta},\cdots,s_R^{\beta}) 
\in \Gamma(\OO_{Y_{\beta}}^{\oplus R} \vert_{Y_{\alpha} \cap Y_{\beta}})$. 
Then they give the same section of 
$E^{\rho} \vert_{Y_{\alpha} \cap Y_{\beta}}$ if 
${}^t (s_1^{\alpha},\cdots,s_R^{\alpha}) = 
\rho(g_{\alpha\beta}) {}^t (s_1^{\beta},\cdots,s_R^{\beta})$. 
In this case, we write $s_{\alpha} \sim s_{\beta}$. 
Now we take $s_{\alpha} \sim s_{\beta}$ and $t_{\alpha} \sim t_{\beta}$. 
Then by a straightforward calculation using \eqref{eqn:condition} and 
$H_{\beta} = {}^t g_{\alpha\beta} H_{\alpha} \overline{g_{\alpha\beta}}$, 
we get $h^{\rho_{\alpha}}(s_{\alpha},t_{\alpha}) =
h^{\rho_{\beta}}(s_{\beta},t_{\beta})$ on $Y_{\alpha} \cap Y_{\beta}$.
\QED

\begin{Remark}
Let $\operatorname{id}_r : \bfGL_r \to \bfGL_r$ be the identity morphism, 
$\rho_1 = (\operatorname{id}_r)^{\otimes k}$, 
$\rho_2 = \Sym^k(\operatorname{id}_r)$, 
and $\rho_3 = \bigwedge^k(\operatorname{id}_r)$. 
Further, let $\rho_4 : \bfGL_r \to \bfGL_r$ be 
the group homomorphism given by $A \mapsto {}^t A^{-1}$. 
Then $\rho_1$, $\rho_2$, $\rho_3$ and $\rho_4$ are 
of degree $k$, $k$, $k$ and $-1$, respectively. 
Let $(E,h)$ be a Hermitian vector bundle of rank $r$. Then 
the associated vector bundles are 
$(E^{\otimes k}, h^{\otimes k})$, $(\Sym^k(E), h^{\rho_2})$, 
$(\bigwedge^k(E), h^{\rho_3})$ and $(E^{\lor}, h^{\lor})$. 
Note, for example, that 
$h^{\rho_2}$ is not the quotient metric $h_{quot}$ given by 
$E^{\otimes k} \to \Sym^k(E)$; Indeed, for a locally orthogonal basis 
$e_1,\cdots,e_r$ of $\overline{E}$ and $\alpha_1,\cdots,\alpha_r \in \ZZ$, 
$h^{\rho_2}(e_1^{\alpha_1}\cdots e_r^{\alpha_r}, 
e_1^{\alpha_1}\cdots e_r^{\alpha_r}) = 1$, while 
$h_{quot}(e_1^{\alpha_1}\cdots e_r^{\alpha_r}, 
e_1^{\alpha_1}\cdots e_r^{\alpha_r}) 
= {\alpha_1 !\cdots\alpha_r !}/ r !$.
\end{Remark}

\renewcommand{\theClaim}{\arabic{section}.\arabic{subsection}.\arabic{Theorem}.\arabic{Claim}}
\renewcommand{\theequation}{\arabic{section}.\arabic{subsection}.\arabic{Theorem}.\arabic{Claim}}

\subsection{Chow forms and their metrics}
\label{subsec:Chow:forms:and:their:metrics}
\renewcommand{\theClaim}{\arabic{section}.\arabic{subsection}.\arabic{Theorem}}
\renewcommand{\theequation}{\arabic{section}.\arabic{subsection}.\arabic{Theorem}}
\setcounter{Theorem}{0}
Let $Y$ be a regular arithmetic variety, and 
$\overline{E}=(E,h)$ a Hermitian vector bundle of rank $r$ on $Y$.  

Let $\rho : \bfGL_r \to \bfGL_R$ be 
a group scheme morphism of degree $k$  
commuting with the transposed morphism 
and 
$\overline{E}^{\rho}=(E^{\rho},h^{\rho})$ 
the associated Hermitian bundle of rank $R$. 
We give the quotient metric on $\OO_{E^{\rho}}(1)$ 
via $\pi^*({E^{\rho}}^{\vee}) \to \OO_{E^{\rho}}(1)$. 
We denote this Hermitian line bundle by 
$\overline{\OO_{E^{\rho}}(1)}$. Further, 
let $\Omega_{\rho} = \cherncl_1(\overline{\OO_{E^{\rho}}(1)})$ 
be the first Chern form. 

Let $X$ be an effective cycle in $\PP (E^{\rho})$ 
such that $X$ is flat over $Y$ 
with the relative dimension $d$ 
and degree $\delta$ on the generic fiber. 
Let $g_X$ be a $\Omega_{\rho}$-normalized Green form for $X$ 
and we set $\widehat{X}= (X, g_X)$. 
Then  
$\widehat{X} \in \aCycle^{R-1-d}_{L^1}(\PP(E^{\rho}))$. 
Thus $\acherncl_1(\overline{\OO_{E^{\rho}}(1)})^{d+1}\cdot\widehat{X}$ 
belongs to $\aChow^{R}_{L^1}(\PP (E^{\rho}))_{\QQ}$. Hence, 
\[
\pi_*\left( 
\acherncl_1(\overline{\OO_{E^{\rho}}(1)})^{d+1}\cdot\widehat{X} \right) 
\in \aChow^{1}_{L^1}(Y)_{\QQ}.
\]
Let us consider elementary properties of
$\pi_*\left( 
\acherncl_1(\overline{\OO_{E^{\rho}}(1)})^{d+1}\cdot\widehat{X} \right)$.

\begin{Proposition}
\label{prop:when:Bost:divisor:smooth}
Let $X = \sum_{k=1}^l a_k X_k$ be the irreducible decomposition
of $X$ as cycles.
Let $\phi_k : \tilde{X}_k \to X_k$
be a generic resolution of singularities of $X_k$
for each $k$, 
i.e., 
$\phi_k$ is a proper birational morphism 
such that $(\tilde{X}_k)_{\QQ}$ is smooth over $\QQ$. 
Let 
$i_k : X_k \hookrightarrow \PP(E^{\rho})$ be the inclusion map and  
$j_k : \tilde{X}_k \to \PP(E^{\rho})$ the composition map $i_k \cdot \phi_k$. 
Also we let 
$f_k : X_k \to Y$ be the composition map $\pi \cdot i_k$ and 
$\tilde{f}_k : \tilde{X}_k \to Y$ the composition map $\pi \cdot j_k$.
Let $Y_0$ be the maximal open set of $Y$ such that
$\tilde{f}_k$ is smooth over there for every $k$.
Then, we have the following.
\begin{enumerate}
\renewcommand{\labelenumi}{(\arabic{enumi})}
\item
${\displaystyle
 \pi_* \left(\acherncl_1(\overline{\OO_{E^{\rho}}(1)})^{d+1} 
             \cdot \widehat{X} \right) =
 \sum_{k=1}^l a_k \tilde{f}_k{}_*
 (\acherncl_1(j_k^* \overline{\OO_{E^{\rho}}(1)})^{d+1}).
}$ \\
In particular,
$\pi_*\left( 
\acherncl_1(\overline{\OO_{E^{\rho}}(1)})^{d+1}\cdot\widehat{X} \right)$ 
is independent of 
the choice of normalized Green forms $g_X$ for $X$, and
$\pi_* \left(\acherncl_1(\overline{\OO_{E^{\rho}}(1)})^{d+1} 
             \cdot \widehat{X} \right)
\in \aBChow^1(Y; Y_0(\CC))$.

\item
Let $y$ be a closed point of $(Y_0)_{\QQ}$, and
$\Gamma'$ the closure of $\{ y \}$ in $Y$.
Here we choose $g_X$ as in 
Proposition~\rom{\ref{prop:normalized:Green:form}}.
Then, there is a representative $(Z, g_Z)$ of
$\acherncl_1(\overline{\OO_{E^{\rho}}(1)})^{d+1} \cdot \widehat{X}$
such that $\pi^{-1}(\Gamma')$ and $Z$ intersect properly, and
$\rest{g_Z}{\pi^{-1}(z)}$ is locally integrable for each 
$z \in O_{\Gal(\overline{\QQ}/\QQ)}(y)$.
\end{enumerate}
\end{Proposition}

\Proof
We may assume that $X$ is reduced and irreducible, so that 
we will omit index $k$ in the following. 

(1)
Let $g_X$ be a $\Omega_{\rho}$-normalized Green form for $X$.
Then, by virtue of Proposition~\ref{prop:formula:restriction:intersection},
\[
\acherncl_1(\overline{\OO_{E^{\rho}}(1)})^{d+1} \cdot \widehat{X}
= j_* \left( \acherncl_1(j^* \overline{\OO_{E^{\rho}}(1)})^{d+1} \right) 
+ a(\Omega_{\rho}^{d+1} \wedge [g_X]).
\]
Therefore, since $\pi_* (g_X \wedge \Omega_{\rho}^{d+1}) = 0$ by 
the definition of $g_X$, 
we get
\[
\pi_* \left(\acherncl_1(\overline{\OO_{E^{\rho}}(1)})^{d+1} 
            \cdot \widehat{X} \right) =
\pi_* j_* \left( \acherncl_1(j^* \overline{\OO_{E^{\rho}}(1)})^{d+1} \right)
= \tilde{f}_* \left( \acherncl_1(j^* \overline{\OO_{E^{\rho}}(1)})^{d+1} \right).
\]

\medskip
(2)
First of all, we need the following lemma.

\begin{Lemma}
\label{lem:intersect:proper}
Let $T$ be a quasi-projective integral scheme over $\ZZ$, 
$L_1, \ldots, L_n$ line bundles on $T$, and
$\Gamma$ a cycle on $T$.
Then, there is a cycle $Z$ on $T$ such that
$Z$ is rationally equivalent to $c_1(L_1) \cdots c_1(L_n)$,
and that $Z$ and $\Gamma$ intersect properly.
\end{Lemma}

\Proof
We prove this lemma by induction on $n$.
First, let us consider the case $n=1$.
Let $\Gamma = \sum_{i=1}^r a_i \Gamma_i$ be the irreducible
decomposition as cycles.
Let $\gamma_i$ be a closed point of 
$\Gamma_i \setminus \bigcup_{j \not= i} \Gamma_j$, and
$m_i$ the maximal ideal at $\gamma_i$.
Let $H$ be an ample line bundle on $X$.
Choose a sufficiently large integer $N$ such that
\[
H^1(T, H^{\otimes N} \otimes m_1 \otimes \cdots \otimes m_r) =
H^1(T, H^{\otimes N} \otimes L_1 \otimes m_1 \otimes \cdots \otimes m_r)
= 0.
\]
Then, the natural homomorphisms
\[
H^0(T, H^{\otimes N}) \to \bigoplus_{i=1}^r H^{\otimes N} \otimes \kappa(\gamma_i)
\quad\text{and}\quad
H^0(T, H^{\otimes N} \otimes L_1) \to \bigoplus_{i=1}^r 
H^{\otimes N} \otimes L_1 \otimes \kappa(\gamma_i)
\]
are surjective. 
Thus, there are sections $s_1 \in H^0(T, H^{\otimes N})$ and
$s_2 \in H^0(T, H^{\otimes N} \otimes L_1)$ such that
$s_1(\gamma_i) \not= 0$ and $s_2(\gamma_i) \not= 0$ for all $i$.
Then, $\zero(s_2) - \zero(s_1) \sim c_1(L_1)$, and
$\zero(s_2) - \zero(s_1)$ and $\Gamma$ intersect properly.

Next we assume $n > 1$. Then, by hypothesis of induction,
there is a cycle $Z'$ such that
$Z' \sim c_1(L_1) \cdots c_1(L_{n-1})$, and
$Z'$ and $\Gamma$ intersect properly.
Let $Z' = \sum_{j} b_j T_j$ be the decomposition as cycles.
We set $\Gamma_j = (T_j \cap \Supp(\Gamma))_{red}$.
Then, using the case $n=1$,
there is a cycle $Z_j$ such that $Z_j \sim c_1(\rest{L_n}{T_j})$, and
$Z_j$ and $\Gamma_j$ intersect properly.
Thus, if we set $Z = \sum_j b_j Z_j$, then
$Z \sim c_1(L_1) \cdots c_1(L_n)$, and
$Z$ and $\Gamma$ intersect properly.
\QED

Let us go back to the proof of (2) of 
Proposition~\ref{prop:when:Bost:divisor:smooth}.
By virtue of Lemma~\ref{lem:intersect:proper},
there is a cycle $Z$ on $X$ such that
$Z \sim c_1\left(i^* \OO_{E^{\rho}}(1) \right)^{d+1}$, and that $Z$ and
$f^{-1}(\Gamma')$ intersect properly.
Then, $Z \sim c_1(\OO_{E^{\rho}}(1))^{d+1} \cdot X$, and
$Z$ and $\pi^{-1}(\Gamma')$ intersect properly.
Let $\phi_X$ be a Green form of logarithmic type for $X$.
Then, since
\[
\acherncl_1(\overline{\OO_{E^{\rho}}(1)})^{d+1} \cdot (X, \phi_X) 
\in \aChow^{R}(\PP(E^{\rho})),
\]
there is a Green form $\phi_Z$ of logarithmic type for $Z$ such that
$(Z, \phi_Z)$ is a representative of 
$\acherncl_1(\overline{\OO_{E^{\rho}}(1)})^{d+1} \cdot (X, \phi_X)$.
Thus, if we set
\[
g_Z = \phi_Z + c_1(\overline{\OO_{E^{\rho}}(1)})^{d+1} \wedge (g_X - \phi_X),
\]
then $(Z, g_Z)$ is a representative of
$\acherncl_1(\overline{\OO_{E^{\rho}}(1)})^{d+1} \cdot \widehat{X}$.
Since $Z$ and $\pi^{-1}(\Gamma')$ intersect properly and
$g_X$ has the property in Proposition~\ref{prop:normalized:Green:form},
we can easily see that
$g_Z$ is locally integrable along $\pi^{-1}(z)$ for each 
$z \in O_{\Gal(\overline{\QQ}/\QQ)}(y)$.
\QED

Here we recall some elementary results of Chow forms. 
Details can be found in \cite{Bo}. 
Consider the incidence subscheme $\Gamma$ in the product
\[
\PP(E^{\rho}) \times_{Y} \PP(E^{\rho}{}^{\lor})^{d+1}
= \PP(E^{\rho}) \times_{Y} \PP(E^{\rho}{}^{\lor}) 
\times_{Y} \cdots \times_{Y} \PP(E^{\rho}{}^{\lor}).
\]
Let $\imath : \Gamma \to \PP(E^{\rho})$ and $\jmath : 
\Gamma \to \PP(E^{\rho}{}^{\lor})^{d+1}$ 
be projection maps. The Chow divisor $\Ch(X)$ of $X$ is defined by 
\[
\Ch(X) = \jmath_* \imath^* (X).
\]
The following facts are well-known: 
\begin{enumerate}
 \item $\Ch(X)$ is an effective cycle of codimension $1$ 
       in $\PP(E^{\rho}{}^{\lor})^{d+1}$;
 \item $\Ch(X)$ is flat over $Y$;
 \item For any $y \in Y$, 
       $\Ch(X)_y$ is a divisor of type $(\delta,\delta,\cdots,\delta)$ 
       in $\PP(E^{\rho}{}^{\lor})^{d+1}_{y}$.
\end{enumerate}
Let $p : \PP(E^{\rho}{}^{\lor})^{d+1} \to Y$ be the canonical morphism, 
and $p_i : \PP(E^{\rho}{}^{\lor})^{d+1} 
\to \PP(E^{\rho}{}^{\lor})$ the projection to the $i$-th factor. 
Then, by the above properties, 
there is a line bundle $L$ on $Y$ and a section $\Phi_{X}$ of 
\[
H^0 \left(\PP(E^{\rho}{}^{\lor})^{d+1},
p^*(L)\otimes\bigotimes_{i=1}^{d+1}p_i^* 
\OO_{E^{\rho}{}^{\lor}}(\delta) \right)
\]
such that $\zero(\Phi_{X}) = \Ch(X)$. Since 
\[
p_* \left(
    p^* (L) \otimes \bigotimes_{i=1}^{d+1} p_i^* 
\OO_{E^{\rho}{}^{\lor}}(\delta)
    \right)
=
L \otimes (\Sym^{\delta}(E^{\rho}))^{{}\otimes {d+1}},
\]
we may view $\Phi_{X}$ as an element of 
$H^0(Y,L \otimes (\Sym^{\delta}(E^{\rho}))^{{}\otimes {d+1}})$. 
We say $\Phi_{X}$ is a {\em Chow form} of $X$. 
Clearly $\Phi_{X}$ is uniquely determined up to $H^0(Y,\OO_{Y}^{\times})$. 

As in \cite[Proposition 1.2 and its remark]{Bo}, we have 
\[
 \cherncl_1(L) 
= \pi_* \left(\cherncl_{1}(\OO_{E^{\rho}}(1))^{d+1}\cdot {X}\right).
\]
We give a generalized metric $h_{L}$ on $L$ so that 
$\overline{L} = (L,h_L)$ satisfies the equality
\addtocounter{Theorem}{1} 
\begin{equation}
\label{eqn:metric:L}
 \acherncl_1(\overline{L}) 
= \pi_* \left(\acherncl_{1}(\overline{\OO_{E^{\rho}}(1)})^
                             {d+1}\cdot\widehat{X}\right) 
\end{equation}
in $\aChow^1_{L^1}(Y)$.

Note that we can also give a metric $L$ by the equation 
\[
\OO_{\PP(E^{\rho}{}^{\lor})^{d+1}}(\Ch(X)) 
= p^*(L)\otimes\bigotimes_{i=1}^{d+1}p_i^* \OO_{E^{\lor}}(\delta)
\]
and by suitably metrizing other terms, 
as is implicitly done in \cite[1.5]{Zh}. 
We do not however pursue this here. 

\renewcommand{\theClaim}{\arabic{section}.\arabic{subsection}.\arabic{Theorem}.\arabic{Claim}}
\renewcommand{\theequation}{\arabic{section}.\arabic{subsection}.\arabic{Theorem}.\arabic{Claim}}

\subsection{Restriction of Chow forms on fibers}
\label{subsec:pullback}
\setcounter{Theorem}{0}
In this section we will consider the restriction of Chow forms
on fibers. 

Let $Y,\overline{E},\rho,X$ 
be as in \S \ref{subsec:Chow:forms:and:their:metrics}. 
Let $y$ be a closed point of $Y_{\QQ}$.
Let $\Gamma'$ be the closure of $\{ y \}$ in $Y$, and
$\Gamma$ the normalization of $\Gamma'$. 
Let $f : \Gamma \to Y$ be the natural morphism. 
We set $E_{\Gamma} = f^*(E)$ and 
$\overline{E_{\Gamma}} = (E_{\Gamma},f^*(h))$. 
Also we put $(E^{\rho})_{\Gamma} = f^*(E^{\rho})$ and 
$\overline{(E^{\rho})_{\Gamma}} = ((E^{\rho})_{\Gamma}, f^*(h^{\rho}))$. 
Then $(\overline{E_{\Gamma}})^{\rho}$ is equal to 
$\overline{(E^{\rho})_{\Gamma}}$, so that
we denote $(E^{\rho})_{\Gamma}$ by $E^{\rho}_{\Gamma}$, and
$\overline{(E^{\rho})_{\Gamma}}$ by $\overline{E^{\rho}_{\Gamma}}$.
Considering the following fiber product  
\begin{equation*}
 \begin{CD}
   \PP(E^{\rho}_{\Gamma})  @>{f '}>> \PP(E^{\rho})   \\
   @VV{\pi'}V           @VV{\pi}V      \\
   \Gamma         @>{f}>>   Y
 \end{CD}
\end{equation*}
we set $X_{\Gamma} = {f'}^*(X)$. 
Then $X_{\Gamma}$ is flat over $\Gamma$ 
with the relative dimension $d$ 
and the degree $\delta$ on the generic fiber. 
For this quadruplet $(\Gamma,\overline{E_{\Gamma}},\rho,X_{\Gamma})$ 
in place of the quadruplet $(Y,\overline{E},\rho,X)$, 
we can define in the same way the Hermitian line bundle 
$\overline{\OO_{E^{\rho}_{\Gamma}}(1)}$ 
on $\PP(E^{\rho}_{\Gamma})$, 
an arithmetic $L^1$-divisor 
$\widehat{X_{\Gamma}} = (X_{\Gamma}, g_{X_{\Gamma}})$ 
on $\PP(E^{\rho}_{\Gamma})$ 
and the arithmetic $L^1$-divisor 
$\pi '_* \left(\acherncl_{1}(
\overline{\OO_{E^{\rho}_{\Gamma}}(1)})^{d+1}\cdot 
\widehat{X_{\Gamma}}\right)$ 
on $\Gamma$.
Further, we have 
$\overline{L_{\Gamma}} = (L_{\Gamma}, h_{\Gamma})$
satisfying 
\[
 \acherncl_1(\overline{L_{\Gamma}}) 
 = \pi '_* \left(\acherncl_{1}(\overline{\OO_{E^{\rho}_{\Gamma}}(1)})^
                             {d+1}\cdot\widehat{X_{\Gamma}}\right) 
\]
in $\aChow_{L^1}^1(\Gamma)$.
We also have $\Ch(X_{\Gamma})$. Moreover, 
letting $p_i ' : \PP((E^{\rho}_{\Gamma})^{\lor})^{d+1} 
\to \PP((E^{\rho}_{\Gamma})^{\lor})$ 
be the $i$-th projection, 
there is a section $\Phi_{X_{\Gamma}}$ of 
\[
 H^0 \left( \PP((E^{\rho}_{\Gamma})^{\lor})^{d+1},
{p '}^*(L_{\Gamma})\otimes
\bigotimes_{i=1}^{d+1} {p_i'}^* 
\OO_{(E^{\rho}_{\Gamma})^{\lor}}(\delta) \right) 
 =
H^0 \left(\Gamma,
L_{\Gamma} \otimes (\Sym^{\delta}
((E^{\rho}_{\Gamma})^{\lor}))^{\otimes d+1} \right),
\] 
such that $\zero(\Phi_{X_{\Gamma}}) = \Ch(X_{\Gamma})$.  
 
Let us consider the following fiber product,  
\begin{equation*}
 \begin{CD}
   \PP((E^{\rho}_{\Gamma})^{\lor})^{d+1} 
 @>{g '} >> \PP((E^{\rho})^{\lor})^{d+1}  \\
   @VV{p '}V
   @VV{p}V \\
   \Gamma @>{f}>> Y
 \end{CD}
\end{equation*}

Then, we have the following proposition.

\begin{Proposition}
\label{proposition:pullback}
\begin{enumerate}
\renewcommand{\labelenumi}{(\roman{enumi})}
\item
${g'}^* \Ch(X) = \Ch(X_{\Gamma})$. 
Moreover, we can choose $\Phi_{X_{\Gamma}}$ to be $f ^* \Phi_X$. 
\item
Let $X_1,\cdots,X_l$ 
be the irreducible components of $X_{red}$. 
Assume that, for every $1 \leq i \leq l$, 
there is a generic resolution of singularities 
$\phi_i : \tilde{X_i} \to X_i$ 
such that the induced map 
$\tilde{X_i} \to Y$ 
is smooth over $y$ for every $i$. 
Then the equality 
\begin{align*}
f^* \pi_* \left(\acherncl_{1}(
\overline{\OO_{E^{\rho}}(1)})^{d+1} \cdot \widehat{X} \right) 
& = \pi '_* f'{}^* \left(\acherncl_{1}(
\overline{\OO_{E^{\rho}}(1)})^{d+1} \cdot \widehat{X} \right) \\
& = \pi '_* \left(\acherncl_{1}(
\overline{\OO_{E^{\rho}_{\Gamma}}(1)})^{d+1}
             \cdot \widehat{X_{\Gamma}} \right). 
\end{align*}
holds in $\aChow_{L^1}^1(\Gamma)$. 
In other words, $f^* (L,h_L) = (L_{\Gamma}, h_{L_{\Gamma}})$.
\end{enumerate}
\end{Proposition}

\Proof
(i) If $f$ is flat, 
then this follows from the base change theorem. 
In the case $f$ is not flat, we refer to the remark \cite[4.3.2(i)]{BGS},
or we can easily see this using Appendix~A. 

(ii) 
We take $g_X$ as in Proposition~\ref{prop:normalized:Green:form}.
Let $\alpha = \acherncl_{1}(
\overline{\OO_{E^{\rho}}(1)})^{d+1} \cdot \widehat{X}
\in \aChow_{L^1}^1(Y)$. 
By Proposition~\ref{prop:when:Bost:divisor:smooth}, 
we can take a representative 
$(Z,g_Z)$ of $\alpha$ such that
$Z$ and $\pi^{-1}(\Gamma')$ intersect properly, and
$g_Z$ is locally integrable along $\pi^{-1}(w)$ for all 
$w \in O_{\Gal(\overline{\QQ}/\QQ)}(y)$.
Now we have 
\begin{align*}
f^* \pi_*(\alpha) 
& = f^* (\pi_*Z, [\pi_* g_Z]) \\
& = \left( f^* \pi_*Z, 
    \sum_{w \in O_{\Gal(\overline{\QQ}/\QQ)}(y)} 
    \left( \int_{\pi^{-1} (w)} g_Z \right) \cdot w \right).
\end{align*}
On the other hand, we have 
\begin{equation*}
\pi '_* f'{}^*(\alpha)
= \left( \pi '_* f'{}^* Z, 
    \sum_{w \in O_{\Gal(\overline{\QQ}/\QQ)}(y)} 
    \left( \int_{\pi^{-1} (w)} g_Z \right) \cdot w \right).
\end{equation*}
Moreover, by Appendix A, $f^* \pi_*Z$ is equal to  $\pi '_* f'{}^* Z$ 
as a cycle. Thus we have proven the first equality. 

Now we will prove the second equality.
Let $\phi_X$ be a Green form of logarithmic type for $X$.
Since 
\[
{f'}^* : \bigoplus_{i \geq 0} \aChow^i(\PP(E^{\rho})) \to
\bigoplus_{i \geq 0} \aChow^i(\PP(E^{\rho}_{\Gamma}))
\]
is a homomorphism of rings (cf. \cite[5) of Theorem in 4.4.3]{GSArInt}).
Thus,
\[
{f'}^* \left(\acherncl_{1}(
\overline{\OO_{E^{\rho}}(1)})^{d+1} \cdot (X, \phi_X) \right) =
{f'}^* \acherncl_{1}(
\overline{\OO_{E^{\rho}}(1)})^{d+1} \cdot {f'}^*(X, \phi_X).
\]
Therefore, since we take $g_X$ as in 
Proposition~\ref{prop:normalized:Green:form}, we can see
\begin{align*}
{f'}^* \left(\acherncl_{1}(
\overline{\OO_{E^{\rho}}(1)})^{d+1} \cdot \widehat{X} \right) & =
{f'}^* \left(\acherncl_{1}(
\overline{\OO_{E^{\rho}}(1)})^{d+1} \cdot ((X, \phi_X) + (0, g_X - \phi_X)) \right) \\
& =
{f'}^* \left(\acherncl_{1}(
\overline{\OO_{E^{\rho}}(1)})^{d+1} \cdot (X, \phi_X) \right) +
a \left( 
{f'}^*(c_{1}(\overline{\OO_{E^{\rho}}(1)})^{d+1} \wedge (g_X - \phi_X))\right) \\
& =
{f'}^* \acherncl_{1}(
\overline{\OO_{E^{\rho}}(1)})^{d+1} \cdot
{f'}^*(X, \phi_X) + 
a \left(
{f'}^* (c_{1}(\overline{\OO_{E^{\rho}}(1)})^{d+1} \wedge (g_X - \phi_X)) \right)\\
& =
{f'}^* \acherncl_{1}(
\overline{\OO_{E^{\rho}}(1)})^{d+1} \cdot
{f'}^*(X, g_X)
\end{align*}
Moreover, as pointed out in Remark~\ref{rem:norm:Green:general:fiber},
${f'}^* g_X$ is a normalized Green form for ${f'}^* X$. 
Thus we have got the second equality. 
\QED

\renewcommand{\theClaim}{\arabic{section}.\arabic{subsection}.\arabic{Theorem}.\arabic{Claim}}
\renewcommand{\theequation}{\arabic{section}.\arabic{subsection}.\arabic{Theorem}.\arabic{Claim}}

\subsection{Chow stability and field extensions}
\label{subsec:Chow:stability:field:ext}
\setcounter{Theorem}{0}

Let $\rho : \bfGL_r \to \bfGL_R$ be 
a group scheme morphism of degree $k$  
commuting with the transposed morphism. 
Let $S$ be a ring (commutative, with the multiplicative identity). 
For a positive integer $\delta$ and $d$, we consider  
$\Sym^{\delta} (S^R)^{\otimes d+1}$. 
Then through the induced group homomorphism 
$\rho(S) : \bfGL_r(S) \to \bfGL_R(S)$, 
$\bfGL_r(S)$ (or $\bfSL_r(S)$) 
acts on 
$\Sym^{\delta} (S^R)^{\otimes d+1}$. 

\begin{Proposition}
\label{prop:stability:and:det}
Let $K$ be an infinite field 
and $L$ an extension field of K. 
Let $P$ be a homogeneous polynomial of degree $e$ on 
$\Sym^{\delta} (L^R)^{\otimes d+1}$, i.e., 
$P \in \Sym^e(\Sym^{\delta} (L^R)^{\otimes d+1}{}^{\lor})$. 
Then we have the following. 
\begin{enumerate}
\renewcommand{\labelenumi}{(\roman{enumi})}
\item
$P$ is $\bfSL_r(K)$-invariant 
if and only if 
$P$ is $\bfSL_r(L)$-invariant. 
\item
If P is $\bfSL_r(K)$-invariant, then
\[
P(v^{\sigma}) ^r = (\det \sigma) ^{ek(d+1)\delta} P(v)^r
\]
in $L$ 
for all $ v \in \Sym^{\delta} (L^R)^{\otimes d+1}$ 
and $\sigma \in \bfGL_r(L)$.
\end{enumerate}
\end{Proposition}

\proof
(i) We only need to prove the `only if' part. 
Let $S_L(P) = \{ \sigma \in \bfSL_r(L) \mid P^{\sigma} = P \}$ 
be the stabilizer of $P$ in $\bfSL_r(L)$. 
$S_L(P)$ is a closed algebraic set of $\bfSL_r(L)$ 
and contains $\bfSL_r(K)$. 
Since $\bfSL_r(K)$ is Zariski dense in $\bfSL_r(L)$, 
$S_L(P)$ must coincide with $\bfSL_r(L)$.  

(ii) Let $M$ be an extension field 
of $L$ such that 
it has an $r$-th root $\xi$ of $\det \sigma$. 
If $\sigma '$ is defined by $\sigma = \xi\sigma '$, 
then $\sigma ' \in \bfSL_r(M)$. 
Since $P$ is $\bfSL_r(M)$-invariant by (i), 
we find 
\begin{align*}
P(v^{\sigma})^r  & = P\left( (\rho(\CC)(\xi\sigma ')) \cdot v\right) ^r
= P\left( (\xi^k \rho(\CC)(\sigma ')) \cdot v \right) ^r \\
 & = \xi^{rek(d+1)\delta} P(v)^r
= (\det \sigma) ^{ek(d+1)\delta} P(v)^r.
\end{align*}
in $M$. 
Hence the equality holds in $L$ because both sides belong to $L$. 
\QED

\begin{Remark}
More strongly, we can show that, 
for any integral domain $S$ of characteristic zero, 
if $P \in \Sym^e(\Sym^{\delta} (S^R)^{\otimes d+1}{}^{\lor})$ 
is $\bfSL_r(\ZZ)$-invariant,
then $P$ is $\bfSL_r(S)$-invariant. 
\end{Remark}

Now let $Y$, $\overline{E}$, $\rho$ and $X$ 
be as in \S \ref{subsec:Chow:forms:and:their:metrics}. 
Recall that for a closed point $y$ of $Y_{\QQ}$, 
$\Ch(X)_y$ is a divisor of type $(\delta,\delta,\cdots,\delta)$ 
in $\PP(E^{\rho}{}^{\lor})^{d+1}_{y}$. 
Hence the Chow form restricted on $y$, i.e., 
$\Phi_X \vert_y = \Phi_{X_y}$ is an element of 
$\Sym^{\delta}(K^R)^{d+1}$. 
We say that $X_y$ is {\em Chow semi-stable} 
if $\Phi_{X_y} \in \Sym^{\delta}(K^R)^{d+1}$ 
is semi-stable under the action of $\bfSL_r(K)$, where
$K$ is the residue field of $y$. 

\begin{Lemma}
\label{lemma:stability:generators:over:Z}
There are a positive integer $e$ 
and 
$\bfSL_r(\QQ)$-invariant homogeneous polynomials 
$P_1,\cdots,P_l \in \Sym^{e}(\Sym^{\delta}(\ZZ^R)^{d+1}{}^{\lor})$, 
which depend only on $\rho$, $d$ and $\delta$, 
with the following property.  
For any closed points $y \in Y_{\QQ}$, 
if $X_y$ is Chow semistable, then
there is a $P_i$ such that 
$P_i (\Phi_{X_y}) \neq 0$.
\end{Lemma} 

\proof
$\bfSL_r(\QQ)$ 
acts linearly on $\Sym^{\delta}(\QQ^R)^{d+1}$. 
Since $\bfSL_r(\QQ)$ is a reductive group, 
we can take $\bfSL_r(\QQ)$-invariant 
homogeneous polynomials $Q_1,\cdots,Q_l$ 
such that they form generators of the algebra of 
$\bfSL_r(\QQ)$ invariant polynomials 
on $\Sym^{\delta}(\QQ^R)^{d+1}$. 
By clearing the denominators, 
we may assume that $Q_1,\cdots,Q_l$ is defined over $\ZZ$. 
Let $e_i$ be the degree of $Q_i$ for $i = 1,\cdots,l$. 
We take a positive integer $e$ such that $e_i \vert e$ 
for $i = 1,\cdots,l$. 
We set $P_i = Q_i^{e/e_i}$.

Let us check that $P_i$'s have the desired property.
Since $X_y$ is Chow semistable, there is a 
$\bfSL_r(K)$-invariant homogeneous
polynomial $F$ on $\Sym^{\delta}(K^R)^{d+1}$
with $F(\Phi_{X_y}) \not= 0$, where
$K$ is the residue field of $y$.
Let us choose $\alpha_1, \ldots, \alpha_n \in K$ and
homogeneous polynomials $F_1, \ldots, F_n$
over $\QQ$ such that
$F = \alpha_1 F_1 + \cdots + \alpha_n F_n$ and that
$\alpha_1, \ldots, \alpha_n$ are linearly independent over $\QQ$.
Here, for $\sigma \in \bfSL_r(\QQ)$,
\[
F^{\sigma} = \alpha_1 F_1^{\sigma} + \cdots + \alpha_n F_{n}^{\sigma}
\]
and $F_{i}^{\sigma}$'s are homogeneous polynomials
over $\QQ$. Thus, we can see that $F_i$'s are $\bfSL_r(\QQ)$-invariant.
Moreover, since
\[
F(\Phi_{X_y}) = \alpha_1 F_1(\Phi_{X_y}) + \cdots + \alpha_n F_{n}(\Phi_{X_y}),
\]
there is $F_i$ with $F_i(\Phi_{X_y}) \not= 0$.
On the other hand, $F_i$ is an element of $\QQ[Q_1, \ldots, Q_l]$.
Thus, we can find $Q_j$ with $Q_j(\Phi_{X_y}) \not= 0$, namely
$P_j(\Phi_{X_y}) \not= 0$.
\QED

\medskip
\section{Semi-stability and positiveness in a relative case}
\label{section:semistability:positiveness:relative}

\subsection{Cornalba-Harris-Bost's inequality in a relative case}
\label{subsec:CHB:inequality}
\setcounter{Theorem}{0}
Let $Y$ be an arithmetic variety and 
$\overline{E} = (E,h)$ a Hermitian vector bundle of rank $r$ on $Y$. 
Let $\rho : \bfGL_r \to \bfGL_R$ 
be a group scheme morphism of degree $k$ 
commuting with the transposed morphism. 

Before we prove the relative Cornalba-Harris-Bost's inequality, 
we need three lemmas. 

\begin{Lemma}
\label{lemma:P:to:sheafhom}
Let $L$ be a line bundle on $Y$. 
Let $P \in \Sym^e(\Sym^{\delta}(\ZZ^R)^{d+1}{}^{\lor}) \backslash \{0\}$ 
and suppose that $P$ is $\bfSL_r(\QQ)$-invariant. 
Then there is a polynomial map of sheaves
\[
\varphi_P : L \otimes \Sym^{\delta} (E^{\rho})^{\otimes d+1} 
\to L^{\otimes er} \otimes (\det E)^{\otimes ek(d+1)\delta}
\]
given by $P^r$, 
namely, 
$\varphi_P$ is locally defined by the evaluation in terms of $P^r$.
\end{Lemma}

\Proof
Let $U$ be a Zariski open set, and
$\phi : \rest{E}{U}
\overset{\sim}{\longrightarrow} \OO_{U}^{\oplus n}$ 
and
$\psi : \rest{L}{U}
\overset{\sim}{\longrightarrow} \OO_{U}$ 
local trivializations of $E$ and $L$ respectively. 
Then, by the construction of $E^{\rho}$, we have 
\[
\phi_{\rho,\delta,d}: 
\rest{\left( \Sym^{\delta} \left(E^{\rho}\right)^{\otimes d+1} \right)}{U}
\overset{\sim}{\longrightarrow} 
\Sym^{\delta} \left(\OO_{U}^{\oplus R}\right)^{\otimes d+1}.
\]
Thus we get 
\[
\psi \otimes \phi_{\rho,\delta,d}: 
\rest{\left( L \otimes \Sym^{\delta} \left(E^{\rho}\right)^{\otimes d+1} 
\right)}{U}
\overset{\sim}{\longrightarrow} 
\Sym^{\delta} \left(\OO_{U}^{\oplus R}\right)^{\otimes d+1}.
\]
Here, we define 
\[
\rest{\varphi_P}{U} : 
\rest{\left( L \otimes \Sym^{\delta} (E^{\rho})^{\otimes d+1} \right)}{U}
\to 
\rest{\left( L^{\otimes er} \otimes (\det E)^{\otimes ek(d+1)\delta} \right)}{U}
\]
such that the following diagram is  commutative.
\[
\begin{CD}
\rest{\left( L \otimes \Sym^{\delta} \left(E^{\rho}\right)^{\otimes d+1} \right)}{U} 
@>{\psi \otimes  \phi_{\rho,\delta,d}}>> 
\Sym^{\delta} \left(\OO_{U}^{\oplus R}\right)^{\otimes d+1} \\
@V{\rest{\varphi_P}{U}}VV @VV{P^r}V \\
\rest{\left( L^{\otimes er} \otimes (\det E)^{\otimes ek(d+1)\delta} \right)}{U}
@>{\psi^{er} \otimes \det(\phi)^{ek(d+1)\delta}}>>
\OO_U,
\end{CD}
\]
where $P^r$ is the map given by the evaluation in terms of the polynomial $P^r$.
In order to see that the local $\rest{\varphi_P}{U}$
glues together on $Y$,
it is sufficient to show that
$\rest{\varphi_P}{U}$ does not depend on the choice of
local trivializations $\phi$ and $\psi$.
Let $\phi' : \rest{E}{U}
\overset{\sim}{\longrightarrow} \OO_{U}^{\oplus n}$ 
and
$\psi' : \rest{L}{U}
\overset{\sim}{\longrightarrow} \OO_{U}$ 
be another local trivializations.
In the same way, we have the following commutative diagram.
\[
\begin{CD}
\rest{\left( L \otimes \Sym^{\delta} \left(E^{\rho}\right)^{\otimes d+1} \right)}{U} 
@>{\psi' \otimes  \phi'_{\rho,\delta,d}}>> 
\Sym^{\delta} \left(\OO_{U}^{\oplus R}\right)^{\otimes d+1} \\
@V{\rest{\varphi'_P}{U}}VV @VV{P^r}V \\
\rest{\left( L^{\otimes er} \otimes (\det E)^{\otimes ek(d+1)\delta} \right)}{U}
@>{{\psi'}^{er} \otimes \det(\phi')^{ek(d+1)\delta}}>>
\OO_U
\end{CD}
\]
We set the transition functions 
$g = \phi \cdot(\phi')^{-1}$ and 
$h = \psi \cdot(\psi')^{-1}$.
Then by a straightforward calculation using 
(ii) of Proposition~\ref{prop:stability:and:det}, 
we get, on $U$, 
\[
P^r \cdot (\psi \otimes  \phi_{\rho,\delta,d})
= h^{re}\det(g)^{ek (d+1)\delta} 
P^r \cdot (\psi' \otimes  \phi'_{\rho,\delta,d}),
\]
which implies
\[
\left({\psi}^{er} \otimes \det(\phi)^{ek(d+1)\delta}\right) \cdot 
\left(\rest{\varphi_P}{U}\right) = 
h^{re}\det(g)^{ek (d+1)\delta} 
\left({\psi'}^{er} \otimes \det(\phi')^{ek(d+1)\delta}\right) \cdot 
\left(\rest{\varphi'_P}{U}\right).
\]
Here note that
\[
h^{re}\det(g)^{ek (d+1)\delta}  
= 
\left({\psi}^{er} \otimes \det(\phi)^{ek(d+1)\delta}\right) \cdot
\left({\psi'}^{er} \otimes \det(\phi')^{ek(d+1)\delta}\right)^{-1}.
\]
Thus, we obtain $\rest{\varphi_P}{U} = \rest{\varphi'_P}{U}$.
\QED

Suppose now $L$ is given a generalized metric $h_L$. 
Since both sides of 
\[
\varphi_P : L \otimes \Sym^{\delta} (E^{\rho})^{\otimes d+1} 
\to L^{\otimes er} \otimes (\det E)^{\otimes ek(d+1)\delta}
\]
in the lemma above are then equipped with metrics, 
we can consider the norm of $\varphi_P$. 
Before evaluating the norm of $\varphi_P$, 
we define the norm of $P$ as follows; 
We first define the metric $\Vert \cdot \Vert_{can}$ 
on $\Sym^{\delta}(\CC^{n})^{\otimes d+1}$ 
induced from the usual Hermitian metric on $\CC$; 
We then define $\VERT P \VERT$ by 
\[
\VERT P \VERT = 
\sup_{v \in \Sym^{\delta}(\CC^{n})^{\otimes d+1} \setminus \{0\}}
\frac{\vert P(v) \vert}{\Vert v \Vert_{can}^e},
\]
where $P$ is regarded as an element of  
$\Sym^e((\Sym^{\delta}(\CC^{m}))^{\otimes d+1})^{\lor})$. 

\begin{Lemma}
\label{lemma:norm:of:morphism:P}
For any section $s \in H^0(Y,L \otimes (\Sym^{\delta}(E^{\rho}))^{\otimes d+1})$ 
and any complex point $y \in Y(\CC)$ 
around which $h_L \otimes (\Sym^{\delta}(h^{\rho}))^{\otimes d+1}$ is 
$C^{\infty}$, 
we have
\[
\Vert \varphi_P(s) \Vert (y) \le 
\VERT P \VERT^r \, \Vert s \Vert^{er} (y). 
\]
\end{Lemma}

\proof
By choosing bases, $\overline{E}(y)$ and in $\overline{L}(y)$ are isometric to 
$\CC^n$ and $\CC$ with the canonical metrics, respectively. 
Then, with respect to these bases, 
$\overline{E}^{\rho}$ is by its construction isomorphic to 
$\CC^R$ with the canonical metric.
Recalling that
$\varphi_P$ is given by the evaluation by $P^r$
once we fix local trivializations of $E$ and $L$, 
the desired inequality follows from the definition of 
$\VERT P \VERT$.
\QED

Now let $X$ be an effective cycle in $\PP(E^{\rho})$ 
such that $X$ is flat over $Y$ 
with the relative dimension $d$ and the degree $\delta$ on the generic fiber. 
In \S\ref{subsec:Chow:forms:and:their:metrics} 
we constructed a Chow form $\Phi_X$ of $X$, 
which is an element of 
$H^0(Y, L \otimes (\Sym^{\delta}(E^{\rho}))^{{}\otimes {d+1}})$. 
Recall that 
$L$ is given a generalized metric by \eqref{eqn:metric:L}.
For each irreducible component $X_i$ of $X_{red}$,
let $\tilde{X_i} \to X_i$ be 
a generic resolution of singularities of $X_i$.
Moreover, let $Y_0$ be the maximal open set of $Y$ such that
the induced morphism $\tilde{X}_i \to Y$ is smooth over $Y_0$
for every $i$.  

Further, we fix terminologies.
Let $T$ be a quasi-projective scheme over $\ZZ$, 
$t$ a closed point of $T_{\QQ}$, and $K$ the residue field
of $t$. By abuse use of notation, let $t : \Spec(K) \to T$
be the induced morphism by $t$.
We say $t$ is {\em extensible in $T$} if $t : \Spec(K) \to T$
extends to $\tilde{t} : \Spec(O_K) \to T$, where
$O_K$ is the ring of integers in $K$.
Note that if $T$ is projective over $\ZZ$, then
every closed point of $T_{\QQ}$ is extensible in $T$.

Let $V$ be a set, $\phi$ a non-negative function on $V$, and
$S$ a finite subset of $V$.
We define the geometric mean $\GM(\phi; S)$ of
$\phi$ over $S$ to be
\[
\GM(\phi; S) = \left( \prod_{s \in S} \phi(s) \right)^{1/\#(S)}.
\]

We will evaluate the norm of $\Phi_X$. 

\begin{Lemma}
\label{lemma:evaluation:of:norm:of:Phi} 
There is a constant $c_1 (R,d,\delta)$ 
depending only on $R,d$ and $\delta$ 
with the following property.
For any closed points $y$ of $(Y_0)_{\QQ}$
with $y$ extensible in $Y$,
\[
\GM\left(
\Vert \Phi_{X} \Vert_{
\overline{L}\otimes (\Sym^{\delta}(\overline{E}^{\rho}))^{{}\otimes {d+1}}}
;\  O_{\Gal(\overline{\QQ}/{\QQ})}(y) \right)
\leq  c_1 (R,d,\delta).
\]
\end{Lemma}

\Proof
Let $K$ be the residue field of $y$.
Let $\Gamma$ be the normalization of the closure of $\{ y \}$ in $Y$. 
Then, since $y$ is extensible in $Y$,
$\Gamma = \Spec(O_K)$.
Thus, by virtue of Proposition~\ref{proposition:pullback}, 
we may assume $Y = \Spec(O_K)$.
In this case, the estimate of the Chow form was already given 
in \cite[Proposition~1.3]{Bo} and \cite[4.3]{BGS}. 
Indeed if we let $k_L$ be the metric on $L$ such that 
\[
\Vert \Phi_{X} \Vert_
{(L,k_L) \otimes (\Sym^{\delta}(\overline{E}^{\rho}))^{{}\otimes {d+1}}} (w) 
= 1
\]
for every $w \in O_{\Gal(\overline{\QQ}/{\QQ})}(y)$, 
then 
$\widehat\deg (L,h_L) = h_{\overline{\OO_E(1)}}(X)$ 
and $\widehat\deg (L,k_L) = h_{Herm}(\Ch(X))$, in the notation of \cite{BGS}.
\QED


Now we will state a relative case of Cornalba-Harris-Bost's inequality. 

\begin{Theorem}
\label{thm:semistability:imply:average:semi-ampleness}
Let $Y$ be a regular arithmetic variety, 
$\overline{E} = (E,h)$ a Hermitian vector bundle of rank $r$ on $Y$, 
$\rho : \bfGL_r \to \bfGL_R$ a group scheme morphism of degree $k$ 
commuting with the transposed morphism. 
Let $X$ be an effective cycle in $\PP (E^{\rho})$ 
such that $X$ is flat over $Y$ 
with the relative dimension $d$ 
and degree $\delta$ on the generic fiber. 
Let $X_1,\ldots,X_l$ 
be the irreducible components of $X_{red}$, and
$\tilde{X}_i \to X_i$ a generic resolution of singularities of $X_i$.
Let $Y_0$ be the maximal open set of $Y$ such that
the induced morphism $\tilde{X}_i \to Y$ is smooth over $Y_0$ for
every $i$.
Let $(B, h_B)$ be a line bundle equipped with
a generalized metric on $Y$ given by the equality:
\[
\acherncl_1(B, h_B) = 
r \pi_* \left(\acherncl_{1}(\overline{\OO_{E^{\rho}}(1)})^{d+1}
                \cdot\widehat{X} \right) 
  + k \delta (d+1) \acherncl_{1}(\overline{E}).
\]
Then, $h_B$ is $C^{\infty}$ over $Y_0$. Moreover,
there are a positive integer
$e=e(\rho,d,\delta)$, a positive integer $l=l(\rho,d,\delta)$,
a positive constant $C=C(\rho,d,\delta)$, and
sections $s_1, \ldots, s_l \in H^0(Y, B^{\otimes e})$
with the following properties.
\begin{enumerate}
\renewcommand{\labelenumi}{(\roman{enumi})}
\item
$e$, $l$, and $C$ depend only on $\rho$, $d$, and $\delta$.

\item
For a closed point $y$ of $Y_{\QQ}$, if $X_y$ is Chow semistable,
then $s_i(y) \not= 0$ for some $i$.

\item
For all $i$ and all closed points $y$ of $(Y_0)_{\QQ}$
with $y$ extensible in $Y$,
\[
\GM\left( \left( h_B^{\otimes e} \right)(s_i, s_i);\ 
          O_{\Gal(\overline{\QQ}/\QQ)}(y)\right)
\leq C.
\]
\end{enumerate}
In particular, if we set
\[
\beta = 
e \left(
r \pi_* \left(\acherncl_{1}(\overline{\OO_{E^{\rho}}(1)})^{d+1}
                \cdot\widehat{X} \right) 
  + k \delta (d+1) \acherncl_{1}(\overline{E})\right)
+ a (\log C),
\]
then, for any closed point $y \in (Y_0)_{\QQ}$ with
$X_y$ Chow semistable,
there is a representative $(D, g)$ of $\beta$
such that $D$ is effective, $y \not\in D$, and that
\[
\sum_{w \in O_{\Gal(\overline{\QQ}/\QQ)}(z)} g(w) \geq 0
\]
for all $z \in (Y_0)_{\QQ}$ with $z$ extensible in $Y$.
\end{Theorem}
 
Note that if $\rho$ is the identity morphism, 
then, by the proof below, 
$C(\rho,d,\delta)$ is depending only on $r,d,\delta$. \\

\Proof
First of all, by
Proposition~\ref{prop:when:Bost:divisor:smooth},
\[
 r \pi_* \left(\acherncl_{1}(\overline{\OO_{E^{\rho}}(1)})^{d+1}
                \cdot\widehat{X} \right) 
  + k \delta (d+1) \acherncl_{1}(\overline{E})
\in \aBChow^1(Y; Y_0(\CC)).
\]
Thus, $h_B$ is $C^{\infty}$ over $Y_0(\CC)$.

By Lemma~\ref{lemma:stability:generators:over:Z},
there are a positive integer $e$ and 
$\bfSL_r(\QQ)$-invariant homogeneous polynomials 
$P_1,\cdots,P_l \in \Sym^{e}(\Sym^{\delta}(\ZZ^R)^{d+1}{}^{\lor})$
depending only on $\rho$, $d$ and $\delta$ such that
if $X_y$ is Chow semistable for a closed point $y$ of $Y_{\QQ}$,
then $P_i(\Phi_{X_{y}}) \ne 0$
for some $P_i$.
For later use, 
we put 
$c_2 (\rho,d,\delta) 
= \max \{ \VERT P_1 \VERT, \cdots, 
\VERT P_l \VERT \}$, 
which is a constant depending only on $\rho$, $d$ and $\delta$. 

Recall that the Chow form $\Phi_X$ is an element of 
$H^0(Y,L \otimes (\Sym^{\delta}(E^{\rho}))^{{}\otimes {d+1}})$ 
and 
by Lemma~\ref{lemma:P:to:sheafhom} 
$P_i$ induces a polynomial map of sheaves 
\[
\varphi_{P_i} : L \otimes \Sym^{\delta} (E^{\rho})^{\otimes d+1} 
\to L^{\otimes er} \otimes (\det E)^{ek(d+1)\delta}. 
\]
Hence we have 
\[
\varphi_{P_i} (\Phi_{X}) 
\in H^0 \left(Y,L^{\otimes e r} 
                \otimes (\det E)^{e k(d+1)\delta} \right) 
= H^0 (Y, B^{\otimes e})
\]
by \eqref{eqn:metric:L}.
Here we set $s_i = \varphi_{P_i} (\Phi_{X})$.
Then, the property (ii) is obvious
by the construction of $\varphi_{P_i}$ 
and (i) of Proposition~\ref{proposition:pullback}. 

Now we will evaluate 
$\Vert s_i \Vert$. 
Let $y$ be a closed point of $(Y_0)_{\QQ}$ with
$y$ extensible in $Y$.
Combining 
Lemma~\ref{lemma:norm:of:morphism:P} and 
Lemma~\ref{lemma:evaluation:of:norm:of:Phi}, 
we have 
\begin{align*}
\GM\left( \Vert s_i \Vert ;\ O_{\Gal(\overline{\QQ}/{\QQ})}(y)
\right) & \leq
\GM\left( \VERT P_i \VERT^r  \, \Vert \Phi_{X} \Vert^{e r};\ 
 O_{\Gal(\overline{\QQ}/{\QQ})}(y) \right) \\
& \leq c_2(\rho,d,\delta) ^{r} c_1 (R,d,\delta) ^{e r}.
\end{align*} 
Now we put
\[
C(\rho,d,\delta) 
= c_1(R,d,\delta)^{2r} c_2(\rho,d,\delta)^{2e r},
\]
which is a positive constant depending only on $\rho$, $d$ and $\delta$.
Thus, we get (iii).
\QED

\begin{Remark}
\label{rem:geom:analog:Cornalba-Harris-Bost}
Here let us consider the geometric analogue of 
Theorem~\ref{thm:semistability:imply:average:semi-ampleness}.
Let $Y$ be an algebraic variety over an algebraically closed
field $k$, $E$ a vector bundle of rank $r$, 
$\rho : \bfGL_r \to \bfGL_R$ a group scheme morphism of degree $l$ 
commuting with the transposed morphism. 
Let $X$ be an effective cycle in $\PP (E^{\rho})$ 
such that $X$ is flat over $Y$ 
with the relative dimension $d$ 
and degree $\delta$ on the generic fiber.
Here we set
\[
b_{X/Y}(E, \rho) =
r \pi_* \left(c_{1}(\OO_{E^{\rho}}(1))^{d+1} \cdot X \right) 
  + l \delta (d+1) c_{1}(E),
\]
which is a divisor on $Y$.
In the same way as in the proof of 
Theorem~\ref{thm:semistability:imply:average:semi-ampleness}
(actually, this case is much easier than
the arithmetic case),
we can show the following.

\begin{quote}
There is a positive integer $e$ depending only on $\rho$, $d$, 
and $\delta$ such that,
if $X_y$ is Chow semi-stable for some $y \in Y$, 
then
\[
H^0(Y, \OO_Y(e b_{X/Y}(E, \rho))) \otimes \OO_Y \to
\OO_Y(e b_{X/Y}(E, \rho))
\]
is surjective at $y$.
\end{quote}
This gives a refinement of \cite[Theorem~3.2]{Bo}. 
\end{Remark}

\subsection{Relationship of two theorems}
\label{section:Bogomolov:to:Bost}
\setcounter{Theorem}{0}
In this subsection we will see some relationship between 
the relative Bogomolov's inequality 
(Theorem~\ref{thm:relative:Bogomolov:inequality:arithmetic:case}) 
and 
the relative Cornalba-Harris-Bost's inequality 
(Theorem~\ref{thm:semistability:imply:average:semi-ampleness}). 
For this purpose, we will first show a more intrinsic version of 
Theorem~\ref{thm:semistability:imply:average:semi-ampleness}. 

\begin{Proposition}
\label{prop:intrinsic:relbost}
Let $f : X \to Y$ be a flat morphism of 
regular projective arithmetic varieties with $\dim f = d$. 
Let $L$ be a relatively very ample line bundle 
such that $E = f_* (L)$ is a vector bundle of rank $r$ on $Y$. 
Let $\eta$ be the generic point of $X$ 
and $\delta = \deg (L_{\eta}^d)$. 
Moreover, let 
$i : X \to \PP(E^{\lor})$ be the embedding over $Y$. 
Assume that $E$ is equipped with an Hermitian structure $h$   
so that $L$ is also endowed with the Hermitian structure 
by $i^* \OO_{E^{\lor}}(1) \simeq L$.
Let $Y_0$ be the maximal open set of $Y$ such that
$f$ is smooth over $Y_0$.
Then, there is a positive integer
$e(r,d,\delta)$ and a positive constant $C(r,d,\delta)$
depending only on $r,d,\delta$ with the following properties.
If we set
\[
\beta = e(r,d,\delta) \left(
r f_* (\acherncl_1 (\overline{L})^{d+1}) 
  - \delta (d+1) \acherncl_{1}(\overline{E})\right)
+ a (\log C(r,d,\delta)),
\]
then, for any closed point $y \in (Y_0)_{\QQ}$ with
$X_y$ Chow semistable,
there is a representative $(D, g)$ of $\beta$
such that $D$ is effective, $y \not\in D$, and
\[
\sum_{w \in O_{\Gal(\overline{\QQ}/\QQ)}(z)} g(w) \geq 0
\]
for all $z \in (Y_0)_{\QQ}$.
\end{Proposition}

\Proof
We identify $X$ with its image by $i$.
Let $\pi : \PP(E) \to Y$ be the projection.
Then, by Proposition~\ref{prop:when:Bost:divisor:smooth}, 
we get 
\[
\pi_* \left(\acherncl_{1}(\overline{\OO_{E^{\lor}}(1)})^{d+1}
                \cdot\widehat{X} \right)
=
f_* (\acherncl_1 (\overline{L})^{d+1})
\]
Thus,
applying Theorem~\ref{thm:semistability:imply:average:semi-ampleness} 
for $(Y,E^{\lor},\operatorname{id},X)$,  
we get our assertion.
\QED

The following proposition will be derived from 
Theorem~\ref{thm:relative:Bogomolov:inequality:arithmetic:case}. 

\begin{Proposition}
\label{prop:Bogomolov:to:Bost}
Let $f : X \to Y$ be a projective morphism of regular arithmetic varieties
such that every fiber of $f_{\CC} : X(\CC) \to Y(\CC)$ is a reduced
and connected curve with only ordinary double singularities.
We assume that the genus $g$ of the generic fiber of $f$
is greater than or equal to $1$. 
Let $L$ be a line bundle on $X$ 
such that \rom{(1)}the degree $\delta$ of
$L$ on the generic fiber is greater than or equal to $2g+1$,
\rom{(2)} $E = f_* (L)$ is a vector bundle of rank $r$ on $Y$
\rom{(}actually $r = \delta + 1 - g$\rom{)},
and that \rom{(3)} $f^*(E) \to L$ is surjective.
Assume that $E$ is equipped with an Hermitian structure $h$   
so that $L$ is also endowed with the quotient metric  
by $f^* (E) \to L$.
Let $Y_0$ be the maximal open set of $Y$ such that
$f$ is smooth over $Y_0$.
Then, for any closed points $y$ of $(Y_0)_{\QQ}$, 
\[
r f_* (\acherncl_1 (\overline{L})^{2}) 
- 2 \delta \acherncl_1 (\overline{E})
\]
is weakly positive at $y$ with respect to
any finite subsets of $Y_0(\CC)$.
\end{Proposition}

Note that if the base space is $\Spec(O_K)$, then 
the second author showed in \cite[Theorem~1.1]{MorFh} 
the above inequality (under weaker assumptions) 
using \cite[Corollary~8.9]{MoBG}. 
Since we can prove Proposition~\ref{prop:Bogomolov:to:Bost} 
in the same way as \cite[Theorem~1.1]{MorFh}, 
we will only sketch the proof.

\Proof
Let $S = \Ker  (f^* (E) \to L)$ and $h_S$ the submetric of $S$ 
induced by $h$.
Then, by \cite{EL},
$S_{z}$ is stable for all $z \in Y_0(\CC)$. 
Applying Theorem~\ref{thm:relative:Bogomolov:inequality:arithmetic:case} 
for $\overline{S} = (S,h_S)$, we obtain that if $y$ is a closed point
of $(Y_0)_{\QQ}$, then 
\[
f_* ( 2 (r-1) \acherncl_2(\overline{S}) - (r-2) \acherncl_1(\overline{S})^2)
\]
is weakly positive at $y$ with respect to
any finite subsets of $Y_0(\CC)$.  
If we set $\rho = \acherncl_2(f^*\overline{E}) 
- \acherncl_2(\overline{S}\oplus \overline{L})$, 
then there is $g \in L^1_{loc}(Y(\CC))$ such that
$f_*(\rho) = a(g)$, $g$ is $C^{\infty}$ over $Y_0(\CC)$, and
$g > 0$ on $Y_0(\CC)$. 
Now by a straightforward calculation, we have
\begin{multline*}
f_* ( 2 (r-1) \acherncl_2(\overline{S}) - (r-2) \acherncl_1(\overline{S})^2) 
+ 2 (r-1) f_*(\rho) \\
= f_* \left( 2 (r-1) \acherncl_2(f^* \overline{E}) 
             - (r-2) \acherncl_1(f^* \overline{E})^2 \right) 
+ f_* \left( r \acherncl_1(\overline{L})^2 
             - 2 \acherncl_1(f^* \overline{E}) 
              \cdot \acherncl_1(\overline{L}) \right) \\
= r f_* (\acherncl_1(\overline{L})^2) 
- 2 \delta \acherncl_1(\overline{E}). 
\end{multline*}
\QED

Let us compare 
Proposition~\ref{prop:intrinsic:relbost} 
with Proposition~\ref{prop:Bogomolov:to:Bost}. 
Both of them give some arithmetic positivity 
of the same divisor 
(although $d=1$ in Proposition~\ref{prop:Bogomolov:to:Bost}), 
under the assumption of some semi-stability 
(of Chow or of vector bundles). 
The former has advantage 
since it treats varieties of arbitrary relative dimension. 
On the other hand, 
the latter has advantage 
since it shows that 
the anonymous constant in the former 
is zero (see also \cite{Zh}). 
Moreover, in the complex case, 
the counterpart of the relative Bogomolov's inequality  
of Theorem~\ref{thm:relative:Bogomolov:inequality:arithmetic:case} 
has a wonderful application to the moduli of stable curves 
(\cite{MoRB}). 

\renewcommand{\thesection}{Appendix \Alph{section}}
\renewcommand{\theTheorem}{\Alph{section}.\arabic{Theorem}}
\renewcommand{\theClaim}{\Alph{section}.\arabic{Theorem}.\arabic{Claim}}
\renewcommand{\theequation}{\Alph{section}.\arabic{Theorem}.\arabic{Claim}}
\setcounter{section}{0}

\section{Commutativity of push-forward and pull-back}
\label{sec::comm:push:pull}

Let $f : X \to Y$ be a smooth proper morphism of regular noetherian schemes, and
$u : Y' \to Y$ a morphism of regular noetherian schemes.
Let $X' = X \times_{Y} Y'$ and  
\[
\begin{CD}
X @<{u'}<< X' \\
@V{f}VV @VV{f'}V \\
Y @<{u}<< Y'
\end{CD}
\]
the induced diagram.
Let $Z$ be a cycle of codimension $p$ and $|Z|$ the support of $Z$.
We assume that $\codim_{X'}({u'}^{-1}(|Z|)) \geq p$.
Then, it is easy to see that $\codim_{Y'}(u^{-1}(|f_*(Z)|)) \geq p - d$,
where $d = \dim X - \dim Y$.
Thus, we can define $f'_*({u'}^*(Z))$ and $u^*(f_*(Z))$ as
elements of $Z^{p-d}(Y')$.
It is well known, we believe, that
$f'_*({u'}^*(Z)) = u^*(f_*(Z))$ in $Z^{p-d}(Y')$.
We could not however find any suitable references for the above fact, so that
in this section, we would like to give the proof of it.

\bigskip
Let $X$ be a regular noetherian scheme, and $T$ a closed subscheme of $X$.
We denote by $K'_T(X)$ the Grothendieck group generated by
coherent sheaves $F$ with $\Supp(F) \subseteq T_{red}$ modulo
the following relation:
$[F] = [F'] + [F'']$ if there is an exact sequence
$0 \to F' \to F \to F'' \to 0$.

Let $p$ be a non-negative integer, and
$X^{(p)}$ the set of all points $x$ of $X$ with
$\codim_X \overline{\{ x \}} = p$.
We define $Z^p_T(X)$ to be
\[
  Z^p_T(X) = \bigoplus_{ x \in X^{(p)} \cap T} \ZZ \cdot \overline{\{ x \}}.
\] 
We assume that $\codim_X T \geq p$. Then, we can define
the natural homomorphism
\[
z^p : K'_T(X) \to Z^p_T(X)
\]
to be
\[
z^p([F]) = \sum_{x \in X^{(p)} \cap T} l_{\OO_{X, x}}(F_x)
\cdot \overline{\{ x \}},
\]
where $l_{\OO_{X, x}}(F_x)$ is the length of $F_x$ as
$\OO_{X,x}$-modules.
Note that if $\codim_X T > p$, then $z^p = 0$.

\medskip
Let $f : X \to Y$ be a proper morphism of regular noetherian schemes, and
$T$ a closed subscheme of $X$. Then, we define
the homomorphism $f_* : K'_T(X) \to K'_{f(T)}(Y)$ to be
\[
f_*([F]) = \sum_{i \geq 0} [R^i f_*(F)].
\]
Here we set $d = \dim X - \dim Y$. 
Let $p$ be a non-negative integer with $\codim_X T \geq p$ and
$p \geq d$. Then, $\codim_Y f(T) \geq p - d$.
First, let us consider the following proposition.

\begin{Proposition}
\label{prop:comm:push:forward}
With notation as above, the diagram
\[
\begin{CD}
K'_T(X) @>{z^p}>> Z^p_T(X) \\
@V{f_*}VV @VV{f_*}V \\
K'_{f(T)}(Y) @>{z^{p-d}}>> Z^{p-d}_{f(T)}(Y)
\end{CD}
\]
is commutative.
\end{Proposition}

\Proof
For a coherent sheaf $F$ on $X$ with
$\Supp(F) \subseteq T_{red}$, there is a filtration
$0 = F_0 \subseteq F_1 \subseteq \cdots \subseteq F_n = F$
with $F_i/F_{i-1} \simeq \OO_X/P_i$ for some
prime ideal sheaves $P_i$ on $X$.
Then,
\[
\begin{cases}
f_*(z^p([F])) = \sum_{i=1}^n f_*(z^p([\OO_X/P_i])) \\
z^{p-d}(f_*([F])) = \sum_{i=1}^n z^{p-d}(f_*([\OO_X/P_i])).
\end{cases}
\]
Thus, we may assume that $F = \OO_X/P$ for some prime ideal
sheaf $P$ with $\Supp(\OO_X/P) \subseteq T_{red}$ and
$\codim_X(\Spec(\OO_X/P)) = p$.
We set $Z = \Spec(\OO_X/P)$. Then, $z^p([\OO_X/P]) = Z$.

First, let us consider the case where $\dim f(Z) < \dim Z$.
In this case, $f_*(Z) = 0$. On the other hand,
since $\Supp(R^i f_*(\OO_X/P))\subseteq f(Z)$,
we can see that $z^{p-d}([R^i f_*(\OO_X/P)]) = 0$ for all $i \geq 0$.
Thus, $z^{p-d}(f_*([\OO/P])) = 0$.

Next, we assume that $\dim f(Z) = \dim Z$.
Then, $Z \to f(Z)$ is generically finite. Thus,
$\Supp(R^i f_*(\OO_X/P))$ is a proper closed subset of $f(Z)$
for each $i \geq 1$. Therefore, we have
\[
z^{p-d}(f_*([\OO_X/P])) = z^{p-d}([f_*(\OO_X/P)]) = f_*(Z).
\]
\QED

Let $g : Z \to X$ be a morphism of regular noetherian schemes, and
$T$ a closed subscheme of $X$. 
Then, we define the homomorphism $g^* : K'_T(X) \to K'_{f^{-1}(T)}(Z)$
to be
\[
g^*([F]) = \sum_{i \geq 0} (-1)^i [L_if^*(F)].
\]
Let $p$ be a non-negative integer with $\codim_X T \geq p$ and
$\codim_Z (g^{-1}(T)) \geq p$. 
Here let us consider 
the following proposition.

\begin{Proposition}
\label{prop:well:def:pull:back}
Let $F$ and $G$ be coherent sheaves on $X$ with
$\Supp(F), \Supp(G) \subseteq T_{red}$.
If $z^p([F]) = z^p([G])$, then $z^{p}(g^*([F])) = 
z^{p}(g^*([G]))$.
\end{Proposition}

\Proof
Let
$0 = F_0 \subseteq F_1 \subseteq \cdots \subseteq F_n = F$
and
$0 = G_0 \subseteq G_1 \subseteq \cdots \subseteq G_m = F$
be filtrations of $F$ and $G$ respectively such that
$F_i/F_{i-1} \simeq \OO_X/P_i$ and
$G_j/G_{j-1} \simeq \OO_X/Q_j$ for some
prime ideal sheaves $P_i$ and $Q_j$ on $X$.
Then,
\[
\begin{cases}
z^{p}(g^*([F])) = \sum_{i=1}^n z^{p}(g^*([\OO_X/P_i])) \\
z^{p}(g^*([G])) = \sum_{j=1}^m z^{p}(g^*([\OO_X/Q_j]))
\end{cases}
\]
Thus, it is sufficient to show that
$z^{p}(g^*([\OO_X/P])) = 0$ for all prime ideals $P$ with 
\[
\text{
$\Supp(\OO_X/P) \subseteq T_{red}$,
$\codim_X (\Supp(\OO_X/P)) > p$ and
$\codim_Z (g^{-1}(\Supp(\OO_X/P))) = p$.
}
\]
This is a consequence of the following lemma.
\QED

\begin{Lemma}
\label{lem:vanshing:tor:alt:sum}
Let $(A, m)$ and $(B, n)$ be regular local rings,
$\phi : A \to B$ a homomorphism of local rings, and
$M$ an $A$-module of finite type.
If $\Supp(M \otimes_A B) = \{ n \}$ and
\[
\codim_{\Spec(B)} (\Supp(M \otimes_A B)) < \codim_{\Spec(A)} (\Supp(M)),
\]
then
\[
\sum_{i \geq 0} (-1)^i l_B(\operatorname{Tor}_i^A(M, B)) = 0.
\]
\end{Lemma}

\Proof
We freely use notations in \cite[Chapter~I]{SoAr}.
Let $f : \Spec(B) \to \Spec(A)$ be a morphism induced by
$\phi : A \to B$.
We set $Y = \Supp(M)$ and $q = \codim_{\Spec(A)} (\Supp(M))$.
Let $P_{\cdot} \to M$ be a free resolution of $M$.
Then, $[P_{\cdot}] \in F^q K_0^Y(\Spec(A))$.
Thus, by \cite[(iii) of Theorem~3 in I.3]{SoAr},
\[
[f^*(P_{\cdot})] = [ P_{\cdot} \otimes_A B ]
\in F^q K_0^{\{ n \}}(\Spec(B))_{\QQ}
\]
because $f^{-1}(Y) = \Supp(M \otimes_A B) = \{ n \}$.
On the other hand, since
\[
 q > \codim_{\Spec(B)} (\Supp(M \otimes_A B)) = \dim B,
\]
we have $F^q K_0^{\{ n \}}(\Spec(B))_{\QQ} = \{ 0 \}$.
Thus, $[ P_{\cdot} \otimes_A B ] = 0$ in $K_0^{\{ n \}}(\Spec(B))$
because 
\[
K_0^{\{ n \}}(\Spec(B)) \simeq \ZZ
\]
has no torsion.
This shows us our assertion.
\QED

As a corollary of Proposition~\ref{prop:well:def:pull:back},
we have the following.

\begin{Corollary}
\label{cor:comm:pull:back}
With notation as in Proposition~\rom{\ref{prop:well:def:pull:back}},
\[
\begin{CD}
K'_T(X) @>{z^p}>> Z^p_T(X) \\
@V{g^*}VV @VV{g^*}V \\
K'_{f^{-1}(T)}(Z) @>{z^{p}}>> Z^{p}_{f^{-1}(T)}(Z)
\end{CD}
\]
is commutative. Note that
$g^* : Z^p_T(X) \to  Z^{p}_{f^{-1}(T)}(Z)$ is defined by
$g^*(Z) = z^{p}(g^*([\OO_Z]))$ for each integral cycle $Z$
in $Z^p_T(X)$.
\end{Corollary}

Let $f : X \to Y$ be a flat proper morphism of regular noetherian schemes, and
$u : Y' \to Y$ a morphism of regular noetherian schemes.
Let $X' = X \times_{Y} Y'$ and  
\[
\begin{CD}
X @<{u'}<< X' \\
@V{f}VV @VV{f'}V \\
Y @<{u}<< Y'
\end{CD}
\]
the induced diagram.
We assume that $X'$ is regular.
Note that if $f$ is smooth, then $X'$ is regular.
We set $d = \dim X - \dim Y = \dim X' - \dim Y'$.
Let $T$ be a closed subscheme of $X$, and $p$ a non-negative
integer with $\codim_X T \geq p$, $\codim_{X'} ({u'}^{-1}(T)) \geq p$
and $p \geq d$.
Note that $\codim_Y f(T) \geq p-d$ and 
$\codim_{Y'} (u^{-1}(f(T))) \geq p - d$ because
$u^{-1}(f(T)) = f'({u'}^{-1}(T))$.
Then, we have the following proposition.

\begin{Proposition}
\label{prop:comm:cycle:push:pull}
The diagram
\[
\begin{CD}
Z^p_T(X) @>{{u'}^*}>> Z^{p}_{{u'}^{-1}(T)}(X') \\
@V{f_*}VV @VV{f'_*}V \\
Z^{p-d}_{f(T)}(Y) @>{u^*}>> Z^{p-d}_{u^{-1}(f(T))}(Y')
\end{CD}
\]
is commutative.
\end{Proposition}

\Proof
Since $f$ is flat, by \cite[Proposition 3.1.0 in IV]{SGA6},
for any coherent sheaves $F$ on $X$,
\[
L_{\cdot}u^* \left( R^{\cdot} f_*(F) \right) \overset{\sim}{\longrightarrow}
R^{\cdot}f'_* \left( L_{\cdot}{u'}^*(F)\right),
\]
which shows that the diagram
\[
\begin{CD}
K'_T(X) @>{{u'}^*}>> K'_{{u'}^{-1}(T)}(X') \\
@V{f_*}VV @VV{f'_*}V \\
K'_{f(T)}(Y) @>{u^*}>> K'_{u^{-1}(f(T))}(Y')
\end{CD}
\]
is commutative.
Thus, by virtue of Proposition~\ref{prop:comm:push:forward} and
Corollary~\ref{cor:comm:pull:back},
we can see our proposition.
\QED

\end{document}